\documentclass[12pt]{article}
\usepackage{amssymb}
\usepackage{amsfonts}
\usepackage{amsmath}
\usepackage{amsthm}
\usepackage{hyperref}
\usepackage{graphicx}
\usepackage{subfig}

\theoremstyle{plain}
\theoremstyle{plain}
\newtheorem{theorem}{Theorem}[section]

\newtheorem*{algorithm}{Globally Convergent Algorithm}

\newtheorem{lemma}{Lemma}[section]

\newtheorem*{nonoproblem}{Problem}

\newtheorem{remark}{Remark}[section]

\input{tcilatex}
\newcommand{\height}{3.5cm}
\newcommand{\width}{3.5cm}
\newcommand{\eight}{3.5cm}
\newcommand{\idth}{3.5cm}
\newcommand{\labela}{True $c(x)$. The $x-$axis indicates $x \in (0, 1)$}
\newcommand{\labelb}{The real (solid line) and imaginary parts of the data $g_0(k)$, defined in \eqref{2.7}.  The $x-$axis indicates $k \in (0.5, 1.5)$.}
\newcommand{\labelc}{The reconstruction of $c$ when $m^0 = 24, 25$. The $x-$axis indicates $x \in (0, 1)$.}
\newcommand{\labelaa}{{The scattering field in the time domain} }
\newcommand{\labelbb}{The scattering field in the frequency domain}
\newcommand{\labelcc}{The scattering field in the frequency domain after cutting off its small information}
\newcommand{\labeldd}{The reconstruction of $R(x)$}
\begin{document}

\date{}
\title{ A globally convergent numerical method for a 1-d inverse medium
problem with experimental data\thanks{
The work of the first two authors was supported by US Army Research
Laboratory and US Army Research Office grant W911NF-15-1-0233 and by the
Office of Naval Research grant N00014-15-1-2330.}}
\author{Michael V. Klibanov\thanks{
Department of Mathematics and Statistics, University of North Carolina at
Charlotte, Charlotte, NC 28213, USA (\texttt{mklibanv{\@@}uncc.edu}, \texttt{%
lnguye50{\@@}uncc.edu}).}, \and Loc H. Nguyen\footnotemark[2], \and Anders
Sullivan\thanks{%
US Army Research Laboratory, 2800 Powder Mill Road, Adelphy, MD 20783-1197,
USA (\texttt{anders.j.sullivan.civ{\@@}mail.mil}, \texttt{lam.h.nguyen2civ{%
\@@}mail.mil})} \and and Lam Nguyen\footnotemark[3]}
\maketitle

\begin{abstract}
In this paper, a reconstruction method for the spatially distributed
dielectric constant of a medium from the back scattering wave field in the
frequency domain is considered. Our approach is to propose a globally
convergent algorithm, which does not require any knowledge of a small
neighborhood of the solution of the inverse problem in advance. The
Quasi-Reversibility Method (QRM) is used in the algorithm. The convergence
of the QRM is proved via a Carleman estimate. The method is tested on both
computationally simulated and experimental data.
\end{abstract}



\numberwithin{equation}{section}

%
%
%
%

\section{Introduction}

\label{sec:1}

In this paper we develop a globally convergent numerical method for a 1-d
inverse medium problem in the frequency domain. The performance of this
method is tested on both computationally simulated and experimental data.
Propagation of electromagnetic waves is used in experiments. A theorem about
the global convergence of our method is proved. Another name for inverse
medium problems is Coefficient Inverse Problems (CIPs). We model the process
as a 1-d CIP due to some specifics of our data collection procedure, see
subsection 6.2. This paper is the first one in which the globally convergent
method of \cite{BK1,BK2,Karch,KSNF1,KSNF2,TBKF1,TBKF2} is extended to the
case of the frequency domain. Indeed, the original version of that method
works with the Laplace transform of the time dependent data. Both the 3-d
and the 1-d versions of the method of \cite{BK1,BK2} were verified on
experimental data, see \cite{BK1,TBKF1,TBKF2} for the 3-d case and \cite%
{Karch,KSNF1,KSNF2} for the 1-d case. The experimental data here are the
same ones as in \cite{Karch,KSNF1,KSNF2}.

Compared with the previous publications \cite{Karch,KSNF1,KSNF2}, two
additional difficulties occurring here are: (1) we now work with complex
valued functions instead of real valued ones and, therefore (2) it is not
immediately clear how to deal with the imaginary part of the logarithm of
the complex valued solution of the forward problem. On the other hand, we
use that logarithm in our numerical procedure. Besides, the previous
numerical scheme of \cite{Karch,KSNF1,KSNF2} is significantly modified here.

We call a numerical method for a CIP \emph{globally convergent} if there is
a rigorous guarantee that this method reaches at least one point in a
sufficiently small neighborhood of the correct solution without any advanced
knowledge about this neighborhood. It is well known that the topic of the
global convergence for CIPs is a highly challenging one. Thus, similarly
with the above cited publications, we prove the global convergence of our
method within the framework of a certain approximation. This is why the term
\textquotedblleft approximate global convergence" was used in above cited
references (\textquotedblleft global convergence" in short). That
approximation is used only on the first iterative step of our algorithm and
it is not used in follow up iterations.\ Besides, this approximation is a
quite natural one, since it amounts to taking into account only the first
term of a certain asymptotic behavior and truncating the rest. The global
convergence property is verified numerically here on both computationally
simulated and experimental data.

CIPs are both highly nonlinear and ill-posed. These two factors cause the
non-convexity of conventional Tikhonov functionals for these problems. It is
well known that typically those functionals have multiple local minima and
ravines, see, e.g. numerical examples of multiple local minima in \cite%
{Scales} and on Figure 5.3 of \cite{KT}. Figure 5.4 of \cite{KT}
demonstrates an example of a ravine. Hence, there is no guarantee of the
convergence of an optimization method to the correct solution, unless the
starting point of iterations is located in a sufficiently small neighborhood
of that solution. However, it is often unclear in practical scenarios how to
reach that neighborhood.

In our inverse algorithm we solve a sequence of linear ordinary differential
equations of the second order on the interval $x\in \left( 0,1\right) $. The
peculiarity here is that we have overdetermined boundary conditions for
these equations: we have the solution and its first derivative at $x=0$ and
we have the first derivative of the solution at $x=1$. Our attempts to use
only boundary conditions at $x=0$ did not lead to good reconstruction
results. The same observation is in place in Remark 5.1 on page 13 of \cite%
{KSNF1}. Thus, similarly with \cite{Karch,KSNF1,KSNF2}, we use here the
Quasi-Reversibility Method (QRM). The QRM is well suitable to solve PDEs
with the overdetermined boundary data. To analyze the convergence rate of
the QRM, we use a Carleman estimate.

The QRM was first introduced by Lattes and Lions \cite{LL}. However, they
have not established convergence rates. These rates were first established
in \cite{KS,KM}, where it was shown that Carleman estimates are the key tool
for that goal. A survey of applications of Carleman estimates to QRM can be
found in \cite{KAPNUM}. Chapter 6 of \cite{BK1} describes the use of the QRM
for numerical solutions of CIPs. We also mention an active work with the QRM
of Bourgeois and Dard\`{e}, see, e.g. \cite{BD1,BD2,BD3} for some samples of
their publications.

Our experimental data are given in the time domain. To obtain the data in
the frequency domain, we apply the Fourier transform. Previously the same
data were treated in works \cite{Karch,KSNF1,KSNF2} of this group. In that
case the Laplace transform was applied. Next, the 1-d version of the
globally convergent method of \cite{BK1,BK2} was used. Indeed, the original
version of this method works with the Laplace transform of the time
dependent data. A similar 1-d inverse problem in the frequency domain was
solved numerically in \cite{Lesnic} by a different method.

The experimental data of this publication were collected by the Forward
Looking Radar which was built in the US Army Research Laboratory \cite{Ng}.
The data were collected in the field (as opposed to a laboratory). Thus,
clutter was present at the measurement site. This certainly adds some
additional difficulties to the imaging problem. The goal of this radar is to
detect and identify shallow explosives, such as, e.g. improvised explosive
devices and antipersonnel land mines. Those explosives can be located either
on the ground surface or a few centimeters below this surface. This radar
provides only a single time dependent curve for a target. Therefore, to
solve a CIP, we have no choice but to model the 3-d process by a 1-d
wave-like PDE. We model targets as inclusions whose dielectric constants are
different from the background.

In terms of working with experimental data, the goal in this paper is not to
image locations of targets, since this is impossible via solving a CIP with
our data, see subsection 6.2 for details. In fact, our goal is to image
maximal values of dielectric constants of targets. Our targets are 3-d
objects, while we use a 1-d model: since we measure only one time resolved
curve for a target. Nevertheless, we show below that our calculated
dielectric constants are well within tabulated limits \cite{Tabl}.

Of course, an estimate of the dielectric constant is insufficient for the
discrimination of explosives from the clutter. On the other hand, the radar
community mostly relies on the intensity of the radar image. Therefore, we
hope that the additional information about values of dielectric constants of
targets might lead in the future to designs of algorithms which would better
discriminate between explosives and clutter.

In the 3-d case the globally convergent method of \cite{BK1,BK2,TBKF1,TBKF2}
works with the data generated either by a single location of the source or
by a single direction of the incident plane wave. The second globally
convergent method for this type of measurements was developed and tested
numerically in \cite{KTSIAP}. We also refer to another globally convergent
methods for a CIPs, which is based on the multidimensional version of the
Gelfand-Levitan method. This version was first developed by Kabanikhin \cite%
{Kab1} and later by Kabanikhin and Shishlenin \cite{Kab2,Kab3}. Unlike \cite%
{BK1,BK2,TBKF1,TBKF2}, the technique of \cite{Kab1,Kab2,Kab3} works with
multiple directions of incident plane waves.

In \cite{Karch} the performance of the method of \cite{KSNF1,KSNF2} was
compared numerically with the performance of the Krein equation \cite{Krein}
for the same experimental data as ones used here. The Krein equation \cite%
{Krein} is close to the Gelfand-Levitan equation \cite{Lev}. It was shown in 
\cite{Krein} that the performance of the Krein equation is inferior to the
performance of the method of \cite{KSNF1,KSNF2} for these experimental data.
This is because the solution of the Krein equation is much more sensitive to
the choice of the calibration factor than the solution obtained by the
technique of \cite{KSNF1,KSNF2}. On the other hand, it is necessary to apply
that factor to those experimental data to make the range of values of the
resulting data comparable with the range of values of computationally
simulated data.

In section 2 we pose forward and inverse problems and analyze some of their
features. In section 3 we study a 1-d version of the QRM. In section 4 we
present our globally convergent method. In section 5 we prove the global
convergence of our method. In section 6 we present our numerical results for
both computationally simulated and experimental data.

\section{ Some Properties of Forward and Inverse Problems}

\label{sec:2}

It was shown numerically in \cite{BMM} that the component of the electric
field, which is incident upon a medium, dominates two other components.\ It
was also shown in \cite{BMM} that the propagation of that component is well
governed by a wave-like PDE. Furthermore, this finding was confirmed by
imaging from experimental data, see Chapter 5 of \cite{BK1} and \cite%
{TBKF1,TBKF2}. Thus, just as in \cite{Karch,KSNF1,KSNF2}, we use a 1-d
wave-like PDE to model the collection of our experimental data of
electromagnetic waves propagation.

\subsection{Formulations of problems}

\label{sec:2.1}

Let $c_{0}<c_{1}$ be two positive numbers. Let the function $c:\mathbb{R}%
\rightarrow \mathbb{R}$ satisfy the following conditions:%
\begin{eqnarray}
c &\in &C^{2}\left( \mathbb{R}\right) ,c\left( x\right) \in \left[
c_{0},c_{1}\right] ,\forall x\in \mathbb{R},  \label{2.1} \\
c\left( x\right) &=&1+\beta \left( x\right) ,\beta \left( x\right)
=0,\forall x\notin \left( 0,1\right) .  \label{2.3}
\end{eqnarray}%
Fix the source position $x_{0}<0$. Consider the generalized Helmholtz
equation in the 1-d case,%
\begin{eqnarray}
u^{\prime \prime }+k^{2}c\left( x\right) u &=&-\delta \left( x-x_{0}\right)
,x\in \mathbb{R},  \label{2.4} \\
\lim_{x\rightarrow \infty }\left( u^{\prime }+iku\right)
&=&0,\lim_{x\rightarrow -\infty }\left( u^{\prime }-iku\right) =0.
\label{2.6}
\end{eqnarray}%
Let $u_{0}\left( x,x_{0},k\right) $ be the solution of the problem (\ref{2.4}%
), (\ref{2.6}) for the case $c\left( x\right) \equiv 1.$ Then%
\begin{equation}
u_{0}\left( x,x_{0},k\right) =\frac{\exp \left( -ik\left\vert
x-x_{0}\right\vert \right) }{2ik}.  \label{2.60}
\end{equation}%
The problem \eqref{2.4}, \eqref{2.6} is the forward problem. Our interest is
in the following inverse problem problem:

\begin{nonoproblem}[Coefficient Inverse Problem (CIP)]
Fix the source position $x_{0}<0.$ Let $[\underline{k},\overline{k}]\subset
\left( 0,\infty \right) $ be an interval of frequencies $k$. Reconstruct the
function $\beta \left( x\right) ,$ assuming that the following function $%
g_{0}\left( k\right) $ is known 
\begin{equation}
g_{0}\left( k\right) =\frac{u(0,x_{0},k)}{u_{0}(0,x_{0},k)},k\in \lbrack 
\underline{k},\overline{k}].  \label{2.7}
\end{equation}
\end{nonoproblem}

\subsection{Some properties of the solution of the forward problem}

\label{sec:2.2}

In this subsection we establish existence and uniqueness of the forward
problem. Even though these results are likely known, we present them here
for reader's convenience. In addition, the techniques using in their proofs
help us to verify that the function $u(x,x_{0},k)$ never vanishes for $%
x>x_{0}$, which plays an important role in our algorithm of solving the to
CIP. Also, this technique justifies the numerical method for solving the
forward problem (\ref{2.4}), (\ref{2.6}). Computationally simulated data are
obtained by numerically solving the problem (\ref{2.4}), (\ref{2.6}).

\begin{theorem}
\label{thm unique forward} Assume that conditions (\ref{2.1})-(\ref{2.3})
hold. Then for each $k>0$ and for each $x_{0}<0$ there exists a single
solution $u\left( x,x_{0},k\right) $ of the problem (\ref{2.4}), (\ref{2.6}%
). Moreover, the function $\widetilde{u}\left( x,x_{0},k\right) =u\left(
x,x_{0},k\right) -u_{0}\left( x,x_{0},k\right) $, called the scattering
field, is in $C^{3}\left( \mathbb{R}\right) $.
\end{theorem}

\textbf{Proof}. We prove uniqueness first. Suppose that ${u}_{1}$ and ${u}%
_{2}$ are two solutions of the problem (\ref{2.4}), (\ref{2.6}). Denote $%
U=u_{1}-u_{2}.$ Then the function $U$ satisfies 
\begin{eqnarray}
U^{\prime \prime }+k^{2}\left( 1+\beta \left( x\right) \right) U &=&0,x\in 
\mathbb{R},  \label{2.8} \\
\lim_{x\rightarrow \infty }\left( U^{\prime }+ikU\right)
&=&0,\lim_{x\rightarrow -\infty }\left( U^{\prime }-ikU\right) =0.
\label{2.9}
\end{eqnarray}%
Since the function $\beta \left( x\right) =0$ outside of the interval $%
\left( 0,1\right) ,$ then (\ref{2.8}) implies that $U^{\prime \prime
}+k^{2}U=0$ for $x\notin \left( 0,1\right) .$ This, together with (\ref{2.9}%
), yields 
\begin{equation}
U\left( x,x_{0},k\right) =\left\{ 
\begin{array}{c}
B_{1}\left( x_{0},k\right) e^{-ikx},x>1, \\ 
B_{2}\left( x_{0},k\right) e^{ikx},x<0,%
\end{array}%
\right.  \label{2.11}
\end{equation}%
where $B_{1}$ and $B_{2}$ are some complex numbers depending on $x_{0},k.$

Let $R>1$ be an arbitrary number. Multiply both sides of (\ref{2.8}) by the
function $\overline{U}$ and integrate over the interval $\left( -R,R\right) $
using integration by parts. We obtain%
\begin{equation}
\left( U^{\prime }\overline{U}\right) \left( R\right) -\left( U^{\prime }%
\overline{U}\right) \left( -R\right) -\dint\limits_{-R}^{R}\left\vert
U^{\prime }\right\vert ^{2}dx+k^{2}\dint\limits_{-R}^{R}\left( 1+\beta
\left( x\right) \right) \left\vert U\right\vert ^{2}dx=0.  \label{2.12}
\end{equation}%
By (\ref{2.11}) $\left( U^{\prime }\overline{U}\right) \left( R\right)
=-ik\left\vert B_{1}\right\vert ^{2}$ and $-\left( U^{\prime }\overline{U}%
\right) \left( -R\right) =-ik\left\vert B_{2}\right\vert ^{2}.$ Hence, the
imaginary part of (\ref{2.12}), $-k\left( \left\vert
B_{1}(x_{0},k)\right\vert ^{2}+\left\vert B_{2}(x_{0},k)\right\vert
^{2}\right) =0,$ vanishes, so do $B_{1}(x_{0},k)$ and $B_{2}(x_{0},k)$.
Using (\ref{2.11}), we obtain $U\left( x,x_{0},k\right) =0$ for $x\notin
\left( 0,1\right) .$ Due to the classical unique continuation principle, $%
U\left( x,x_{0},k\right) =0,\forall x\in \mathbb{R}.$ Thus, $u_{1}=u_{2}$.

We now prove existence. Consider the 1-d analog of the Lippman-Schwinger
equation with respect to a function $P$, 
\begin{equation}
P\left( x,x_{0},k\right) =\frac{\exp \left( -ik\left\vert x-x_{0}\right\vert
\right) }{2ik}+\frac{k}{2i}\dint\limits_{0}^{1}\exp \left( -ik\left\vert
x-\xi \right\vert \right) \beta \left( \xi \right) P\left( \xi
,x_{0},k\right) d\xi ,x\in \mathbb{R}.  \label{2.13}
\end{equation}%
Fix $x_{0}<0$ and $k>0$. Consider equation (\ref{2.13}) only for $x\in
\left( 0,1\right) .$ Assume that there exist two functions $P_{1},P_{2}$
satisfying (\ref{2.13}) for $x\in \left( 0,1\right) .$ Consider their
extensions on the set $\mathbb{R}\diagdown \left( 0,1\right) $ via the right
hand side of (\ref{2.13}). Then so defined functions $P_{1},P_{2}$ satisfy (%
\ref{2.13}) for all $x\in $ $\mathbb{R}.$ Hence, both of them are solutions
of the problem (\ref{2.4}), (\ref{2.6}). Hence, the above established
uniqueness result implies that $P_{1}\equiv P_{2}.$

Consider again the integral equation (\ref{2.13}) for $x\in \left(
0,1\right) $ and apply the Fredholm alternative. This alternative, combined
with the discussion in the previous paragraph, implies that there exists
unique solution $P\in C\left[ 0,1\right] $ of equation (\ref{2.13}).
Extending this function for $x\in \mathbb{R}\diagdown \left( 0,1\right) $
via the right hand side of (\ref{2.13}), we obtain that there exists unique
solution $P\in C\left( \mathbb{R}\right) $ of equation (\ref{2.13}).
Furthermore, the function $P-u_{0}\in C^{3}\left( \mathbb{R}\right) $ and
also this function $P$ is the required solution of the problem (\ref{2.4}), (%
\ref{2.6}). $\square $

Below we consider only such a solution $u$ of the problem (\ref{2.4}), (\ref%
{2.6}) that $u-u_{0}\in C^{3}\left( \mathbb{R}\right) $: as in Theorem 2.1.

\subsection{The asymptotic behavior of the function $u(x,x_{0},k)$ as $%
k\rightarrow \infty $}

\label{sec:2.3}

Consider the function $\phi \left( x\right) $ defined as 
\begin{equation}
\phi \left( x\right) =-\frac{c^{\prime \prime }\left( x\right) }{c^{2}\left(
x\right) }+\frac{7}{16}\frac{\left( c^{\prime }\left( x\right) \right) ^{2}}{%
c^{3}\left( x\right) }.  \label{200}
\end{equation}%
Note that by (\ref{2.3}) $\phi \left( x\right) =0$ for $x\notin \left[ 0,1%
\right] .$

\begin{theorem}
\label{thm assymptotic u} Assume that the function $\phi \left( x\right)
\leq 0,\forall x\in \left[ 0,1\right] .$ Then for every pair $(x,x_{0})\in 
\mathbb{R}\times (-\infty ,0)$, the asymptotic behavior of the function $%
u\left( x,x_{0},k\right) $ is 
\begin{equation}
u\left( x,x_{0},k\right) =\frac{1}{2ikc^{1/4}\left( x\right) }\exp \left[
-ik\left\vert \dint\limits_{x_{0}}^{x}\sqrt{c\left( \xi \right) }d\xi
\right\vert \right] \left( 1+O\left( \frac{1}{k}\right) \right)
,k\rightarrow \infty .  \label{201}
\end{equation}%
Furthermore, for any finite interval $\left( a,b\right) \subset \mathbb{R}$
there exists a number $\gamma =\gamma \left( a,b,\phi \right) >0$ such that
for all $x\in \left( a,b\right) $ the function $u\left( x,x_{0},k\right) $
can be analytically extended with respect $k$ from $\left\{ k:k>0\right\} $
in the half plane $\mathbb{C}_{\gamma }=\left\{ z\in \mathbb{C}:\func{Im}%
z<\gamma \right\} $.
\end{theorem}

\textbf{Proof}. Consider the following Cauchy problem%
\begin{eqnarray}
c\left( x\right) \widehat{u}_{tt} &=&\widehat{u}_{xx},x\in \mathbb{R},t>0,
\label{2.15} \\
\widehat{u}\left( x,0\right) &=&0,\widehat{u}_{t}\left( x,0\right) =\delta
\left( x-x_{0}\right) .  \label{2.16}
\end{eqnarray}%
This proof is based on the connection between the solution $u$ of the
problem (\ref{2.4}), (\ref{2.6}) and the solution $\widehat{u}$ of the
problem (\ref{2.15}), (\ref{2.16}) via the Fourier transform with respect to 
$t$. Here and below the Fourier transform is understood in terms of
distributions, see, e.g. the book \cite{Vlad} for the Fourier transform of
distributions.

We now obtain a hyperbolic equation with potential from equation (\ref{2.15}%
). To do this, we use a well known change of variables $x\Leftrightarrow y$ 
\cite{Rom}%
\begin{equation}
y=\dint\limits_{x_{0}}^{x}\sqrt{c\left( \xi \right) }d\xi .  \label{2.17}
\end{equation}%
Denote $b\left( y\right) =c\left( x\left( y\right) \right) .$ By (\ref{2.17}%
) $db/dy=c^{\prime }\mid _{x=x\left( y\right) }c^{-1/2}\left( x\left(
y\right) \right) .$ Hence, (\ref{2.15}) and (\ref{2.16}) become%
\begin{eqnarray}
\widehat{u}_{tt} &=&\widehat{u}_{yy}+\frac{b^{\prime }\left( y\right) }{%
2b\left( y\right) }\widehat{u}_{y},\text{ }y\in \mathbb{R},t>0,  \label{2.18}
\\
\widehat{u}\left( y,0\right) &=&0,\widehat{u}_{t}\left( y,0\right) =\delta
\left( y\right) .  \label{2.19}
\end{eqnarray}%
Consider now a new function $\widehat{v}\left( y,t\right) =\widehat{u}\left(
y,t\right) /S\left( y\right) ,$ where the function $S\left( y\right)
=b^{-1/4}\left( y\right) $ is chosen in such a way that the coefficient at $%
\widehat{v}_{y}$ becomes zero after the substitution $\widehat{u}=S\widehat{v%
}$ in equation (\ref{2.18}). Then (\ref{2.18}) and (\ref{2.19}) become%
\begin{eqnarray}
\widehat{v}_{tt} &=&\widehat{v}_{yy}+p\left( y\right) \widehat{v},
\label{2.22} \\
\widehat{v}\left( y,0\right) &=&0,\widehat{v}_{t}\left( y,0\right) =\delta
\left( y\right) ,  \label{2.23} \\
p\left( y\right) &=&-\frac{b^{\prime \prime }\left( y\right) }{4b\left(
y\right) }+\frac{5}{16}\frac{\left( b^{\prime }\left( y\right) \right) ^{2}}{%
b^{2}\left( y\right) }.  \label{2.24}
\end{eqnarray}
It follows from (\ref{2.1}), (\ref{2.3}), (\ref{2.17}) and (\ref{2.24}) that
the function $p\in C\left( \mathbb{R}\right) $ and has a finite support, 
\begin{equation}
p\left( y\right) =0\text{ for }y<-x_{0}\text{ and for }y>\dint%
\limits_{x_{0}}^{1}\sqrt{c\left( \xi \right) }d\xi .  \label{2.25}
\end{equation}%
Recall that $\phi \left( x\right) $ is the function defined in (\ref{200}).
Expression (\ref{2.24}) in the $x-$coordinate becomes 
\begin{equation}
p\left( y\left( x\right) \right) =\phi \left( x\right) .  \label{2.250}
\end{equation}

Let $H\left( z\right) ,z\in \mathbb{R}$ be the Heaviside function. It is
well known that the solution $\widehat{v}\left( y,t\right) $ of the problem (%
\ref{2.22}), (\ref{2.23}) has the following form, see, e.g. Chapter 2 in 
\cite{Rom}%
\begin{eqnarray}
\widehat{v}\left( y,t,x_{0}\right) &=&\frac{1}{2}H\left( t-\left\vert
y\right\vert \right) +\widetilde{v}\left( y,t\right) H\left( t-\left\vert
y\right\vert \right) ,  \label{2.26} \\
\widetilde{v}\left( y,t\right) &\in &C^{2}\left( t\geq \left\vert
y\right\vert \right) ,\lim_{t\rightarrow \left\vert y\right\vert ^{+}}%
\widetilde{v}\left( y,t\right) =0.  \label{2.27}
\end{eqnarray}

The backwards substitution $y\rightarrow x$ transforms the function $p\left(
y\right) $ in the function $\phi \left( x\right) .$ Hence, (\ref{2.250})
implies that the function $p\left( y\right) $ in (\ref{2.24}) is non
positive, 
\begin{equation}
p\left( y\right) \leq 0,\forall y\in \mathbb{R}.  \label{2.28}
\end{equation}%
Consider now the operator $\mathcal{F}$ of the Fourier transform,%
\begin{equation}
\mathcal{F}\left( f\right) \left( k\right) =\dint\limits_{0}^{\infty
}f\left( t\right) e^{-ikt}dt,\forall f\in L_{1}\left( 0,\infty \right) .
\label{2.2727}
\end{equation}%
In the sense of distributions we have%
\begin{equation}
\mathcal{F}\left( \frac{1}{2}H\left( t-\left\vert y\right\vert \right)
\right) =\frac{\exp \left( -ik\left\vert y\right\vert \right) }{2ik},k>0,
\label{2.29}
\end{equation}%
see section 6 in \S 9 of Chapter 2 of \cite{Vlad}. Consider an arbitrary
finite interval $\left( a,b\right) \subset \mathbb{R}.$ Let $s,m\geq 0$ be
two arbitrary integers such that $s+m\leq 2$. Since by (\ref{2.25}) the
function $p\left( y\right) $ has a finite support, then (\ref{2.28}), Lemma
6 of Chapter 10 of the book \cite{V} as well as Remark 3 after that lemma
guarantee that functions $\partial _{y}^{s}\partial _{t}^{m}\widetilde{v}%
\left( y,t\right) $ decay exponentially with respect to $t,$ as long as $%
y\in \left( a,b\right) .$ Hence, one can apply the operator $\mathcal{F}$ to
the functions $\partial _{y}^{s}\partial _{t}^{m}\widetilde{v}\left(
y,t\right) $ in the regular sense. The assertion about the analytic
extension follows from \eqref{2.2727} and the exponential decay of the
functions $\partial _{y}^{s}\partial _{t}^{m}\widetilde{v}\left( y,t\right) $

Let $\widehat{V}\left( y,k,x_{0}\right) =\mathcal{F}\left( \widehat{v}%
\right) .$ The next question is whether the function $\widehat{V}$ satisfies
analogs of conditions (\ref{2.4}), (\ref{2.6}). Theorem 3.3 of \cite{V1} and
theorem 6 of Chapter 9 of \cite{V} guarantee that the function $\widehat{V}$
\ satisfies the following conditions 
\begin{eqnarray}
\widehat{V}^{\prime \prime }+k^{2}\widehat{V}+p\left( y\right) \widehat{V}
&=&-\delta \left( y\right) ,y\in \mathbb{R},  \label{2.31} \\
\lim_{x\rightarrow \infty }\left( \widehat{V}^{\prime }+ik\widehat{V}\right)
&=&0,\lim_{x\rightarrow -\infty }\left( \widehat{V}^{\prime }-ik\widehat{V}%
\right) =0.  \label{2.32}
\end{eqnarray}

\ Using (\ref{2.27}) and the integration by parts, we obtain for $%
k\rightarrow \infty $ 
\begin{equation}
\mathcal{F}\left( H\left( t-\left\vert y\right\vert \right) \widetilde{v}%
\right) =\exp \left( -ik\left\vert y\right\vert \right)
\dint\limits_{0}^{\infty }\widetilde{v}\left( y,t+\left\vert y\right\vert
\right) e^{-ikt}dt=\exp \left( -ik\left\vert y\right\vert \right) O\left( 
\frac{1}{k^{2}}\right) .  \label{2.33}
\end{equation}%
Two backwards substitution: $y\rightarrow x$ via (\ref{2.17}) and $\widehat{V%
}\rightarrow u\left( x,x_{0},k\right) =S\left( y\left( x\right) \right) 
\widehat{V}\left( y\left( x\right) ,x_{0},k\right) $ imply that conditions (%
\ref{2.31}) and (\ref{2.32}) for the function $\widehat{V}$ turn into
conditions (\ref{2.4}), (\ref{2.6}) for the function $u$. Thus, (\ref{2.29})
and (\ref{2.33}) imply (\ref{201}). $\square $

\subsection{Some properties of the solution of the inverse problem}

\label{sec:2.4}

Since the source position $x_{0}<0$ is fixed, we drop everywhere below the
dependence on $x_{0}$ in notations of functions. First, we show that, having
the function $g_{0}\left( k\right) $ in (\ref{2.7}), one can uniquely find
the function $u_{x}\left( 0,k\right) .$ Indeed, for $x<0$ conditions (\ref%
{2.4}) and (\ref{2.6}) become%
\begin{equation}
u^{\prime \prime }+k^{2}u=-\delta \left( x-x_{0}\right) ,x\in \left( -\infty
,0\right) ,  \label{2.34}
\end{equation}%
\begin{equation}
\lim_{x\rightarrow -\infty }\left( u_{x}-iku\right) =0.  \label{2.35}
\end{equation}%
Let $\widetilde{u}\left( x,k\right) =u-u_{0}.$ Then (\ref{2.34}) and (\ref%
{2.35}) imply that%
\begin{equation*}
\widetilde{u}^{\prime \prime }+k^{2}\widetilde{u}=0\text{ for }x\in \left(
-\infty ,0\right) \text{ and }\lim_{x\rightarrow -\infty }\left( \widetilde{u%
}_{x}-ik\widetilde{u}\right) =0.
\end{equation*}%
Hence, 
\begin{equation}
\widetilde{u}(x,k)=B\left( k\right) e^{ikx},\quad x<0,  \label{def B}
\end{equation}%
for a complex number $B\left( k\right) $. Note that by (\ref{2.60}) and (\ref%
{2.7}) 
\begin{equation*}
B(k)=\left( u-u_{0}\right) \left( 0,k\right) =\left( g_{0}\left( k\right)
-1\right) u_{0}\left( 0,k\right) =\frac{\left( g_{0}\left( k\right)
-1\right) e^{ikx_{0}}}{2ik}.
\end{equation*}%
%
%
%
%
%
%
%
Hence, by \eqref{def B} and the definition of the function $\widetilde{u}$,
we have 
\begin{equation}
u_{x}\left( 0,k\right) =\left( \frac{g_{0}\left( k\right) }{2}-1\right)
e^{ikx_{0}}.  \label{2.36}
\end{equation}

Let 
\begin{equation}
w\left( x,k\right) =\frac{u\left( x,k\right) }{u_{0}\left( x,k\right) }.
\label{2.37}
\end{equation}%
Thus, \eqref{2.7}, \eqref{2.36} and (\ref{2.37}) provide us with an
additional data $g_{1}(k)$ to solve our CIP, where 
\begin{equation}
g_{1}(k):=w_{x}(0,k)=2ik\left( g_{0}\left( k\right) -1\right) ,k\in \lbrack 
\underline{k},\overline{k}].  \label{2.38}
\end{equation}

\textbf{Theorem 2.3 (}uniqueness of our CIP\textbf{)}. \emph{Let }$\phi
\left( x\right) $\emph{\ be the function defined in (\ref{200}). Assume that 
}$\phi \left( x\right) \leq 0,\forall x\in \left[ 0,1\right] .$\emph{\ Then
our CIP has at most one solution.}

\textbf{Proof}. By Theorem 2.2 one can analytically extend the function $%
g_{0}\left( k\right) $ from the interval $k\in \lbrack \underline{k},%
\overline{k}]$ in the half plane $\mathbb{C}_{\gamma }=\left\{ z\in \mathbb{C%
}:\func{Im}z<\gamma \right\} ,$ where $\gamma =\gamma \left( \underline{k},%
\overline{k},\phi \right) $ is a positive number$.$ Hence, it follows from (%
\ref{2.37}) and (\ref{2.38}) that we know functions $u\left( 0,k\right)
,u_{x}\left( 0,k\right) $ for all $k\in \mathbb{R}.$ Consider the inverse
Fourier transform of the function $u$ with respect to $k$, $\mathcal{F}%
^{-1}\left( u\right) .$ It was established in the proof of Theorem 2.2 that
this transform can indeed be applied to the function $u$ and $\mathcal{F}%
^{-1}\left( u\right) =\widehat{u}\left( x,t\right) $, where the function $%
\widehat{u}\left( x,t\right) $ is the solution of the problem (\ref{2.15}), (%
\ref{2.16}). \ Functions 
\begin{equation}
\widehat{u}\left( 0,t\right) =\mathcal{F}^{-1}\left( u\left( 0,k\right)
\right) \text{ and }\widehat{u}_{x}\left( 0,t\right) =\mathcal{F}^{-1}\left(
u_{x}\left( 0,k\right) \right)  \label{2.39}
\end{equation}%
are known. Hence, we have obtained the inverse problem for equation (\ref%
{2.15}) with initial conditions (\ref{2.16}) and the data (\ref{2.39}). It
is well known that this inverse problem has at most one solution, see, e.g.
Chapter 2 of \cite{Rom}. $\square $

\section{A Version of the Quasi-Reversibility Method}

\label{sec:3}

As it was pointed out in section 1, we solve ordinary differential equations
with over determined boundary conditions using the QRM. In this section, we
develop the QRM for an arbitrary linear ordinary differential equation of
the second order with over determined boundary conditions. Below for any
Hilbert space $H$ its scalar product is $\langle \cdot ,\cdot \rangle _{H}$
. For convenience, we use in this section the notation \textquotedblleft $w$%
" for a generic complex valued function, which is irrelevant to the function 
$w$ in (\ref{2.37}).

Let functions $a\left( x\right) $ and $b\left( x\right) $ be in $C([0,1],%
\mathbb{C})$ and the function $d\left( x\right) $ be in $L^{2}([0,1],\mathbb{%
C})$. Let $p_{0}$ and $p_{1}$ be complex numbers. In this section we
construct an approximate solution of the following problem:%
\begin{equation}
\left\{ 
\begin{array}{c}
Lw=w^{\prime \prime }+a\left( x\right) w^{\prime }+b\left( x\right)
w=d\left( x\right) ,x\in \left( 0,1\right) , \\ 
w\left( 0\right) =p_{0},w^{\prime }\left( 0\right) =p_{1},w^{\prime }\left(
1\right) =0.%
\end{array}%
\right.  \label{3.1}
\end{equation}%
The QRM for problem (\ref{3.1}) amounts to the minimization of the following
functional 
\begin{equation*}
J_{\alpha }(w)=\frac{1}{2}\left( \Vert Lw-d\Vert _{L^{2}(\Omega
)}^{2}+\alpha \Vert w\Vert _{H^{3}(0,1)}^{2}\right)
\end{equation*}%
where 
\begin{equation*}
w\in W:=\{w\in H^{3}(0,1):w(0)=p_{0}\mbox{~and }w^{\prime
}(0)=p_{1},w^{\prime }(1)=0\}.
\end{equation*}
Below in this section we will establish existence and uniqueness of the
minimizer of $J_{\alpha }$. We will also show how close that minimizer is to
the solution of \eqref{3.1} if it exists. We start from the Carleman
estimate for the operator $d^{2}/dx^{2}.$

\subsection{Carleman estimate for the operator $d^{2}/dx^{2}$}

\label{sec:3.1}

\begin{lemma}[Carleman estimate]
\label{lem Car} For any complex valued function $u\in H^{2}\left( 0,1\right) 
$ with $u(0)=u^{\prime }(0)=0$ and for any parameter $\lambda >1$ the
following Carleman estimate holds 
\begin{equation}
\dint\limits_{0}^{1}\left\vert u^{\prime \prime }\right\vert
^{2}e^{-2\lambda x}dx\geq C\left[ \lambda \dint\limits_{0}^{1}|u^{\prime
}|^{2}e^{-2\lambda x}dx+\lambda ^{3}\dint\limits_{0}^{1}|u|^{2}e^{-2\lambda
x}dx\right] ,  \label{3.10}
\end{equation}%
where $C$ is a constant independent of $u$ and $\lambda .$
\end{lemma}

\textbf{Proof}. Since $C^{2}([0,1])$ is dense in $H^{2}(0,1)$, it is
sufficient to prove (\ref{3.10}) only for functions $u\in C^{2}\left[ 0,1%
\right] .$ Moreover, without loss of the generality, we can assume that $u$
is real valued.

Introduce the function $v=ue^{-\lambda x},\quad x\in (0,1).$ We have 
\begin{equation*}
u^{\prime }=(v^{\prime }+\lambda v),\quad u^{\prime \prime }=(v^{\prime
\prime }+2\lambda v^{\prime 2}v)e^{\lambda x}.
\end{equation*}%
A simple calculation yields 
\begin{align}
\frac{\left( u^{\prime \prime }\right) ^{2}e^{-2\lambda x}}{x+1}& =\frac{%
\left[ 2\lambda v^{\prime }+\left( v^{\prime \prime }+\lambda ^{2}v\right) %
\right] ^{2}}{x+1}\geq \frac{4\lambda v^{\prime }\left( v^{\prime \prime
}+\lambda ^{2}v\right) }{x+1}  \notag \\
& \geq \frac{d}{dx}\left( \frac{2\lambda \left( v^{\prime }\right)
^{2}+2\lambda ^{3}v^{2}}{x+1}\right) +\frac{2\lambda }{(x+1)^{2}}(\left(
v^{\prime }\right) ^{2}+\lambda ^{2}v^{2}).  \label{3.3}
\end{align}%
One the other hand, we have 
\begin{equation}
\left( v^{\prime }\right) ^{2}=\left( u^{\prime }-\lambda u\right)
^{2}e^{-2\lambda x}=\left[ \left( u^{\prime }\right) ^{2}-2\lambda u^{\prime
}u+\lambda ^{2}u^{2}\right] e^{-2\lambda x}\geq \left( \frac{1}{3}\left(
u^{\prime }\right) ^{2}-\frac{\lambda ^{2}}{2}u^{2}\right) e^{-2\lambda x}.
\label{3.44}
\end{equation}%
Combine \eqref{3.3} and \eqref{3.44} and then integrate the result. We
obtain 
\begin{equation*}
\int_{0}^{1}\frac{\left( u^{\prime \prime }\right) ^{2}e^{-2\lambda x}}{x+1}%
dx\geq \int_{0}^{1}\frac{2\lambda }{(x+1)^{2}}\left( \frac{1}{3}\left(
u^{\prime }\right) ^{2}+\frac{\lambda ^{2}}{2}u^{2}\right) e^{-2\lambda x}dx.
\end{equation*}%
Inequality \eqref{3.10} follows. $\square $

\subsection{The existence and uniqueness of the minimizer of $J_\protect%
\alpha$}

\label{sec:3.2}

We first establish the existence and uniqueness for the minimizer of $%
J_{\alpha}.$

\begin{theorem}
\label{thm unique} For every $\alpha \in \left( 0,1\right) $ there exists a
unique minimizer $w_{\alpha }\in W$ of the functional $J_{\alpha }.$
Furthermore, the following estimate holds%
\begin{equation}
\left\Vert w_{\alpha }\right\Vert _{H^{3}\left( 0,1\right) }\leq \frac{C_{1}%
}{\sqrt{\alpha }}\left( \left\vert p_{0}\right\vert +\left\vert
p_{1}\right\vert +\left\Vert d\right\Vert _{L^{2}\left( 0,1\right) }\right) ,
\label{3.4}
\end{equation}%
where the number $C_{1} >0$ depends only on $\|a\|_{L^{\infty}(0, 1)}$ and $%
\|b\|_{L^{\infty}(0, 1)}$.
\end{theorem}

Let $\{w_{n}\}_{n\geq 1}\subset W$ be a minimizing sequence of $J_{\alpha }$%
. That means, 
\begin{equation}
J_{\alpha }(w_{n})\rightarrow \inf_{W}J_{\alpha },\quad \mbox{as }%
n\rightarrow \infty .
\end{equation}%
It is not hard to see that $\{w_{n}\}_{n\geq 1}$ is bounded in $H^{3}(0,1)$.
In fact, if $\{w_{n}\}_{n\geq 1}$ has a unbounded subsequence then 
\begin{equation}
\inf_{W}J_{\alpha }\geq \limsup_{n\rightarrow \infty }\frac{\alpha }{2}\Vert
w_{n}\Vert _{H^{3}(0,1)}^{2}=\infty .
\end{equation}%
Without the loss of generality, we can assume that $\{w_{n}\}_{n\geq 1}$
weakly converges in $H^{3}(0,1)$ to a function $w_{\alpha }$ and strongly
converges to $w_{\alpha }$ in $H^{2}(0,1)$. The function $w_{\alpha }$
belongs to $W$ because $W$ is close and convex. We have 
\begin{align*}
J(w_{\alpha })& =\frac{1}{2}\Vert L(w_{\alpha })-d\Vert _{L^{2}(0,1)}^{2}+%
\frac{\alpha }{2}\Vert w_{\alpha }\Vert _{H^{3}(0,1)}^{2} \\
& \leq \overline{\lim_{n\rightarrow \infty }}\left( \frac{1}{2}\Vert
L(w_{n})-d\Vert _{L^{2}(0,1)}^{2}+\frac{\alpha }{2}\Vert w_{n}\Vert
_{H^{3}(0,1)}^{2}\right) =\inf_{W}J_{\alpha }.
\end{align*}%
The uniqueness of $w_{\alpha }$ is due to the strict convexity of $J_{\alpha
}$.

Inequality \eqref{3.4} can be verified by the fact that $J_{\alpha
}(w_{\alpha })\leq J_{\alpha }(v)$ where $v(x)=\chi (x)(p_{0}+xp_{1})\in W$
and the function $\chi \in C^{3}\left[ 0,1\right] $ is such that%
\begin{equation*}
\chi (x)=\left\{ 
\begin{array}{ll}
1 & x\in \lbrack 0,1/2], \\ 
0 & x\in \lbrack 3/4,1].%
\end{array}%
\right.
\end{equation*}%
$\square $

Let $r_{1}$ and $r_{2}$ be the real and imaginary parts respectively of the
complex valued function $r$. Without confusing, we identify $r$ with the
pair of real valued functions $(r_{1},r_{2}).$ Define 
\begin{align*}
L_{1}(w)& =w_{1}^{\prime \prime }+a_{1}(x)w_{1}^{\prime
}+b_{1}(x)w_{1}+(-a_{2}(x)w_{2}^{\prime }-b_{2}(x)w_{2}), \\
L_{2}(w)& =w_{2}^{\prime \prime }+a_{2}(x)w_{2}^{\prime
}+b_{2}(x)w_{2}+(a_{2}(x)w_{1}^{\prime }+b_{2}(x)w_{1}).
\end{align*}%
Rewrite $J_{\alpha }$ as 
\begin{equation}
J_{\alpha }(w_{1},w_{2})=\frac{1}{2}\left( \Vert
L_{1}(w_{1},w_{2})-d_{1}\Vert _{L^{2}(0,1)^{2}}^{2}+\Vert
L_{2}(w_{1},w_{2})-d_{2}\Vert _{L^{2}(0,1)^{2}}^{2}+\alpha \Vert
(w_{1},w_{2})\Vert _{H^{3}(0,1)^{2}}^{2}\right) .
\end{equation}

In order to find $w_{\alpha }$, we find the zero of the Fr\'{e}chet
derivative $DJ_{\alpha }$ of $J_{\alpha }(w_{\alpha })$. Since the operators 
$L_{1}$ and $L_{2}$ are linear, $DJ_{\alpha }$ is given by 
\begin{multline*}
DJ_{\alpha }(w_{1},w_{2})(h_{1},h_{2})=\langle L_{1}^{\ast
}(L_{1}(w_{1},w_{2})-d_{1}),(h_{1},h_{2})\rangle _{L^{2}(0,1)} \\
+\langle L_{2}^{\ast }(L_{2}(w_{1},w_{2})-d_{2}),(h_{1},h_{2})\rangle
_{L^{2}(0,1)}+\alpha \langle (w_{1},w_{2}),(h_{1},h_{2})\rangle _{H^{3}(0,1)}
\end{multline*}%
for all $h=h_{1}+ih_{2}\in H^{3}(0,1)$ with $h(0)=h^{\prime }(0)=h^{\prime
}(1)=0.$ The existence of a zero of $DJ_{\alpha }$ follows from the above
existence of $w_{\alpha }$ and the uniqueness is, again, deduced from the
convexity of $J_{\alpha }$ and $W$. In our computations, we use the finite
difference method to approximate the equation $DJ_{\alpha }(w_{1},w_{2})=0$,
together with the condition $w_{1}+iw_{2}\in W$, as a linear system for $%
w_{\alpha }$. The minimizer $w_{\alpha }$ of $J_{\alpha }$ is called the
regularized solution of \eqref{3.1} \cite{BK1,T}.

\subsection{Convergence of regularized solutions}

\label{sec:3.3}

While Theorem 3.1 claims the existence and uniqueness of the regularized
solution $w_{\alpha }$ of problem (\ref{3.1}), we now prove convergence of
regularized solutions to the exact solution of this problem, provided that
the latter solution exists, see \cite{BK1,T} for the definition of the
regularized solution. It is well known that one of concepts of the Tikhonov
regularization theory is the \emph{a priori} assumption about the existence
of an exact solution of an ill-posed problem, i.e. solution with noiseless
data \cite{BK1,T}. Estimate (\ref{3.4}) is valid for the $H^{3}\left(
0,1\right) -$norm and it becomes worse as long as $\alpha \rightarrow 0.$
However, Theorem 3.2 provides an estimate for the $H^{2}\left( 0,1\right) -$%
norm and the latter estimate is not worsening as $\alpha \rightarrow 0.$ To
prove Theorem 3.2, we use the Carleman estimate of subsection 3.1.

Suppose that there exists the exact solution $w^{\ast }$ of the problem (\ref%
{3.1}) with the exact (i.e. noiseless) data $d^{\ast } \in L^{2}$, $p_0^*,
p_1^* \in \mathbb{C}$. Let the number $\delta \in \left( 0,1\right) $ be the
level of the error in the data, i.e.%
\begin{equation}
\max\left\{\|d - d^*\|_{L^2(0, 1)}, |p_0 - p_0^*|, |p_1 - p_1^*|\right\}
\leq \delta  \label{3.12}
\end{equation}
Let $w_{\alpha }\in W$ be the unique minimizer of the functional $J_\alpha$,
which is guaranteed by Theorem \ref{thm unique}.

\begin{theorem}[convergence of regularized solutions]
Assume that (\ref{3.12}) holds. Then there exists a constant $C_{2}>0$
depending only on $\Vert a\Vert _{L^{\infty }(0,1)}$ and $\Vert b\Vert
_{L^{\infty }(0,1)}$ such that the following estimate holds%
\begin{equation}
\left\Vert w_{\alpha }-w^{\ast }\right\Vert _{H^{2}\left( 0,1\right) }\leq
C_{2}\left( \delta +\sqrt{\alpha }+\sqrt{\alpha }\left[ w^{\ast }\right]
\right) .  \label{3.13}
\end{equation}%
In particular, if $\alpha =\delta ^{2},$ then the following convergence rate
of regularized solutions $w_{\alpha }$ takes place (with a different
constant $C_{2}$)%
\begin{equation}
\left\Vert w_{\alpha }-w^{\ast }\right\Vert _{H^{2}\left( 0,1\right) }\leq
C_{2}\left( 1+\left[ w^{\ast }\right] \right) \delta .  \label{30}
\end{equation}%
\label{thm convergence of regularized solutions}
\end{theorem}

\textbf{Proof}. In this proof, $C_{2}$ is a generic constant depending only
on $\Vert a\Vert _{L^{\infty }(0,1)}$ and $\Vert b\Vert _{L^{\infty }(0,1)}.$
Let the function $\chi \in C^{2}\left[ 0,1\right] $ satisfies the following
condition 
\begin{equation*}
\chi (x)=\left\{ 
\begin{array}{ll}
1 & x\in \lbrack 0,1/2], \\ 
0 & x\in \lbrack 3/4,1].%
\end{array}%
\right.
\end{equation*}%
Define the \textquotedblleft error" function 
\begin{equation*}
\mathcal{E}(x)=\chi (x)\left( (p_{1}-p_{1}^{\ast })x+(p_{0}-p_{0}^{\ast
})\right) .
\end{equation*}%
Obviously, 
\begin{equation}
\Vert L\mathcal{E}\Vert _{L^{2}(0,1)}\leq C_{2}\delta ,\quad \mbox{and }%
\Vert \mathcal{E}\Vert _{H^{3}(0,1)}\leq C_{2}\delta .  \label{3.1313}
\end{equation}%
Since $w_{\alpha }$ is the minimizer of $J_{\alpha },$ then we have for all $%
h\in H^{3}(0,1)$ with $h(0)=h^{\prime }(0)=h^{\prime }(1)=0$ 
\begin{equation}
\langle Lw_{\alpha }-d,Lh\rangle _{L^{2}(0,1)}+\langle Lh,Lw_{\alpha
}-d\rangle _{L^{2}(0,1)}+\alpha \langle w_{\alpha },h\rangle
_{H^{3}(0,1)}+\alpha \langle h,w_{\alpha }\rangle _{H^{3}(0,1)}=0.
\label{3.131313}
\end{equation}%
Since $w^{\ast }$ is a solution of \eqref{3.1}, then 
\begin{multline}
\langle Lw^{\ast }-d^{\ast },Lh\rangle_{L^{2}(0,1)}+\langle Lh,Lw^{\ast
}-d^{\ast }\rangle _{L^{2}(0,1)}+\alpha \langle w^{\ast },h\rangle
_{H^{3}(0,1)}+\alpha \langle h,w^{\ast }\rangle _{H^{3}(0,1)} \\
= \alpha \langle w^{\ast },h\rangle _{H^{3}(0,1)}+\alpha \langle h,w^{\ast
}\rangle _{H^{3}(0,1)}.  \label{3.1414}
\end{multline}
%
Denoting $v=w_{\alpha }-w^{\ast }-\mathcal{E}$, using the test function $h=v$%
, and using the Cauchy-Schwarz inequality, we derive in a standard way from (%
\ref{3.12}), \eqref{3.131313} and \eqref{3.1414} that 
\begin{equation}
\Vert L(v)\Vert _{L^{2}(0,1)}^{2}+\alpha \Vert v\Vert _{H^{3}(0,1)}^{2}\leq
C_{2}\left( \delta ^{2}+\alpha \Vert w^{\ast }\Vert _{H^{3}(0,1)}^{2}\right)
.  \label{3.16}
\end{equation}%
On the other hand, Lemma \ref{lem Car} gives 
\begin{align*}
\Vert L(v)\Vert _{L^{2}(0,1)}^{2}& \geq \int_{0}^{1}|L(v)|^{2}e^{-2\lambda
x}dx\geq \int_{0}^{1}|v^{\prime \prime 2}e^{-2\lambda
x}dx-C_{2}\int_{0}^{1}(|v^{\prime 2}+|v|^{2})e^{-2\lambda x}dx \\
& \geq \frac{1}{2}\int_{0}^{1}\left\vert v^{\prime \prime }\right\vert
^{2}e^{-2\lambda x}dx+C_{2}\left[ (\lambda -1)\dint\limits_{0}^{1}|v^{\prime
}|^{2}e^{-2\lambda x}dx+(\lambda
^{3}-1)\dint\limits_{0}^{1}|v|^{2}e^{-2\lambda x}dx\right] .
\end{align*}%
Choosing $\lambda $ sufficiently large, we obtain 
\begin{equation*}
\Vert L(v)\Vert _{L^{2}(0,1)}^{2}\geq C_{2}e^{-2\lambda }\Vert v\Vert
_{H^{2}(0,1)}\geq C_{2}e^{-2\lambda }(\Vert w_{\alpha }-w^{\ast }\Vert
_{H^{2}(0,1)}^{2}-\delta ^{2}).
\end{equation*}%
Combining this and \eqref{3.16} completes the proof. $\square $

\section{Globally Convergent Numerical Method}

\label{sec:4}

\subsection{Integral differential equation}

\label{sec:4.1}

\begin{lemma}
\label{lem nonzero} Fix $x_{0}<0$ and $k>0$. Let $u\left( x,x_{0},k\right) $
be the solution of the problem (\ref{2.4})-(\ref{2.6}). Then $u\left(
x,x_{0},k\right) \neq 0,\forall x>x_{0}.$ In particular, $g_{0}\left(
k\right) \neq 0,\forall k\in \lbrack \underline{k},\overline{k}],$ where $%
g_{0}\left( k\right) $ is the function defined in (\ref{2.7}).
\end{lemma}

\textbf{Proof}. The proof is similar to that of the uniqueness of the
problem \eqref{2.4}, \eqref{2.6} in Theorem \ref{thm unique forward}. $%
\square$

Below $w\left( x,k\right) $ is the function defined in (\ref{2.37}). Since
by Theorem 2.2 $w\left( x,k\right) \neq 0,\forall x,$ then we can consider $%
\log w\left( x,k\right) .$ Since $\log z=\ln \left\vert z\right\vert +i\arg
z,\forall z\in \mathbb{C},z\neq 0,$ then the natural question is about $\arg
w\left( x,k\right) =\func{Im}w\left( x,k\right) .$ Hence, we consider the
asymptotic behavior at $k\rightarrow \infty $ of the function $w\left(
x,k\right) $. Using (\ref{2.60}) and Theorem 2.2, we obtain for $x>x_{0}$ 
\begin{equation}
w(x,k)=c^{-1/4}(x)\exp \left[ -ik\left( \displaystyle\int_{x_{0}}^{x}\sqrt{%
c(\xi )}d\xi -(x-x_{0})\right) \right] \left( 1+O\left( \frac{1}{k}\right)
\right) .  \label{4.0}
\end{equation}%
Obviously $\arg \left( 1+O(1/k)\right) \in \left[ -\pi ,\pi \right] $ for
sufficiently large $k>0.$ Hence, we set for sufficiently large $k>0$ and for 
$x>0$%
\begin{equation}
\log w\left( x,k\right) =\ln \left\vert w\left( x,x_{0},k\right) \right\vert
-ik\left( \dint\limits_{x_{0}}^{x}\sqrt{c\left( \xi \right) }d\xi
-x+x_{0}\right) +i\arg \left( 1+O\left( \frac{1}{k}\right) \right) .\text{ }
\label{4.01}
\end{equation}%
The function $\log w\left( x,k\right) $ is defined via (\ref{4.01}) for
sufficiently large $k.$ On the other hand, for not large values of $k$ it
would be better to work with derivatives of $\log w\left( x,k\right) .$
Indeed, we would not have problems then with defining $\arg w\left(
x,k\right) .$ Hence, taking $\overline{k}$ sufficiently large, we define the
function $\phi(x, k) = \log w(x, k)$ as 
\begin{equation}
\phi(x,k)=-\dint\limits_{k}^{\overline{k}}\frac{\partial_k w(x,\kappa )}{%
w(x,\kappa )}d\kappa +\log w(x,\overline{k}).  \label{4.2222}
\end{equation}%
Differentiate \eqref{4.2222} with respect to $k$. We have 
\begin{equation*}
\partial_kw(x, k) - w(x, k) \partial_k\phi(x, k) = 0.
\end{equation*}
Multiplying both sides of the equation above by $\exp(-\phi(x, k))$ gives 
\begin{equation*}
\partial_k\left( e^{-\phi(x,k)}w(x,k)\right) =0.
\end{equation*}%
Since $\phi(x,\overline{k})=\log w(x,\overline{k})$, then 
\begin{equation}
w(x,k)=e^{\phi(x,k)}.  \label{def log w}
\end{equation}%
The function $\phi(x,k)$, therefore, defines $\log w(x,k)$.

For each $k > 0$, define 
\begin{equation}
v(x, k) = \frac{\log w(x, k)}{k^2}.  \label{4.2}
\end{equation}

%

\begin{remark}
It follows from (\ref{2.7}), (\ref{2.37}) and (\ref{4.0}) that $g_{0}\left(
k\right) =1+O\left( 1/k\right) $ as $k\rightarrow \infty .$ Hence, $\arg
g_{0}\left( k\right) \in \left[ -\pi ,\pi \right] $ for sufficiently large $%
k>0.$ Hence, the function $\log g_{0}\left( k\right) $ can be defined
similarly with the function $v\left( x,k\right) $ in (\ref{4.2}).
\end{remark}

Let 
\begin{equation}
q\left( x,k\right) =\partial _{k}v\left( x,k\right) ,\quad x\in \lbrack
0,1],k\in \lbrack \underline{k},\overline{k}].  \label{4.1}
\end{equation}%
Hence,%
\begin{equation}
v\left( x,k\right) =-\dint\limits_{k}^{\overline{k}}q\left( x,\tau \right)
d\tau +v\left( x,\overline{k}\right) .  \label{4.20}
\end{equation}%
Denote%
\begin{equation}
V\left( x\right) =v\left( x,\overline{k}\right) .  \label{4.21}
\end{equation}%
We call $V$ the \textquotedblleft tail function" and this function is
unknown. Note that we do not use below the function $V$. Rather we use only
its $x-$derivatives. Hence, when using these derivatives, we are not
concerned with $\arg w\left( x,\overline{k}\right) .$

It easily follows from (\ref{2.3}), (\ref{2.4}), (\ref{2.60}), (\ref{2.37})
and (\ref{2.38}) that 
\begin{eqnarray}
w^{\prime \prime }-2ikw^{\prime }+k^{2}\beta \left( x\right) w &=&0,x\in
\left( 0,1\right) ,  \label{4.3} \\
w\left( 0,k\right) &=&g_{0}\left( k\right) ,w_{x}\left( 0,k\right)
=g_{1}\left( k\right) ,w_{x}\left( 1,k\right) =0.  \label{4.4}
\end{eqnarray}%
Using (\ref{def log w}) and (\ref{4.3}), we obtain 
\begin{equation}
v^{\prime \prime }+k^{2}\left( v^{\prime }\right) ^{2}-2ikv^{\prime }=-\beta
\left( x\right) .  \label{4.40}
\end{equation}%
Therefore, (\ref{4.1})-(\ref{4.40}) imply that 
\begin{equation}
q^{\prime \prime }-2ikq^{\prime }+2k^{2}q^{\prime }\left( -\dint\limits_{k}^{%
\overline{k}}q^{\prime }\left( x,\tau \right) d\tau +V^{\prime }\right)
-2i\left( -\dint\limits_{k}^{\overline{k}}q^{\prime }\left( x,\tau \right)
d\tau +V^{\prime }\right)  \label{4.5}
\end{equation}%
\begin{equation*}
+2k\left( -\dint\limits_{k}^{\overline{k}}q^{\prime }\left( x,\tau \right)
d\tau +V^{\prime }\right) ^{2}=0,\text{ }x\in \left( 0,1\right) ,
\end{equation*}%
\begin{equation}
q\left( 0,k\right) =\frac{\partial }{\partial k}\left( \frac{\log g_{0}(k)}{%
k^{2}}\right) ,q_{x}\left( 0,k\right) =\frac{\partial }{\partial k}\left[ 
\frac{2i}{k}\left( 1-\frac{1}{g_{0}\left( k\right) }\right) \right]
,q_{x}\left( 0,k\right) =0.  \label{4.6}
\end{equation}

We have obtained an integral differential equation (\ref{4.5}) for the
function $q$ with the overdetermined boundary data (\ref{4.6}). The tail
function in (\ref{4.5}) is also unknown. Hence, to approximate both
functions $q$ and $V$, we need to use not only conditions (\ref{4.5}), (\ref%
{4.6}) but something else as well. Thus, in our iterative procedure, we
solve problem (\ref{4.5}), (\ref{4.6}), assuming that $V$ is known, and
update the function $q$ this way. Then we update the unknown coefficient $%
\beta \left( x\right) .$ Next, we solve problem (\ref{4.3}), (\ref{4.4}) for
the function $w$ at $k:=\overline{k}$ and update the tail function via (\ref%
{4.1}), (\ref{4.20}) and (\ref{4.21}).

\subsection{Initial approximation $V_{0}\left( x \right) $ for the tail
function}

\label{sec:4.2}

It is important for the above iterative process to properly choose the
initial approximation $V_{0}\left( x\right) $ for the tail function. Since
we want to construct a globally convergent method, this choice must not use
any advanced knowledge of a small neighborhood of the exact solution $%
c^{\ast }\left( x\right) $ of our inverse problem.

We now describe how do we choose the initial tail $V_{0}\left( x\right) .$
It follows from Theorem \ref{thm assymptotic u} and the definition of $V(x)$
via (\ref{4.0})-\eqref{4.21} that there exists a function $p\left( x\right)
\in C^{2}\left[ 0,1\right] $ such that 
\begin{equation}
V\left( x,k\right) =\frac{p\left( x\right) }{k}+O\left( \frac{1}{k^{2}}%
\right) ,\text{ }q\left( x,k\right) =-\frac{p\left( x\right) }{k^{2}}%
+O\left( \frac{1}{k^{3}}\right) ,k\rightarrow \infty ,x>0,  \label{4.60}
\end{equation}%
Hence, assuming that the number $\overline{k}$ is sufficiently large, we
drop terms $O\left( 1/\overline{k}^{2}\right) $ and $O\left( 1/\overline{k}%
^{3}\right) $ in \eqref{4.60} and set%
\begin{equation}
V_{0}\left( x,k\right) =\frac{p\left( x\right) }{k},q^{0}\left( x,k\right) =-%
\frac{p\left( x\right) }{k^{2}},k\geq \overline{k},x>0.  \label{4.7}
\end{equation}%
Set $k:=\overline{k}$ in (\ref{4.5}) and (\ref{4.6}) and then substitute (%
\ref{4.7}) there. We obtain $p^{\prime \prime }=0.$ It follows from this, %
\eqref{2.38}, \eqref{4.1}, \eqref{4.4} and \eqref{4.7} that 
\begin{align}
V_{0}^{\prime \prime }& =0, \quad \mbox{in } (0, 1),  \label{guess v''} \\
V_{0}(0)& =\frac{\log g_{0}(\overline{k})}{\overline{k}^{2}},V_{0}^{\prime
}(0)=\frac{g_{1}(\overline{k})}{\overline{k}^{2}g(\overline{k})}%
,V_{0}^{\prime }(1)=0.  \label{V0 boundary condition}
\end{align}%
We solve the problem (\ref{guess v''}), (\ref{V0 boundary condition}) via
the QRM. By the embedding theorem $H^{2}\left( 0,1\right) \subset C^{1}\left[
0,1\right] $ and $\left\Vert f\right\Vert _{C^{1}\left[ 0,1\right] }\leq
C\left\Vert f\right\Vert _{H^{2}\left( 0,1\right) },\forall f\in H^{2}\left(
0,1\right) ,$ where $C>0$ is a generic constant$.$ Recall that the function $%
g_{1}\left( k\right) $ in (\ref{V0 boundary condition}) is linked with the
function $g_{0}\left( k\right) $ as in (\ref{2.38}). Thus, Theorems \ref{thm
unique} and \ref{thm convergence of regularized solutions} lead to Theorem %
\ref{thm 4.1}. In this theorem, we use the entire interval $\left[ 
\underline{k},\overline{k}\right] $ rather than just $k=\overline{k}$ (in (%
\ref{4.11})) for brevity: since we will use this interval below.

\begin{theorem}
\label{thm 4.1} Let $c^{\ast }\left( x\right) $ satisfying conditions (\ref%
{2.1})-(\ref{2.3}) be the exact solution of our CIP. For $k\geq \overline{k}%
, $ let the exact tail $V^{\ast }\left( x,k\right) $ have the form (\ref{4.7}%
). Assume that for $k\in \left[ \underline{k},\overline{k}\right] $ 
\begin{equation}
\left\vert \log g_{0}\left( k\right) -\log g_{0}^{\ast }\left( k\right)
\right\vert \leq \delta ,\left\vert g_{0}\left( k\right) -g_{0}^{\ast
}\left( k\right) \right\vert \leq \delta ,\left\vert g_{0}^{\prime }\left(
k\right) -\left( g_{0}^{\ast }\right) ^{\prime }\left( k\right) \right\vert
\leq \delta ,  \label{4.11}
\end{equation}%
where $\delta >0$ is a sufficiently small number, which characterizes the
level of the error in the boundary data. Let the function $V_{0,a}\left(
x\right) \in H^{3}\left( 0,1\right) $ be the approximate solution of the
problem (\ref{guess v''})-(\ref{V0 boundary condition}) obtained via the QRM
with $\alpha =\delta ^{2}$. Then there exists a constants $C_{3}=C_{3}\left( 
\overline{k},c^{\ast }\right) >0$ depending only on $\overline{k}$ and $%
c^{\ast }$ such that 
\begin{equation*}
\left\Vert V_{0,a}\left( x,\overline{k}\right) -V^{\ast }\left( x,\overline{k%
}\right) \right\Vert _{C^{1}\left[ 0,1\right] }\leq C\left\Vert
V_{0,a}\left( x,\overline{k}\right) -V^{\ast }\left( x,\overline{k}\right)
\right\Vert _{H^{2}\left( 0,1\right) }\leq C_{3}\delta .
\end{equation*}
\end{theorem}

\begin{remark}
\label{remark 4.2}~

\begin{enumerate}
\item \textrm{Theorem \ref{thm 4.1} is valid only within the framework of a
quite natural approximation (\ref{4.7}), in which small terms }$O\left(
1/k^{2}\right) ,O\left( 1/k^{3}\right) $\textrm{\ of formulae (\ref{4.60})
are ignored. We use this approximation only to find the first tail and do
not use it in follow up iterations. We believe that the use of this
approximation is justified by the fact that the topic of the globally
convergent numerical methods for CIPs is a very challenging one. }

\item \textrm{Thus, it follows from Theorem \ref{thm 4.1} that our initial
tail function }$V_{0,a}\left( x,\overline{k}\right) $\textrm{\ provides a
good approximation for the exact tail }$V^{\ast }\left( x,\overline{k}%
\right) $\textrm{\ already at the start of our iterative process. Hence,
setting in (\ref{4.40}) }$k=\overline{k}$ and recalling \textrm{(\ref{4.21}),%
} we conclude that\textrm{\ the target coefficient }$c^{\ast }\left(
x\right) $\textrm{\ is also reconstructed with a good accuracy at the start
of our iterative process. The error of the approximation of both }$V^{\ast }$%
\textrm{\ and }$c^{\ast }$\textrm{depends only on the level }$\delta $%
\textrm{\ of the error in the boundary data. The latter is exactly what is
usually required when solving inverse problems. It is important that when
obtaining this approximation for }$V^{\ast }$\textrm{, we have not used any
advanced knowledge about a small neighborhood of the exact solution }$%
c^{\ast }$\textrm{$.$ In other words, the requirement of the global
convergence is in place (see Introduction for this requirement). }

\item \textrm{Even though we obtain good approximations for }$V^{\ast
}\left( x,k\right) $\textrm{\ and }$c^{\ast }\left( x\right) $\textrm{\ from
the start, our numerical experience tells us that results improve with
iterations in our iterative process described below. A similar observation
took place in the earlier above cited works of this group, where the Laplace
transform of the time dependent data was used. This is of course due to the
approximate nature of (\ref{4.7}). }

\item \textrm{Even though it is possible to sort of \textquotedblleft unite
in one" first two conditions (\ref{4.11}), we are not doing this here for
brevity. }

\item \textrm{In the convergence analysis, we use the form (\ref{4.7}) for
the functions }$V^{\ast }$ and $q^{\ast }$\textrm{\ only on the first
iteration, since this form of functions }$V$\textrm{$,$}$q$\textrm{\ is used
only on the first iteration of our algorithm. }
\end{enumerate}
\end{remark}

Below we consider the error parameter $\eta $ defined as%
\begin{equation}
\eta =h+\delta .  \label{4.110}
\end{equation}

\subsection{Numerical method}

\label{sec:4.3}

\subsubsection{Equations for $q_{n}$}

\label{sec:4.3.1}

Consider a partition of the frequency interval $\left[ \underline{k},%
\overline{k}\right] $ in $N$ subintervals with the step size $h$,%
\begin{equation*}
k_{N}=\underline{k}<k_{N-1}<...<k_{1}<k_{0}=\overline{k},k_{j-1}-k_{j}=h,
\end{equation*}%
where the number $h>0$ is sufficiently small. We assume that the function $%
q\left( x,k\right) $ is piecewise constant with respect to $k$, $q\left(
x,k\right) =q\left( x,k_{n}\right) $ for $k\in \left[ k_{n},k_{n-1}\right) .$
For each $n=1,\cdots ,N$ and for all $x\in (0,1)$ define 
\begin{align}
q_{0}(x)& =0,q_{n}(x)=q(x,k_{n}),  \label{4.12} \\
Q_{n-1}(x)& =\int_{k_{n-1}}^{\overline{k}}q(x,\kappa )d\kappa
=h\sum_{j=0}^{n-1}q_{j}(x).  \label{4.14}
\end{align}%
Hence, by (\ref{4.20})%
\begin{equation}
v^{\left( s\right) }\left( x,k_{n}\right) =-hq_{n}^{\left( s\right) }\left(
x\right) -Q_{n-1}^{\left( s\right) }\left( x\right) +V^{\left( s\right)
},s=1,2.  \label{4.13}
\end{equation}

Then (\ref{4.5}) and (\ref{4.12})-(\ref{4.13}) imply that for all $%
n=1,\cdots ,N$ 
\begin{multline*}
q_{n}^{\prime \prime }+2k\left( -hq_{n}^{\prime }-Q_{n-1}^{\prime
}+V^{\prime }\right) ^{2}+2k^{2}\left( -hq_{n}^{\prime }-Q_{n-1}^{\prime
}+V^{\prime }\right) q_{n}^{\prime } \\
-2i\left( -hq_{n}^{\prime }-Q_{n-1}^{\prime }+V^{\prime }\right)
-2ikq_{n}^{\prime }=0.
\end{multline*}%
Choose the step size $h$ sufficiently small and ignore terms with $h$ and $%
h^{2}$. Note that $k-k_{n}<h$ for $k\in k\in \left[ k_{n},k_{n-1}\right) .$
Also, we keep in mind that we will iterate with respect to tail functions
for each $n$ as well as with respect to $n$. Thus, we rewrite the last
equation as 
\begin{equation}
q_{n,j}^{\prime \prime }+\left[ 2k_{n}^{2}\left( -Q_{n-1}^{\prime
}+V_{n,j}^{\prime }\right) -2ik_{n}\right] q_{n,j}^{\prime }=-2k_{n}\left(
-Q_{n-1}^{\prime }+V_{n,j}^{\prime }\right) ^{2}+2i\left( -Q_{n-1}^{\prime
}+V_{n,j}^{\prime }\right)  \label{4.15}
\end{equation}%
for all $x\in (0,1),j=1,\cdots ,m$ for some $m>0$. The boundary conditions
for $q_{n,j}$ in \eqref{4.15} are taken according to those for $q$ in %
\eqref{4.6}. Precisely, 
\begin{equation}
\begin{array}{rcl}
q_{n,j}\left( 0\right) & = & \displaystyle\frac{1}{h}\left( \frac{\log
g_{0}(k_{n-1})}{k_{n-1}^{2}}-\frac{\log g_{0}(k_{n})}{k_{n}^{2}}\right)
=\psi _{n}^{0}, \\ 
q_{n,j}^{\prime }(0) & = & \displaystyle\frac{1}{h}\left( \frac{2i}{k_{n-1}}%
\left( 1-\frac{1}{g_{0}\left( k_{n-1}\right) }\right) -\frac{2i}{k_{n}}%
\left( 1-\frac{1}{g_{0}\left( k_{n}\right) }\right) \right) =\psi _{n}^{1},
\\ 
q_{n,j}^{\prime }(1) & = & \displaystyle0.%
\end{array}
\label{4.16}
\end{equation}

\subsubsection{The algorithm}

\label{sec:4.3.2}

The procedure to solve the CIP is described below:

\begin{algorithm}
We reconstruct a set $\{\beta_1, \cdots, \beta_N\}$ of approximations for $%
\beta^*.$

\begin{enumerate}
\item Set $q_{0}\equiv 0.$ Find the first approximation $V_{0}$ for the tail
function solving the problem (\ref{guess v''})-(\ref{V0 boundary condition})
via the QRM.

\item \label{step 2} For an integer $n\in \left[ 1,N\right] $, suppose that
functions $q_{0},\cdots ,q_{n-1}$, $V_{0}^{\prime },\cdots ,V_{n-1}^{\prime
} $, $\beta _{0},\cdots ,\beta _{n-1}$ are known. Therefore, $Q_{n-1}$ is
known. We calculate the function $\beta _{n}$ as follows.

\begin{enumerate}
\item Set $V_{n,1}^{\prime }=V_{n-1}^{\prime }$, $V_{n,1}^{\prime \prime
}=V_{n-1}^{\prime \prime },$ $\beta _{n,0}=\beta _{n-1}$.

\item For $j=1,\cdots ,m$:

\begin{enumerate}
\item Solve the problem \eqref{4.15}, \eqref{4.16} for $q_{n,j}$ via the QRM.

\item For $s=1,2,$ let $v_{n,j}^{\left( s\right) }=-hq_{n,j}^{\left(
s\right) }\left( x\right) -Q_{n-1}^{\left( s\right) }\left( x\right)
+V_{n,j}^{\left( s\right) }$ due to an analog of \eqref{4.13}.

\item \label{step 2b3}Calculate $\beta _{n,j}$ by \eqref{4.2323}, which will
be explained later.

\item Solve the problem (\ref{4.3}), (\ref{4.4}) via the QRM with $k:=%
\overline{k}$ and $\beta \left( x\right) :=\beta _{n,j}\left( x\right) $.
Let $w_{n,j}\left( x,\overline{k}\right) $ be its solution. Next, using (\ref%
{4.2}) and (\ref{4.20}), set 
\begin{equation}
V_{n,j+1}^{\prime }\left( x\right) =\frac{1}{\overline{k}^{2}}\frac{%
w_{n,j}^{\prime }(x,\overline{k})}{w_{n,j}(x,\overline{k})},\text{ }%
V_{n,j+1}^{\prime \prime }\left( x\right) =\frac{1}{\overline{k}^{2}}\frac{%
w_{n,j}^{\prime \prime }(x,\overline{k})}{w_{n,j}(x,\overline{k})}-\frac{1}{%
\overline{k}^{2}}\frac{\left( w_{n,j}^{\prime }(x,\overline{k})\right) ^{2}}{%
w_{n,j}^{2}(x,\overline{k})}.  \label{4.17}
\end{equation}
\end{enumerate}

\item Set $\beta _{n}=\beta _{n,j^{0}}$ where 
\begin{equation*}
j^{0}=\mathrm{argmin}\left\{ \frac{\Vert \beta _{n,j}-\beta _{n,j-1}\Vert
_{L^{2}(0,1)}}{\Vert \beta _{n,j}\Vert _{L^{2}(0,1)}},j=1,\cdots ,m\right\} .
\end{equation*}
\end{enumerate}

\item \label{step 3} Chose $\beta =\beta _{n^{0}}$ where 
\begin{equation*}
n^{0}=\mathrm{argmin}\left\{ \frac{\Vert \beta _{n}-\beta _{n-1}\Vert
_{L^{2}(0,1)}}{\Vert \beta _{n}\Vert _{L^{2}(0,1)}},n=1,\cdots ,N\right\} .
\end{equation*}
\end{enumerate}
\end{algorithm}

In the algorithm, all differential equations are solved via the QRM. Thus,
we keep for those \textquotedblleft QRM solutions" the same notations for
brevity. For simplicity, we assume here that $\beta \left( x\right) \geq 0,$
although we also work in one case of experimental data with a non-positive
function $\beta .$ Thus, in the algorithm above, we update the function $%
\beta \left( x\right) $ using (\ref{2.1}), (\ref{2.3}), (\ref{4.40}) and (%
\ref{4.13}) as%
\begin{equation}
\beta _{n,j}=\min \left\{ \max \left\{ \left\vert -v_{n,j}^{\prime \prime
}-k_{n}^{2}\left( v_{n,j}^{\prime }\right) ^{2}+2ik_{n}v_{n,j}^{\prime
}\right\vert ,c_{0}-1\right\} ,c_{1}-1\right\} .  \label{4.2323}
\end{equation}%
This truncation helps us to get a better accuracy in the reconstructed
function $\beta $. In fact, it follows from (\ref{2.1}), (\ref{2.3}), %
\eqref{4.40} and (\ref{4.2323}) that 
\begin{equation}
\left\vert \beta _{n,j}\left( x\right) -\beta ^{\ast }\left( x\right)
\right\vert \leq \left\vert -v_{n,j}^{\prime \prime }-k_{n}^{2}\left(
v_{n,j}^{\prime }\right) ^{2}+2ik_{n}v_{n,j}^{\prime }-\beta ^{\ast }\left(
x\right) \right\vert .  \label{4.22}
\end{equation}%
where $\beta ^{\ast }=c^{\ast }-1$ is the exact solution of the CIP.

\section{Global Convergence}

\label{sec:5}

In this section we prove our main result about the global convergence of the
algorithm of the previous section. This method actually has the
approximately global convergence property, see the third paragraph of
Section \ref{sec:1} and Remarks \ref{remark 4.2}. For brevity we assume in
this section that $j^{0}=1$ in the above algorithm. In other words, we
assume that we do not perform inner iterations. The case $j^{0}>1$ can be
done similarly.

First, we need to introduce some assumptions about the exact solution.
Everywhere below the superscript \textquotedblleft $\ast $" denotes
functions which correspond to the exact coefficient $c^{\ast }\left(
x\right) $. Denote $q_{n}^{\ast }\left( x\right) =q^{\ast }\left(
x,k_{n}\right) .$ Then $q^{\ast }\left( x,k\right) =q_{n}^{\ast }\left(
x\right) +O\left( h\right) $ for $h\rightarrow 0$ and for $k\in \left[
k_{n},k_{n-1}\right) .$ Set $q_{0}^{\ast }\left( x\right) \equiv 0.$Let the
function $Q_{n-1}^{\ast }\left( x\right) $ be the same as in (\ref{4.14}),
except that functions $q_{j}$ are replaced with $q_{j}^{\ast }.$ Also, let 
\begin{equation}
\left( {v_{n}^{\ast }}\right) ^{\left( s\right) }=-h\left( {q_{n}^{\ast }}%
\right) ^{\left( s\right) }\left( x\right) -\left( {Q_{n-1}^{\ast }}\right)
^{\left( s\right) }\left( x\right) +\left( {V^{\ast }}\right) ^{\left(
s\right) },s=1,2.  \label{5.1}
\end{equation}%
Then (\ref{4.40}) and implies that 
\begin{equation}
\beta ^{\ast }\left( x\right) =-\left( v_{n}^{\ast }\right) ^{\prime \prime
}-k_{n}^{2}\left( \left( v_{n}^{\ast }\right) ^{\prime }\right)
^{2}+2ik_{n}\left( v_{n}^{\ast }\right) ^{\prime }+F_{n}^{\ast }\left(
x\right) ,  \label{5.2}
\end{equation}%
Also, by (\ref{4.15}) and (\ref{4.16}) we have for $x\in \left( 0,1\right) $%
\begin{eqnarray}
&&\left( q_{n}^{\ast }\right) ^{\prime \prime }+\left[ 2k_{n}^{2}\left(
-\left( Q_{n-1}^{\ast }\right) ^{\prime }+\left( V^{\ast }\right) ^{\prime
}\right) -2ik_{n}\right] \left( q_{n}^{\ast }\right) ^{\prime }  \notag \\
&=&-2k_{n}\left( -\left( Q_{n-1}^{\ast }\right) ^{\prime }+\left( V^{\ast
}\right) ^{\prime }\right) ^{2}+2i\left( -\left( Q_{n-1}^{\ast }\right)
^{\prime }+\left( V^{\ast }\right) ^{\prime }\right) +G_{n}^{\ast }\left(
x\right) ,  \label{5.3} \\
q_{n}^{\ast }\left( 0\right) &=&\psi _{n}^{\ast ,0},\left( q_{n}^{\ast
}\right) ^{\prime }\left( 0\right) =\psi _{n}^{\ast ,1},\left( q_{n}^{\ast
}\right) ^{\prime }\left( 1\right) =0.  \notag
\end{eqnarray}%
Since the number $\delta $ characterizes the error in the boundary data and
since $\eta >\delta $ is the error parameter introduced in (\ref{4.110}),
then, taking into account (\ref{4.16}), we set%
\begin{equation}
\left\vert \psi _{n}^{0}-\psi _{n}^{\ast ,0}\right\vert \leq \eta ,\text{ }%
\left\vert \psi _{n}^{1}-\psi _{n}^{\ast ,1}\right\vert \leq \eta .
\label{5.4}
\end{equation}%
In (\ref{5.2}) and (\ref{5.3}) $F_{n}^{\ast }\left( x\right) $ and $%
G_{n}^{\ast }\left( x\right) $ are error functions, which can be estimated
as 
\begin{equation}
\left\Vert F_{n}^{\ast }\right\Vert _{L^{2}\left( 0,1\right) }\leq M\eta ,%
\text{ }\left\Vert G_{n}^{\ast }\right\Vert _{L^{2}\left( 0,1\right) }\leq
M\eta ,  \label{5.5}
\end{equation}%
where $M>0$ is a constant. We also assume that 
\begin{eqnarray}
\left\Vert q_{n}^{\ast }\right\Vert _{C^{1}\left[ 0,1\right] } &\leq
&M,\left\Vert q_{n}^{\ast }\right\Vert _{H^{2}\left( 0,1\right) }\leq
M,\left\Vert \left( V^{\ast }\right) ^{\prime }\right\Vert _{C\left[ 0,1%
\right] }\leq M,\left\Vert \left( V^{\ast }\right) ^{\prime \prime
}\right\Vert _{L^{2}\left( 0,1\right) }\leq M,  \label{5.50} \\
\left\Vert \left( w^{\ast }\right) ^{\prime }\right\Vert _{C\left[ 0,1\right]
} &\leq &M,\left\Vert \left( w^{\ast }\right) ^{\prime \prime }\right\Vert
_{L^{2}\left( 0,1\right) }\leq M.  \label{5.501}
\end{eqnarray}

\textbf{Theorem 5.1}. \emph{Let conditions of Theorem 4.1 hold. In
procedures (i) and (iv) of the algorithm set in the QRM }$\alpha =\eta ^{2}.$%
\emph{\ Assume that the number }$\underline{k}>1$\emph{\ and the number }$%
\overline{k}$ \emph{is so large that in (\ref{4.0}) }$\left\vert O\left(
1/k\right) \right\vert <1/2$\emph{\ for }$k\geq \overline{k}$\emph{\ for }$%
c=c^{\ast }.$\emph{\ Let the function }$w^{\ast }\left( x,\overline{k}%
\right) \in C^{2}\left[ 0,1\right] $ \emph{be the solution of the problem (%
\ref{4.3}), (\ref{4.4}) with the exact coefficient }$\beta ^{\ast }\left(
x\right) =c^{\ast }\left( x\right) -1$\emph{\ and the exact data }$%
g_{0}^{\ast }\left( \overline{k}\right) ,g_{1}^{\ast }\left( \overline{k}%
\right) .$ \emph{Let $N_{1}$ be an integer in $[1,N].$ Then there exists a
sufficiently large constant }$M=M\left( c_{0},c_{1},\overline{k}\right) >1$%
\emph{\ for which estimates (\ref{5.5})-(\ref{5.501}) are valid and which
also satisfies }%
\begin{equation}
M>24\overline{k}^{2},M>16\sqrt{c_{1}},M>C_{3}  \label{5.6}
\end{equation}%
\emph{as well as a constant }$M_{1}=M_{1}\left( M\overline{k}^{2}\right) >0$%
\emph{\ such that if the error parameter }$\eta $\emph{\ is so small that }%
\begin{equation}
\eta \in \left( 0,\eta _{0}\right) \emph{,\ where}\text{ }\eta _{0}\leq 
\frac{1}{4c_{1}^{1/4}\left( M_{1}M^{14}\right) ^{2N_{1}}},\text{ }
\label{5.7}
\end{equation}%
\emph{then for }$n=1,2...,N_{1}$\emph{\ the following estimate holds true} 
\begin{equation}
\left\Vert \beta _{n}-\beta ^{\ast }\right\Vert _{L^{2}\left( 0,1\right)
}\leq \left( M_{1}M^{16}\right) ^{n}\eta <\sqrt{\eta }.  \label{5.10}
\end{equation}

\textbf{Remark 6.1. }\emph{Thus, this theorem claims that our iteratively
found functions }$\beta _{n}$\emph{\ are located in a sufficiently small
neighborhood of the exact solution }$\beta ^{\ast }$\emph{\ as long as }$%
n\in \left[ 1,N_{1}\right] .$\emph{\ Since this is achieved without any
advanced knowledge of a small neighborhood of the exact solution }$\beta
^{\ast }$\emph{, then Theorem 5.1 implies the global convergence of our
algorithm, see Introduction. On the other hand, this is achieved within the
framework of the approximation of subsection 4.2. Hence, to be more precise,
this is the approximate global convergence property, see section 1.1.2 of 
\cite{BK1} and section 4 of \cite{KSNF1} for the definition of this
property. Recall that the number of iterations (}$N_{1}$\emph{\ in our case)
can be considered sometimes as a regularization parameter in the theory of
ill-posed problems \cite{BK1,T}.}

\textbf{Proof}. In addition to (\ref{5.10}), we will also prove that for $%
n=1,2...,N_{1}$ 
\begin{eqnarray}
\left\Vert q_{n}-q_{n}^{\ast }\right\Vert _{C^{1}\left[ 0,1\right]
},\left\Vert q_{n}-q_{n}^{\ast }\right\Vert _{H^{2}\left( 0,1\right) } &\leq
&\left( M_{1}M^{16}\right) ^{n}\eta <\sqrt{\eta },  \label{5.8} \\
\left\Vert V_{n+1}^{^{\prime }}-\left( V^{\ast }\right) ^{\prime
}\right\Vert _{C\left[ 0,1\right] },\left\Vert V_{n+1}^{\prime \prime
}-\left( V^{\ast }\right) ^{\prime \prime }\right\Vert _{L^{2}\left(
0,1\right) } &\leq &\left( M_{1}M^{16}\right) ^{n}\eta <\sqrt{\eta }.
\label{5.9}
\end{eqnarray}%
To simplify and shorten the proof, we assume in this proof that we work only
with real valued functions. Hence, we replace in two terms of (\ref{4.15}) $%
``i"$ with \textquotedblleft 1" and similarly in two terms of (\ref{5.3}).
The case of complex valued functions is very similar. However, it contains
some more purely technical details and is, therefore, more space consuming.
We use the mathematical induction method. Denote%
\begin{equation}
\widetilde{q}_{n}=q_{n}-q_{n}^{\ast },\widetilde{V}_{n}=V_{n}-V^{\ast },%
\widetilde{v}_{n}=v_{n}-v_{n}^{\ast },\widetilde{Q}_{n-1}=Q_{n-1}-Q_{n-1}^{%
\ast },\widetilde{\beta }_{n}=\beta _{n}-\beta ^{\ast }.  \label{5.11}
\end{equation}%
Using Theorem 4.1, we obtain%
\begin{equation}
\left\Vert V_{1}^{\prime }-\left( V^{\ast }\right) ^{\prime }\right\Vert _{C 
\left[ 0,1\right] }\leq M\eta <M^{9}\eta ,\text{ }\left\Vert V_{1}^{\prime
\prime }-\left( V^{\ast }\right) ^{\prime \prime }\right\Vert _{L^{2}\left(
0,1\right) }\leq M\eta <M^{9}\eta .  \label{5.12}
\end{equation}%
Hence, by (\ref{5.50}) and (\ref{5.7})%
\begin{equation}
\left\Vert V_{1}^{\prime }\right\Vert _{C\left[ 0,1\right] }\leq
2M,\left\Vert V_{1}^{\prime \prime }\right\Vert _{L^{2}\left( 0,1\right)
}\leq 2M.  \label{5.13}
\end{equation}%
Following notations of section 3, denote%
\begin{eqnarray}
L_{n}\left( y\right) &=&y^{\prime \prime }+\left[ 2k_{n}^{2}\left(
-Q_{n-1}^{\prime }+V_{n}^{\prime }\right) -2k_{n}\right] y^{\prime },
\label{5.14} \\
L_{n,\ast }\left( y\right) &=&y^{\prime \prime }+\left[ 2k_{n}^{2}\left(
-\left( Q_{n-1}^{\prime }\right) ^{\ast }+V^{\ast }\right) -2ik_{n}\right]
y^{\prime }.  \notag
\end{eqnarray}%
Define 
\begin{equation}
W_{0} = \{\phi \in H^3(0, 1): \phi(0) = \phi^{\prime }(0) = \phi^{\prime
}(1) = 0\}.  \label{def W0}
\end{equation}
Since all functions $q_{n}$ are QRM solutions of corresponding problems with 
$\alpha =\eta ^{2}$, then, using (\ref{5.3}) and (\ref{4.15}), we obtain for
all functions $z\in W_0$,%
\begin{eqnarray}
\left( L_{n}q_{n},L_{n}z\right) +\eta ^{2}\left[ q_{n},z\right] &=&\left(
-2k_{n}\left( -Q_{n-1}^{\prime }+V_{n}^{\prime }\right) ^{2}+2\left(
-Q_{n-1}^{\prime }+V_{n}^{\prime }\right) ,L_{n}z\right) ,  \notag \\
\left( L_{n,\ast }q_{n}^{\ast },L_{n}z\right) +\eta ^{2}\left[ q_{n}^{\ast
},z\right] &=&\left( -2k_{n}\left( -\left( Q_{n-1}^{\ast }\right) ^{\prime
}+\left( V^{\ast }\right) ^{\prime }\right) ^{2},L_{n}z\right)  \label{5.15}
\\
+\left( 2\left( -\left( Q_{n-1}^{\ast }\right) ^{\prime }+\left( V^{\ast
}\right) ^{\prime }\right) ,L_{n}z\right) &&+\left( G_{n}^{\ast },z\right)
+\eta ^{2}\left[ q_{n}^{\ast },z\right] .  \notag
\end{eqnarray}%
Subtracting the second equality (\ref{5.15}) from the first one and using (%
\ref{5.14}) and, we obtain%
\begin{equation*}
\left( L_{n}\left( \widetilde{q}_{n}\right) ,L_{n}z\right) +\eta ^{2}\left[ 
\widetilde{q}_{n},z\right] =\left( \left( -2k_{n}^{2}\left( -\widetilde{Q}%
_{n-1}^{\prime }+\widetilde{V}_{n}^{\prime }\right) \left( q_{n}^{\ast
}\right) ^{\prime }\right) ,L_{n}z\right)
\end{equation*}%
\begin{eqnarray}
&&-\left( \left( 2k_{n}\left( -\widetilde{Q}_{n-1}^{\prime }+\widetilde{V}%
_{n}^{\prime }\right) \left( -Q_{n-1}^{\prime }-\left( Q_{n-1}^{\ast
}\right) ^{\prime }+V_{n}^{\prime }+\left( V^{\ast }\right) ^{\prime }-\frac{%
1}{k_{n}}\right) \right) ,L_{n}z\right)  \label{5.16} \\
&&-\left( G_{n}^{\ast },L_{n}z\right) -\eta ^{2}\left[ q_{n}^{\ast },z\right]
.  \notag
\end{eqnarray}%
In addition, by (\ref{5.4})%
\begin{equation}
\left\vert \widetilde{q}_{n}\left( 0\right) \right\vert \leq \eta
,\left\vert \widetilde{q}_{n}^{\prime }\left( 0\right) \right\vert \leq \eta
,\widetilde{q}_{n}^{\prime }\left( 1\right) =0.  \label{5.17}
\end{equation}

We now explain the meaning of the constant $M_{1}.$ Since the constant $%
C_{2} $ in Theorem 3.3 depends on $C-$norms of coefficients of the operator $%
L$ in (\ref{3.1}), we need to estimate from the above the $C-$norm of the
coefficient of the operator $L_{n}$ in (\ref{5.14}). If (\ref{5.8}) is true,
then using (\ref{5.50}) and (\ref{5.7}) and noting that by (\ref{5.11}) $%
\left\Vert q_{n}\right\Vert _{C^{1}\left[ 0,1\right] }\leq \left\Vert 
\widetilde{q}_{n}\right\Vert _{C^{1}\left[ 0,1\right] }+\left\Vert
q_{n}^{\ast }\right\Vert _{C^{1}\left[ 0,1\right] }\leq 1+M\leq 2M,$ we
obtain%
\begin{equation}
\left\Vert q_{n}\right\Vert _{C^{1}\left[ 0,1\right] }\leq 2M,\left\Vert
q_{n}\right\Vert _{H^{2}\left( 0,1\right) }\leq 2M.  \label{5.18}
\end{equation}%
Hence, by (\ref{4.14}) $\left\vert Q_{n-1}^{\prime }\right\vert \leq
2MNh=2M\left( \overline{k}-\underline{k}\right) \leq 2M\overline{k}.$ Hence,
if (\ref{5.9}) is also true, then the coefficient of the operator $L_{n}$
can be estimated as 
\begin{equation}
\left\vert 2k_{n}^{2}\left( -Q_{n-1}^{\prime }+V_{n}^{\prime }\right)
-2k_{n}\right\vert \leq 10M\overline{k}^{2}.  \label{5.19}
\end{equation}%
Thus, in the case of the operator $L_{n}$ in the analog of estimate (\ref{30}%
) for the QRM, the constant $C_{2}$ should be replaced with another constant 
$M_{1}=M_{1}\left( M\overline{k}^{2}\right) >0.$

First, consider the case $n=1$ and estimate functions $\widetilde{q}%
_{1},q_{1}.$ In this case $Q_{n-1}^{\prime }=Q_{0}=0$ and so (\ref{5.19}) is
an over-estimate of course. Still, to simplify the presentation, we use $%
M_{1}$ in this case. Estimate first two terms in the right hand side of (\ref%
{5.16}) at $n=1$. By Theorem 4.1, (\ref{5.5})-(\ref{5.6}), (\ref{5.11}) and (%
\ref{5.13}) 
\begin{equation}
\left\vert -2k_{1}^{2}\left( -\widetilde{Q}_{0}^{\prime }+\widetilde{V}%
_{1}^{\prime }\right) \left( q_{1}^{\ast }\right) ^{\prime }\right\vert \leq
2\overline{k}^{2}M^{2}\eta ,  \label{5.20}
\end{equation}%
\begin{equation}
\left\vert -Q_{0}^{\prime }-Q_{0}^{\ast \prime }+V_{1}^{\prime }+\left(
V^{\ast }\right) ^{\prime }-\frac{1}{k_{1}}\right\vert \leq
2M+1<3M,\left\Vert G_{1}^{\ast }\right\Vert _{L^{2}\left( 0,1\right) }+\eta
^{2}\left\Vert q_{n}^{\ast }\right\Vert _{H^{2}\left( 0,1\right) }\leq
2M\eta .  \label{5.21}
\end{equation}%
Hence, Theorem 3.3, (\ref{5.16}), (\ref{5.17}), (\ref{5.20}) and (\ref{5.21}%
) lead to%
\begin{equation}
\left\Vert \widetilde{q}_{1}\right\Vert _{H^{2}\left( 0,1\right)
},\left\Vert \widetilde{q}_{1}\right\Vert _{C^{1}\left[ 0,1\right] }\leq
M_{1}\left( 16\overline{k}^{2}M^{2}\right) \eta <M_{1}M^{5}\eta .
\label{5.22}
\end{equation}

Hence, (\ref{5.8}) is true for $n=1$. We now estimate derivatives of the
function $v_{1}.$ By (\ref{5.13}) and (\ref{5.18}) 
\begin{equation}
\left\vert v_{1}^{\prime }\right\vert \leq \left\vert -hq_{1}^{\prime
}\right\vert +\left\vert V_{1}^{\prime }\right\vert \leq 2M\eta +2M\leq
3M,\left\Vert v_{1}^{\prime \prime }\right\Vert _{L^{2}\left( 0,1\right)
}\leq 3M.  \label{5.23}
\end{equation}%
Next, by (\ref{5.1}), (\ref{5.12}) and (\ref{5.22})%
\begin{equation}
\left\Vert \widetilde{v}_{1}^{\prime }\right\Vert _{C\left[ 0,1\right] }\leq
\eta \left\Vert \widetilde{q}_{1}^{\prime }\right\Vert _{C\left[ 0,1\right]
}+\left\Vert \widetilde{V}_{1}^{\prime }\right\Vert _{C\left[ 0,1\right]
}\leq M_{1}M^{5}\eta +M\eta \leq 2M_{1}M^{5}\eta .  \label{5.24}
\end{equation}%
We now estimate $\left\Vert \widetilde{\beta }_{1}\right\Vert _{L^{2}\left(
0,1\right) }.$ Subtracting (\ref{5.2}) from (\ref{4.40}) and using (\ref%
{4.22}), (\ref{5.51}), (\ref{5.6}), (\ref{5.12}), (\ref{5.22}) and (\ref%
{5.24}), we obtain 
\begin{eqnarray}
\left\Vert \widetilde{\beta }_{1}\right\Vert _{L^{2}\left( 0,1\right) }
&\leq &\eta \left\Vert \widetilde{q}_{1}\right\Vert _{H^{2}\left( 0,1\right)
}+\left\Vert \widetilde{V}_{1}^{\prime \prime }\right\Vert _{L^{2}\left(
0,1\right) }+\overline{k}^{2}\left\Vert \widetilde{v}_{1}^{\prime
}\right\Vert _{L^{2}\left( 0,1\right) }\left\Vert v_{1}^{\prime }+\left(
v_{1}^{\ast }\right) ^{\prime }\right\Vert _{C\left[ 0,1\right] }
\label{5.25} \\
&\leq &M_{1}M^{5}\eta +M\eta +2M_{1}M^{5}\overline{k}^{2}\left(
3M+M^{2}\right) \eta \leq M_{1}M^{8}\eta .  \notag
\end{eqnarray}

We now estimate functions $\widetilde{w}_{1}\left( x,\overline{k}\right)
=\left( w_{1}-w^{\ast }\right) \left( x,\overline{k}\right) $ and $%
w_{1}\left( x,\overline{k}\right) .$ Let 
\begin{equation*}
A_{n}\left( y\right) =y^{\prime \prime }-2\overline{k}y^{\prime }+\overline{k%
}^{2}\beta _{n}\left( x\right) y.
\end{equation*}
Recall that we find the function $w_{n}$ via solving the problem (\ref{4.3}%
), (\ref{4.4}) with $k=\overline{k}$ using the QRM. Also, it follows from (%
\ref{2.38}) and (\ref{4.11}) that $\left\vert g_{0}\left( \overline{k}%
\right) -g_{0}^{\ast }\left( \overline{k}\right) \right\vert \leq \eta $ and 
$\left\vert g_{1}\left( \overline{k}\right) -g_{1}^{\ast }\left( \overline{k}%
\right) \right\vert \leq M\eta .$ Hence, we obtain similarly with (\ref{5.15}%
) and (\ref{5.16}) 
\begin{eqnarray}
\left( A_{1}\widetilde{w}_{1},A_{1}z\right) +\eta ^{2}\left[ \widetilde{w}%
_{1},z\right] &=&-\left( \overline{k}^{2}\widetilde{\beta }_{1}w^{\ast
},A_{1}z\right) -\eta ^{2}\left[ w^{\ast },z\right] ,\forall z\in W_{0},
\label{5.26} \\
\left\vert \widetilde{w}_{1}\left( 0\right) \right\vert &\leq &\eta
,\left\vert \widetilde{w}_{1}^{\prime }\left( 0\right) \right\vert \leq
M\eta ,\widetilde{w}_{1}^{\prime }\left( 1\right) =0.  \notag
\end{eqnarray}%
The function $\widetilde{w}_{n}=w_{n}-w^{\ast }$ is the solution of a QRM
problem, which is completely similar with (\ref{5.26}). Since by (\ref%
{4.2323}) functions $\left\vert \beta _{n}\right\vert $ are uniformly
bounded for all $n$, $\left\vert \beta _{n}\right\vert \leq \max \left(
\left\vert c_{0}-1\right\vert ,\left\vert c_{1}-1\right\vert \right) ,$ then
there exists an analog of the constant $C_{2}$ of Theorem 3.2, which
estimates functions $\widetilde{w}_{n}$ for all $n$ as solutions of analogs
of problems (\ref{5.26}). Hence, we can assume that this constant equals $%
M_{1}$. Using Theorem 3.2, (\ref{5.25}) and (\ref{5.26}), we obtain%
\begin{equation}
\left\Vert \widetilde{w}_{1}\right\Vert _{C^{1}\left[ 0,1\right]
},\left\Vert \widetilde{w}_{1}\right\Vert _{H^{2}\left( 0,1\right) }\leq
M_{1}M^{10}\eta .  \label{5.27}
\end{equation}%
Hence, using (\ref{5.7}), (\ref{5.27}) and $w_{1}=\widetilde{w}_{1}+w^{\ast
},$ we obtain%
\begin{equation*}
\left\Vert w_{1}\right\Vert _{C^{1}\left[ 0,1\right] },\left\Vert
w_{1}\right\Vert _{H^{2}\left( 0,1\right) }\leq 2M.
\end{equation*}

The final step of the proof for the case $n=1$ is to estimate derivatives of
the second tail, i.e. functions $\widetilde{V}_{2}^{\prime },\widetilde{V}%
_{2}^{\prime \prime },V_{2}^{\prime },V_{2}^{\prime \prime }.$ By (\ref{4.17}%
)%
\begin{equation}
\widetilde{V}_{2}^{\prime }=\left( \frac{w_{1}^{\prime }}{\overline{k}%
^{2}w_{1}}-\frac{\left( w^{\ast }\right) ^{\prime }}{\overline{k}^{2}w^{\ast
}}\right) \left( x,\overline{k}\right) =\left( \frac{\widetilde{w}%
_{1}^{\prime }w^{\ast }-\widetilde{w}_{1}\left( w^{\ast }\right) ^{\prime }}{%
\overline{k}^{2}w_{1}w^{\ast }}\right) \left( x,\overline{k}\right) .
\label{5.29}
\end{equation}%
Estimate the denominator in (\ref{5.29}). Since $\left\vert O\left( 1/%
\overline{k}\right) \right\vert <1/2$ in (\ref{4.0}) for $c=c^{\ast },$ then
(\ref{2.1}) and (\ref{4.0}) imply that $\left\vert w^{\ast }\left( x,%
\overline{k}\right) \right\vert \geq c_{1}^{-1/4}/2$ for $x\in \left[ 0,1%
\right] .$ Hence, using (\ref{5.7}) and (\ref{5.27}), we obtain 
\begin{equation*}
\left\vert w_{1}\left( x,\overline{k}\right) \right\vert =\left\vert w^{\ast
}+\widetilde{w}_{1}\right\vert \left( x,\overline{k}\right) \geq \frac{%
c_{1}^{-1/4}}{2}-\left\vert \widetilde{w}_{1}\left( x,\overline{k}\right)
\right\vert \geq \frac{c_{1}^{-1/4}}{2}-\frac{c_{1}^{-1/4}}{4}=\frac{%
c_{1}^{-1/4}}{4}.
\end{equation*}%
Hence,%
\begin{equation}
\frac{1}{\overline{k}^{2}\left\vert w_{1}w^{\ast }\right\vert }\leq \frac{8%
\sqrt{c_{1}}}{\overline{k}^{2}}.  \label{5.30}
\end{equation}%
We now estimate from the above the modulus of the nominator in each of two
formulas of (\ref{5.29}). Using (\ref{5.501}), (\ref{5.6}) and (\ref{5.27})-(%
\ref{5.30}), we obtain for $x\in \left[ 0,1\right] $%
\begin{equation}
\left\vert \widetilde{V}_{2}^{\prime }\right\vert =\left\vert \frac{%
\widetilde{w}_{1}^{\prime }w^{\ast }-\widetilde{w}_{1}\left( w^{\ast
}\right) ^{\prime }}{\overline{k}^{2}w_{1}w^{\ast }}\right\vert \left( x,%
\overline{k}\right) \leq 16\sqrt{c_{1}}M_{1}M^{11}\eta \leq M_{1}M^{12}\eta .
\label{5.31}
\end{equation}%
Next, 
\begin{equation*}
\widetilde{V}_{2}^{\prime \prime }=\left( \frac{\widetilde{w}_{1}^{\prime
\prime }w^{\ast }-\widetilde{w}_{1}\left( w^{\ast }\right) ^{\prime \prime }%
}{\overline{k}^{2}w_{1}w^{\ast }}\right) \left( x,\overline{k}\right)
-\left( \frac{w_{1}^{\prime }}{\overline{k}^{2}w_{1}}-\frac{\left( w^{\ast
}\right) ^{\prime }}{\overline{k}^{2}w^{\ast }}\right) \left( \frac{%
w_{1}^{\prime }}{w_{1}}+\frac{\left( w^{\ast }\right) ^{\prime }}{w^{\ast }}%
\right) \left( x,\overline{k}\right) .
\end{equation*}%
Hence, we obtain similarly with (\ref{5.31}) 
\begin{equation}
\left\Vert \widetilde{V}_{2}^{\prime \prime }\right\Vert _{L^{2}\left(
0,1\right) }\leq M_{1}M^{14}\eta .  \label{5.32}
\end{equation}%
It can be easily derived from (\ref{5.7}), (\ref{5.31}) and (\ref{5.32})
that $\left\Vert V_{2}^{\prime }\right\Vert _{C\left[ 0,1\right] }\leq 2M$
and $\left\Vert V_{2}^{\prime \prime }\right\Vert _{L^{2}\left( 0,1\right)
}\leq 2M.$

Thus, in summary (\ref{5.22}), (\ref{5.25}), (\ref{5.31}) and (\ref{5.32})
imply that 
\begin{equation}
\left\Vert \widetilde{q}_{1}\right\Vert _{H^{2}\left( 0,1\right)
},\left\Vert \widetilde{q}_{1}\right\Vert _{C^{1}\left[ 0,1\right]
},\left\Vert \widetilde{\beta }_{1}\right\Vert _{L^{2}\left( 0,1\right)
},\left\Vert \widetilde{V}_{2}^{\prime }\right\Vert _{C\left[ 0,1\right]
},\left\Vert \widetilde{V}_{2}^{\prime \prime }\right\Vert _{L^{2}\left(
0,1\right) }\leq M_{1}M^{14}\eta .  \label{5.33}
\end{equation}%
In other words, estimates (\ref{5.8})-(\ref{5.10}) are valid for $n=1$.
Assume that they are valid for $n-1$ where $n\geq 2$. Denote $K_{n-1}=\left(
M_{1}M^{14}\right) ^{n-1}.$ Then, similarly with the above, one will obtain
estimates (\ref{5.33}) where \textquotedblleft 1" in first three terms is
replaced with $n$, \textquotedblleft 2" in the fourth and fifth terms is
replaced with $n+1$ and the right hand side is $M_{1}M^{14}K_{n-1}\eta
=\left( M_{1}M^{14}\right) ^{n}\eta .$ $\square $

\section{Numerical results}

\label{sec:6}

In all our computation, $x_{0}=-1$ and $k\in \left[ 0.5,1.5\right] $. We
have observed in our computationally simulated data as well as in
experimental data that the function $\left\vert u(x,k)\right\vert $ becomes
very small for $k>2.$ On the other hand, the largest values of $\left\vert
u(x,k)\right\vert $ were observed in some points of the interval $k\in \left[
0.5,1.5\right] $. Thus, we assign in all computations $\overline{k}=1.5,%
\underline{k}=0.5.$ We have used $h= 0.02.$

%

Actually in all our computations we go along the interval $k\in \left[
0.5,1.5\right] $ several times. More precisely, let $\beta ^{\left( 1\right)
}\left( x\right) $ be the result obtained in Step \ref{step 3} of the above
globally convergent algorithm. Set $c^{\left( 1\right) }\left( x\right)
=1+\beta ^{\left( 1\right) }\left( x\right) $. Next, solve the problem the
problem (\ref{4.3}), (\ref{4.4}) with $k:=\overline{k}$ and $\beta \left(
x\right) :=\beta ^{\left( 1\right) }\left( x\right) $. Let the function $%
w^{\left( 1\right) }\left( x,\overline{k}\right) $ be its solution. Then we
define derivatives of the new tail function $V_{0}^{\left( 1\right) }$ as in
(\ref{4.17}) where $w_{n,j}(x,\overline{k})$ is replaced with $w^{\left(
1\right) }\left( x,\overline{k}\right) .$ Next, we go to Step 2 and repeat.
The process is repeated $K=50$ times in our computational program. We choose 
$m^{0}$ such that%
\begin{equation*}
\frac{\left\Vert c_{m^{0}}-c_{m^{0}-1}\right\Vert _{L_{2}\left( 0,1\right) }%
}{\left\Vert c_{m^{0}}\right\Vert _{L_{2}\left( 0,1\right) }}=\min_{2\leq
m\leq K}\left\{ \frac{\left\Vert c_{m}-c_{m-1}\right\Vert _{L_{2}\left(
0,1\right) }}{\left\Vert c\right\Vert _{L_{2}\left( 0,1\right) }}\right\} .
\end{equation*}%
\textrm{\ Our final solution of the inverse problem is }$c\left( x\right)
=c_{m^{0}}\left( x\right) .$ It is also worth mentioning that in Step \ref%
{step 2b3}, after calculating $\beta _{n,j}$, we replace it by 
\begin{equation*}
\beta _{n,j}(x):=\frac{1}{l_{x}} \int_{U_{x}\cap (0,1)}\beta _{n,j}(y)dy
\end{equation*}%
where $U_{x}$ is a small neighborhood of $x$, $x\in (0,1)$ and $l_{x}$ is
the length of the interval $U_{x}\cap (0,1).$ We also use the truncation
technique to improve the accuracy of the reconstruction function $\beta_{n,
j}$ (see \eqref{4.2323}).

%
%

\subsection{Computationally simulated data}

\label{sec:6.1}

In this section, we show the numerical reconstruction of the spatially
distributed dielectric constant from computationally simulated data. Let the
function $c(x)$ has the form 
\begin{equation*}
c(x)=\left\{ 
\begin{array}{ll}
c_{\mathrm{target}} & \mbox{in }(1/4,1/3), \\ 
1 & \mbox{otherwise}%
\end{array}%
\right.
\end{equation*}%
Let the function $u(x,k)$ be the solution of problem \eqref{2.4}, \eqref{2.6}%
. As mentioned in the proof of Theorem \ref{thm unique forward} (see %
\eqref{2.13}), the function $u(x,k)$ satisfies the Lippman-Schwinger
equation, 
\begin{equation*}
u(x,k)=\frac{\exp (-ik|x-x_{0}|)}{2ik}+k^{2}\int_{0}^{1}\frac{\exp
(-ik|x-\xi |)}{2ik}(c(\xi )-1)u(\xi ,x_{0},k)d\xi .
\end{equation*}%
This equation can be approximated as a linear system. We have solved that
system numerically to computationally simulate the data.



Our numerical results are displayed in Figures 1. In the top row $c_{\mathrm{%
target}}=4$ and in the bottom row $c_{\mathrm{target}}=7.$ In each figure,
we show the true function $c(x)$, the data obtained by solving the forward
problem with that true $c(x)$ and the solution of the CIP. In both cases we
had two values of $m^{0}:m^{0}=24$ and $m^{0}=25.$ These figures confirm
that both the target/background contrast and the position of the target are
computed with small errors.

\begin{figure}[h!]
\hfill%
\subfloat[\labela]{
		\includegraphics[width=\width, height = \height ]{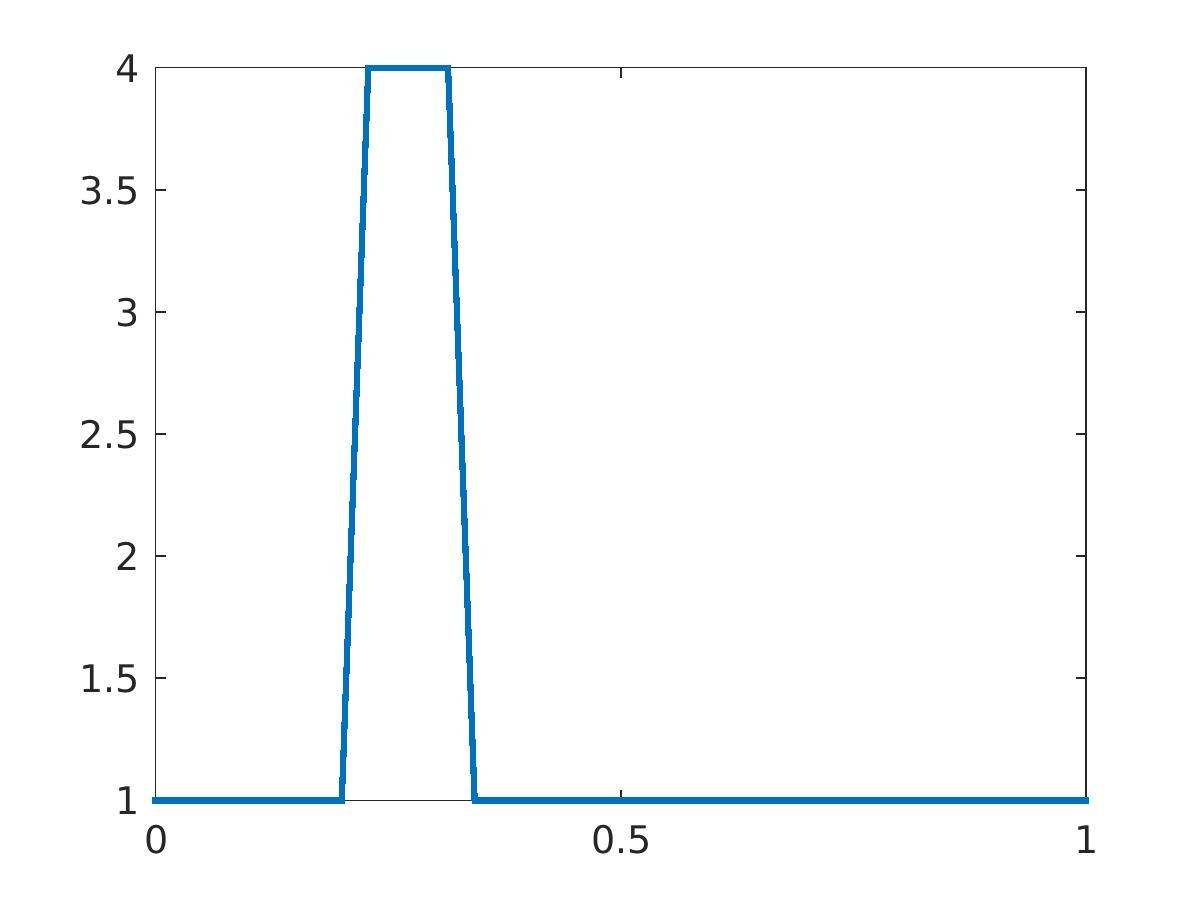} 
} \hfill%
\subfloat[\labelb]{
		\includegraphics[width=\width, height = \height ]{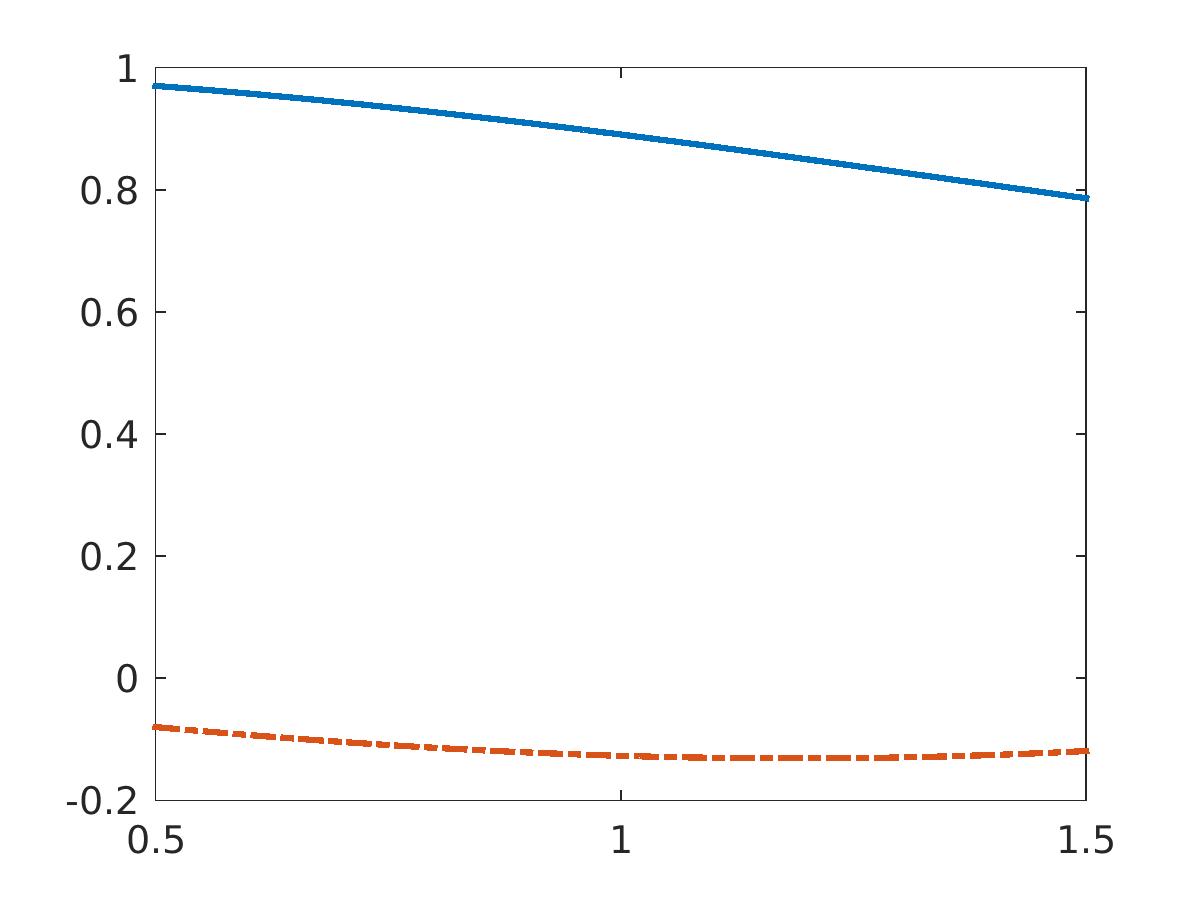} 
} \hfill%
\subfloat[\labelc]{
		\includegraphics[width=\width, height = \height ]{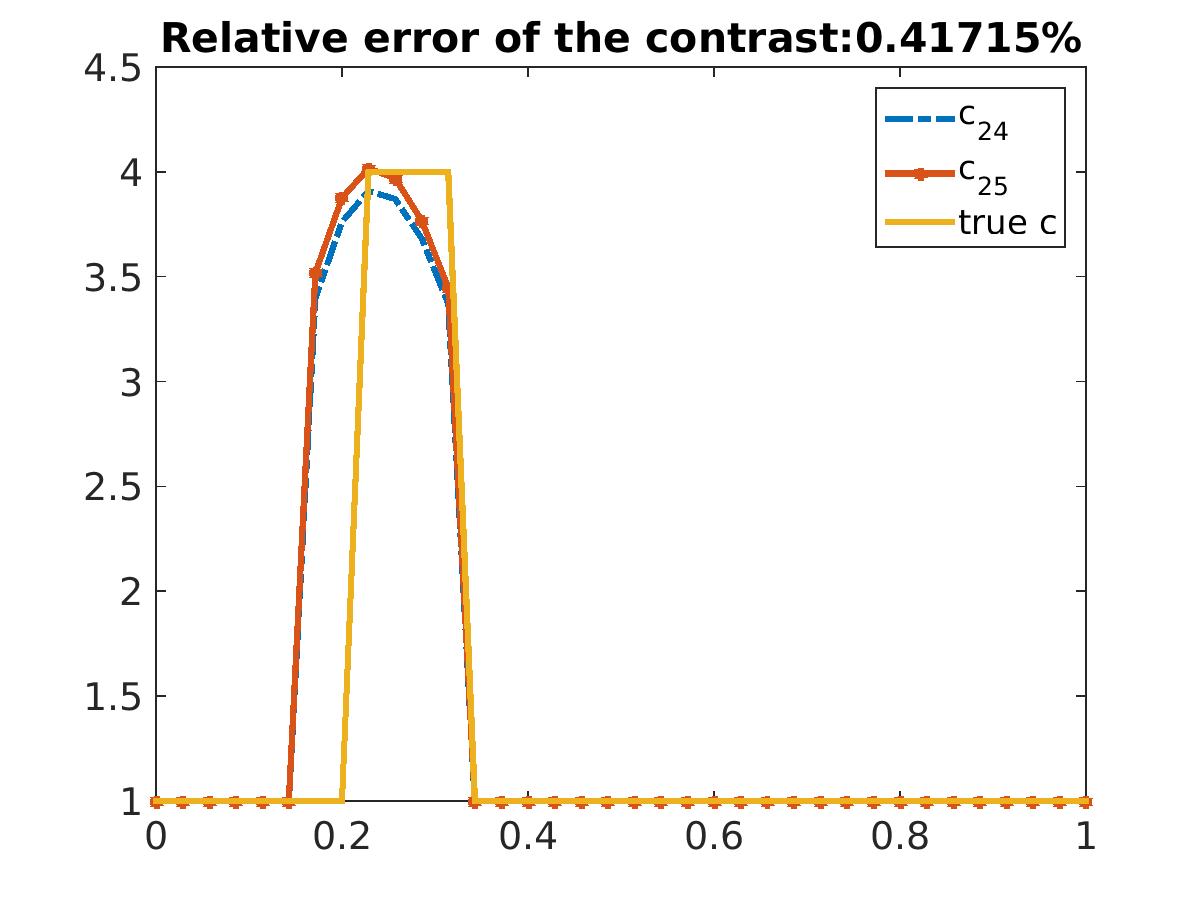} 
}
\par
\hfill%
\subfloat[\labela]{
		\includegraphics[width=\width, height = \height ]{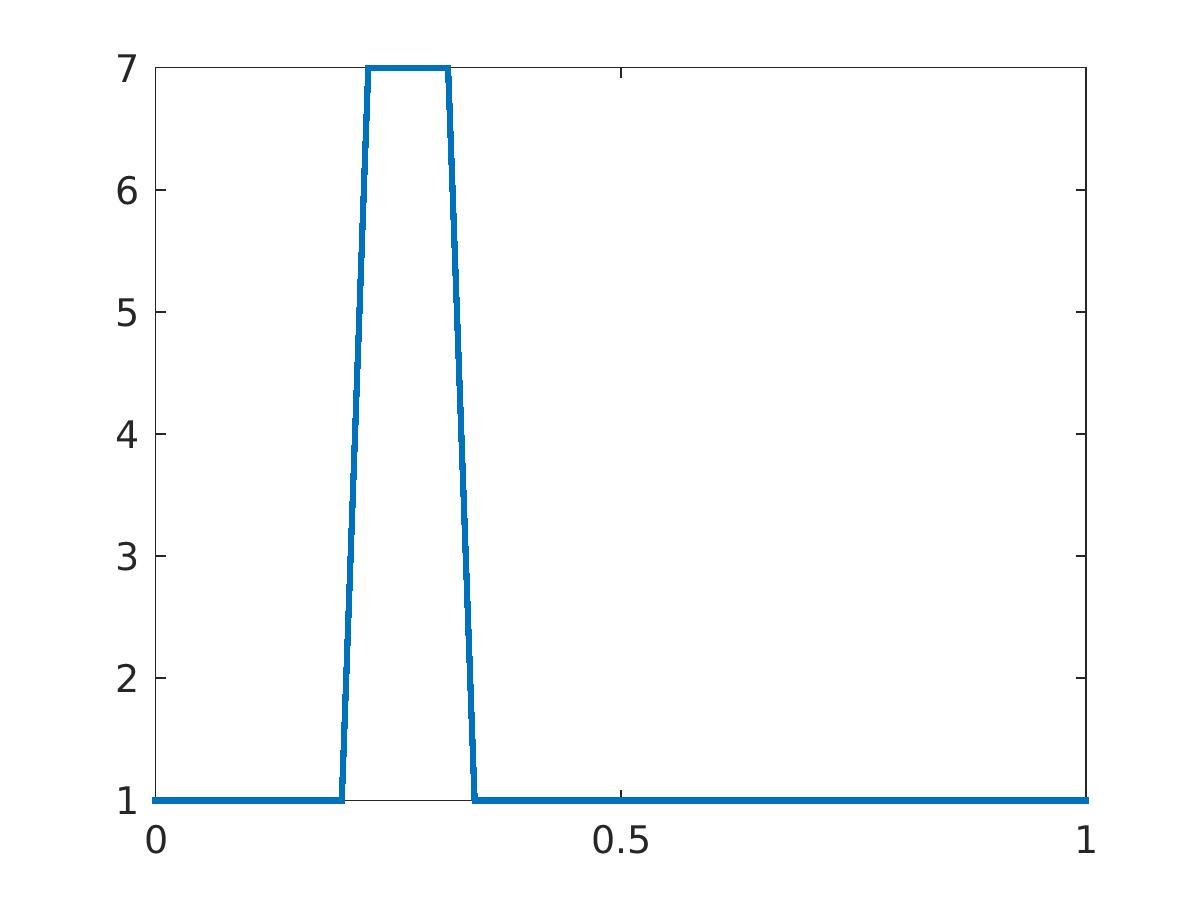} 
} \hfill%
\subfloat[\labelb]{
		\includegraphics[width=\width, height = \height ]{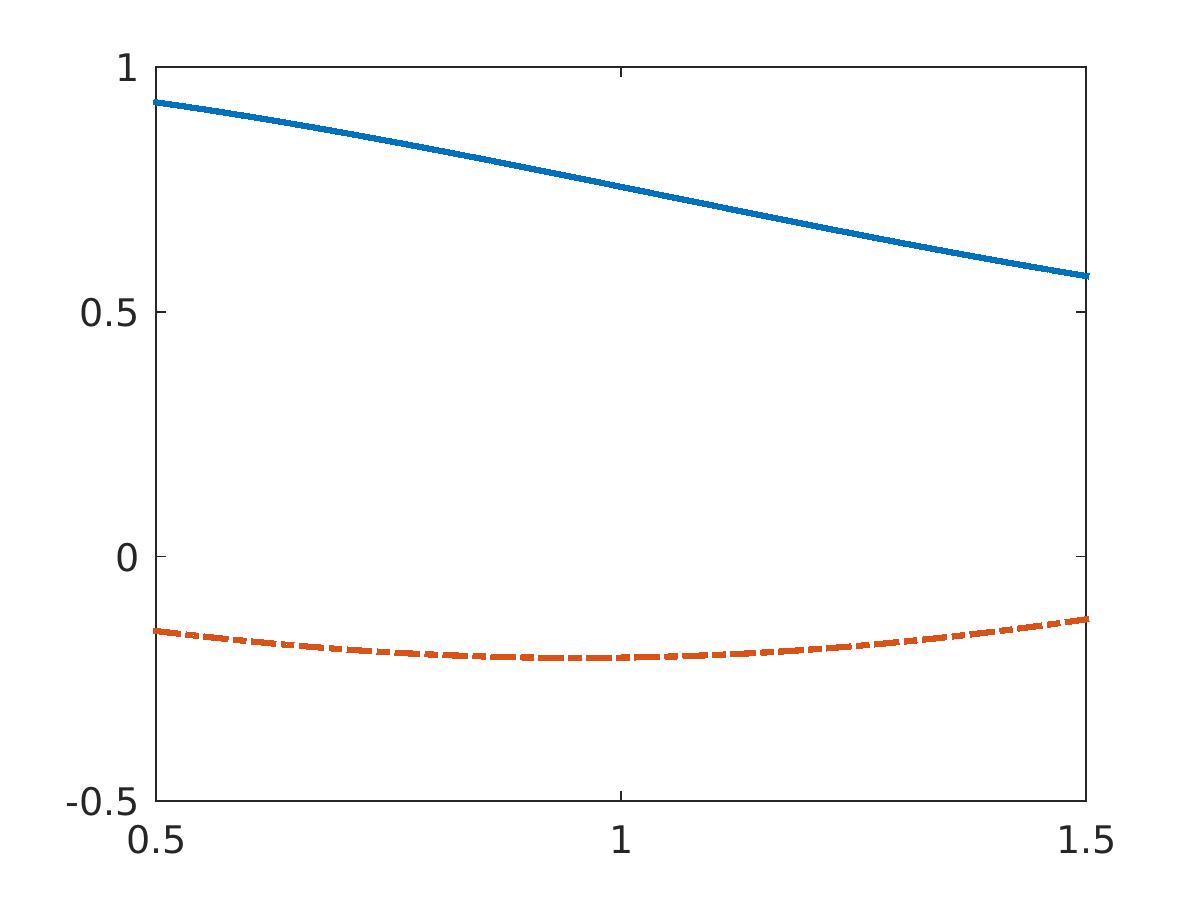} 
} \hfill%
\subfloat[\labelc]{
		\includegraphics[width=\width, height = \height ]{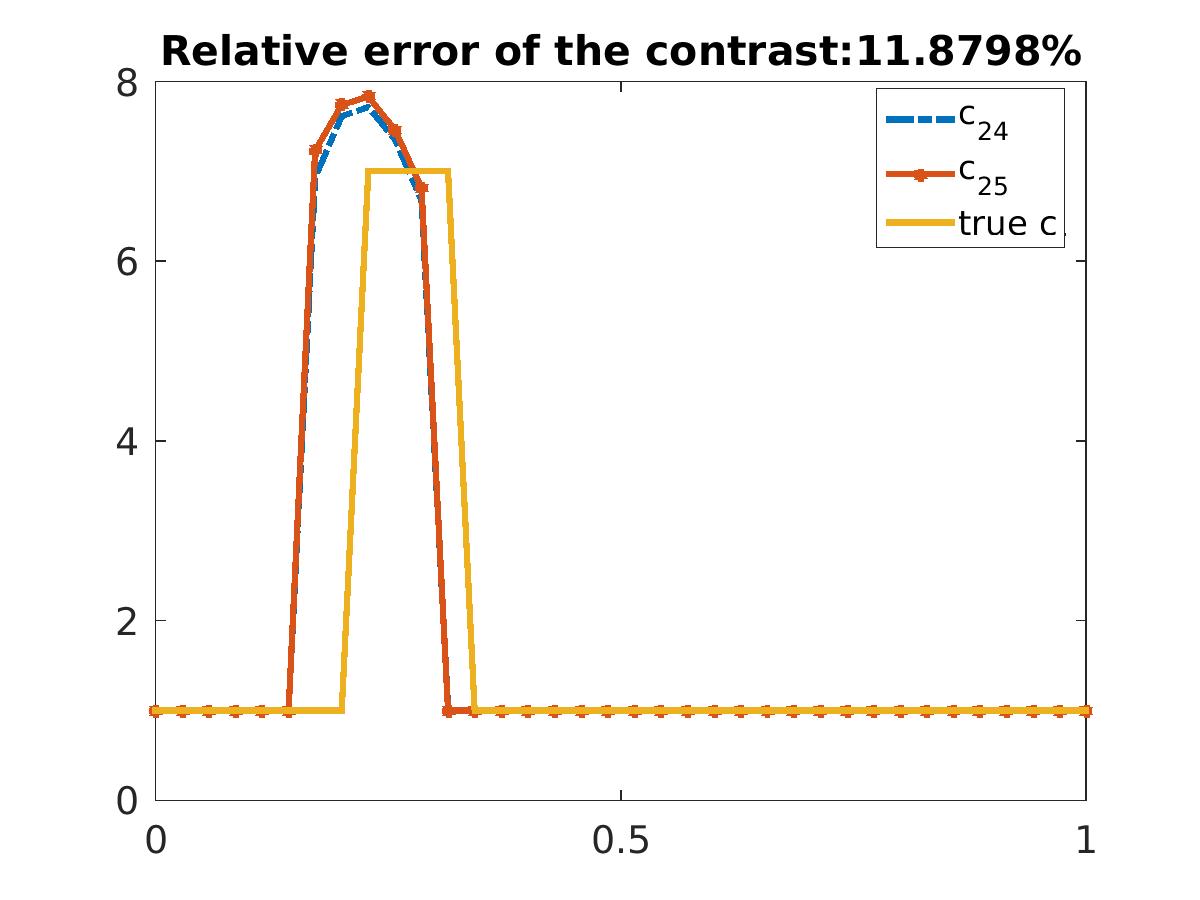} 
}
\caption{\textit{Numerical results when $c_{\mathrm{target}} $ is 4 (in row
1) and 7 (in row 2). The relative error in figures (c) and (f) is $%
\left|\|c_{25}\|_{L^{\infty}(0, 1)} - \|c^*\|_{L^{\infty}(0,
1)}\right|/\|c^*\|_{L^{\infty}(0, 1)}$.}}
\label{fig 2}
\end{figure}

\begin{remark}
\textrm{In our computer program, we use the linear algebra package, named as
Armadillo \cite{San} to solve linear systems. The software is very helpful
to speed up the program and to simplify the codes. }
\end{remark}

\subsection{Experimental data}

\label{sec:6.2}

The experimental data were collected by the Forward Looking Radar which was
built in the US Army Research Laboratory \cite{Ng}. The device consists of
two main parts. The first one, emitter, generates the time resolved electric
pulses. The emitter sends out only one component of the electric field. The
second part involves 16 detectors. These detectors collect the time resolved
backscattering electric signal (voltage) in the time domain. The same
component of the electric field is measured as the one which is generated by
the emitter: see our comment about this in the beginning of section 2. The
step size of time is 0.133 nanosecond. The backscattering data in the time
domain are collected when the distance between the radar and the target
varies from 20 to 8 meters. The average of these data with respect to both
the position of the radar and those 16 detectors is the data on which we
have tested our algorithm. To identify horizontal coordinates of the
position of the target, Ground Positioning System (GPS) is used. The error
in each of horizontal coordinates does not exceed a few centimeters. When
the target is under the ground, the GPS provides the distance between the
radar and a point on the ground located above the target. As to the depth of
a buried target, it is not of a significant interest, since horizontal
coordinates are known and it is also known that the depth does not exceed 10
centimeters. We refer to \cite{Ng} for more details about the data
collection process. Publications \cite{Karch,KSNF1,KSNF2} contain schematic
diagrams of the measurements.

Our interest is in computing maximal values of dielectric constants of
targets. In one target (plastic cylinder below) we compute the minimal value
of its dielectric constant, since its value was less than the dielectric
constant of the ground. For each target, the only information the
mathematical team (MVK,LHN) had, in addition to just a single experimental
curve, was whether it was located in air or below the ground.

We calculate $R(x),$ the relative spatial dielectric constant of the whole
structure including the background (air or the ground) and the target. More
precisely, 
\begin{equation}
R(x)=\left\{ 
\begin{array}{ll}
\displaystyle\frac{c_{\mathrm{target}}}{c_{\mathrm{bckgr}}} & x\in D, \\ 
1 & \mathrm{otherwise}%
\end{array}%
\right.  \label{6.1}
\end{equation}%
where $D$ is the subinterval of the interval $\left( 0,1\right) ,$ which is
occupied by the target. Here $c_{\text{target}}$ and $c_{\text{bckgr}}$ are
values of the function $c\left( x\right) $ in target and background
respectively. We assume that $c_{\text{bckgr}}=const.>0$ for each set of
experimental data. Hence, $c_{\mathrm{bckgr}}=1$ if the target is located in
air. The ground was dry sand. It is well known that the dielectric constant
of the dry sand varies between 3 and 5 \cite{Tabl}. Hence, in the case of
buried targets $c_{\mathrm{bckgr}}\in \lbrack 3,5]$.

In our mathematical model the time resolved electric signal $\widehat{u}%
(x,x_{0},t)$ collected by the detectors satisfies the equations \eqref{2.15}%
, \eqref{2.16} with $c(x)$ being replaced by $R(x),$ where the position $%
x_{0}$ of the source is actually unknown. The latter is one of the
difficulties of working with these data. Thus, we set in all our tests $%
x_{0}=-1,$ which is the same as in \cite{Karch,KSNF1,KSNF2}.

Let $R_{\text{comp}}(x)$ be the function $R(x)$ which we compute. Following 
\cite{KSNF1}\textrm{, we define the computed target/background contrast as}%
\begin{equation}
\widetilde{R}=\left\{ 
\begin{array}{c}
\max R\left( x\right) \text{ if }R\left( x\right) \geq 1,\forall x\in \left[
0,1\right] , \\ 
\min R\left( x\right) \text{ if }0<R\left( x\right) \leq 1,\forall x\in %
\left[ 0,1\right] .%
\end{array}%
\right.  \label{6.2}
\end{equation}%
\textrm{\ }Since the dielectric constant of air equals 1, then we have $%
R\left( x\right) \geq 1$ for targets located in air. As to the buried
targets, we have developed a procedure of the analysis of the original time
resolved data, which provides us with the information on which of two cases (%
\ref{6.2}) takes place. We refer to Case 1 and Case 2 on page 2944 of \cite%
{KSNF2} for this procedure. In addition, since we had a significant mismatch
of magnitudes of experimental and computationally simulated data, we have
multiplied, before computations, our experimental data by the calibration
number $10^{-7},$ see \cite{KSNF1,KSNF2} for details of our choice of this
number.

There is a significant discrepancy between computationally simulated and
experimental data, which was noticed in our earlier publications \cite%
{Karch,KSNF1,KSNF2}. This discrepancy is evident from, comparison of, e.g.
Figure 1b with Figure 2b and other similar ones. Therefore, to at least
somehow mitigate this discrepancy, we perform a data pre-processing
procedure. Besides of the Fourier transform of the time resolved data, we
multiply them by a calibration factor, truncate a certain part of the data
in the frequency domain and shift the data in the frequency domain, see
details below.

The function $u(0,x_{0},k),$ which we have studied in the previous sections,
is the Fourier transform of the time resolved data $\widehat{u}(0,x_{0},t)$.
The function $u(0,x_{0},k)$ is called \textquotedblleft the data in the
frequency domain". Observing that $|u(0,x_{0},k)|$ is small when $k$ belongs
to a certain interval, we do not analyze $u(0,x_{0},k)$ on that interval.
Rather, we only focus on such a frequency interval which contains the major
part of the information. Now, to keep the consistency with our study of
computationally simulated data, we always force the frequency interval to be 
$\left[ 0.5,1.5\right] $. To do so, we simply shift our data in the
frequency domain: compare Figures 2b and 2c, Figures 2f and 2g, Figures 3b
and 3c, Figures 3f and 3g and Figures 3j and 3k.

\begin{figure}[h!]
\hfill 
\subfloat[\labelaa]{
		\includegraphics[width=\idth, height = \eight ]{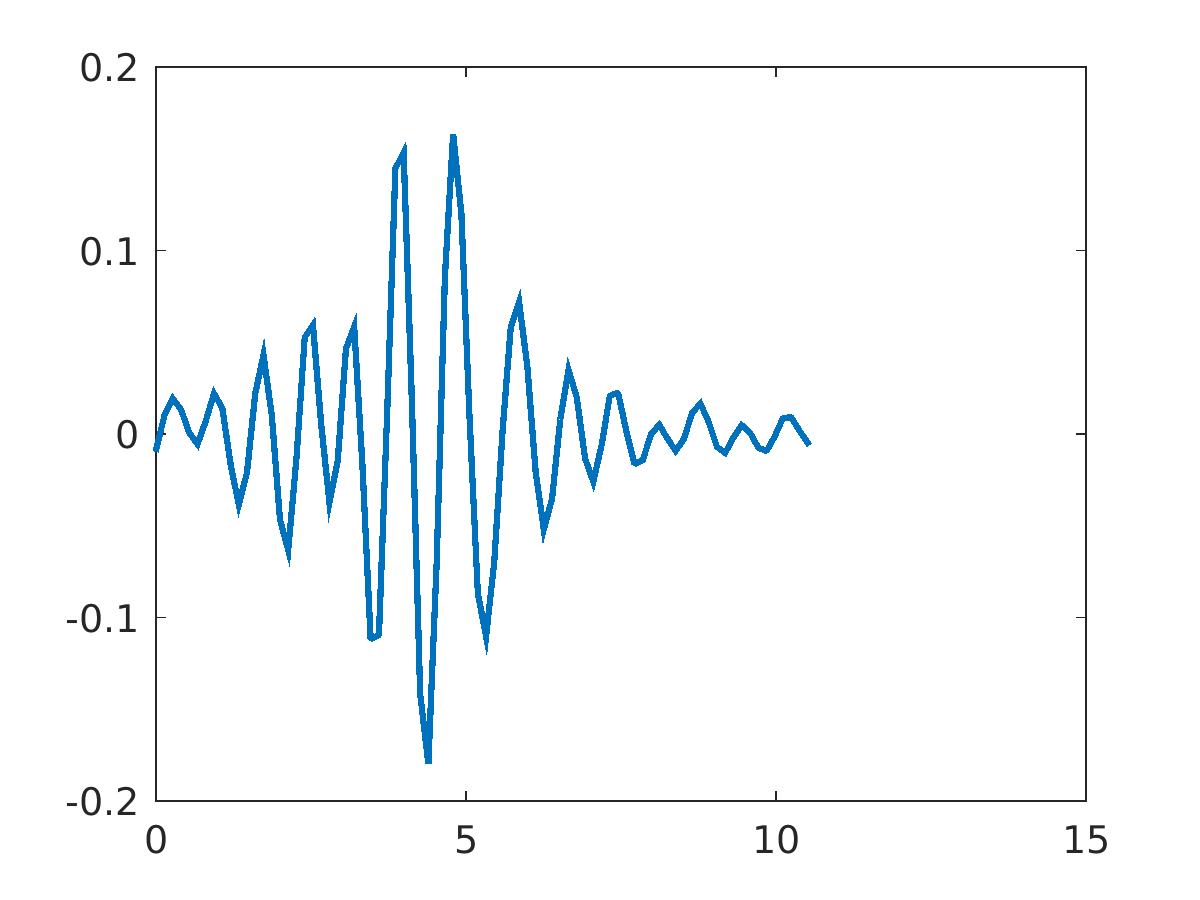} 
} \hfill 
\subfloat[\labelbb]{
		\includegraphics[width=\idth, height = \eight ]{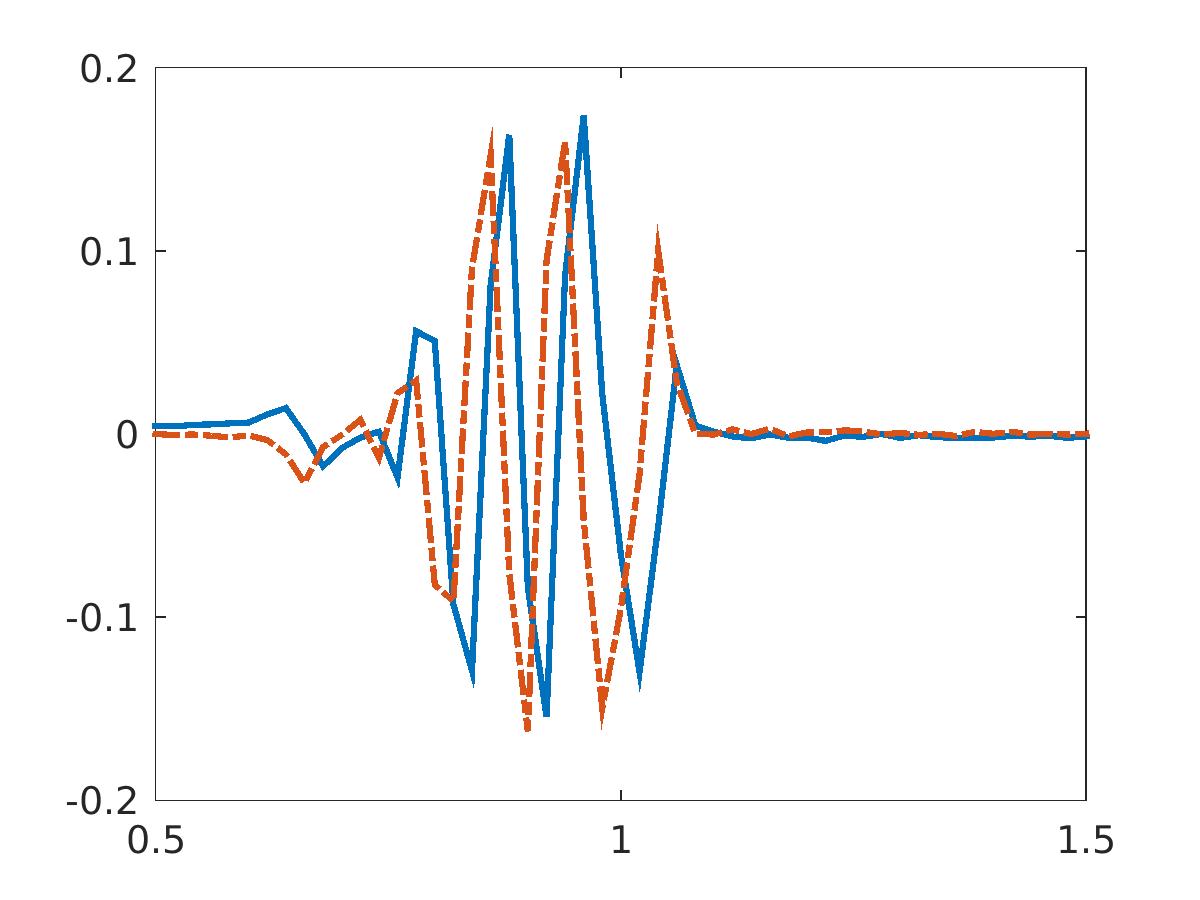} 
} \hfill 
\subfloat[\labelcc]{
		\includegraphics[width=\idth, height = \eight ]{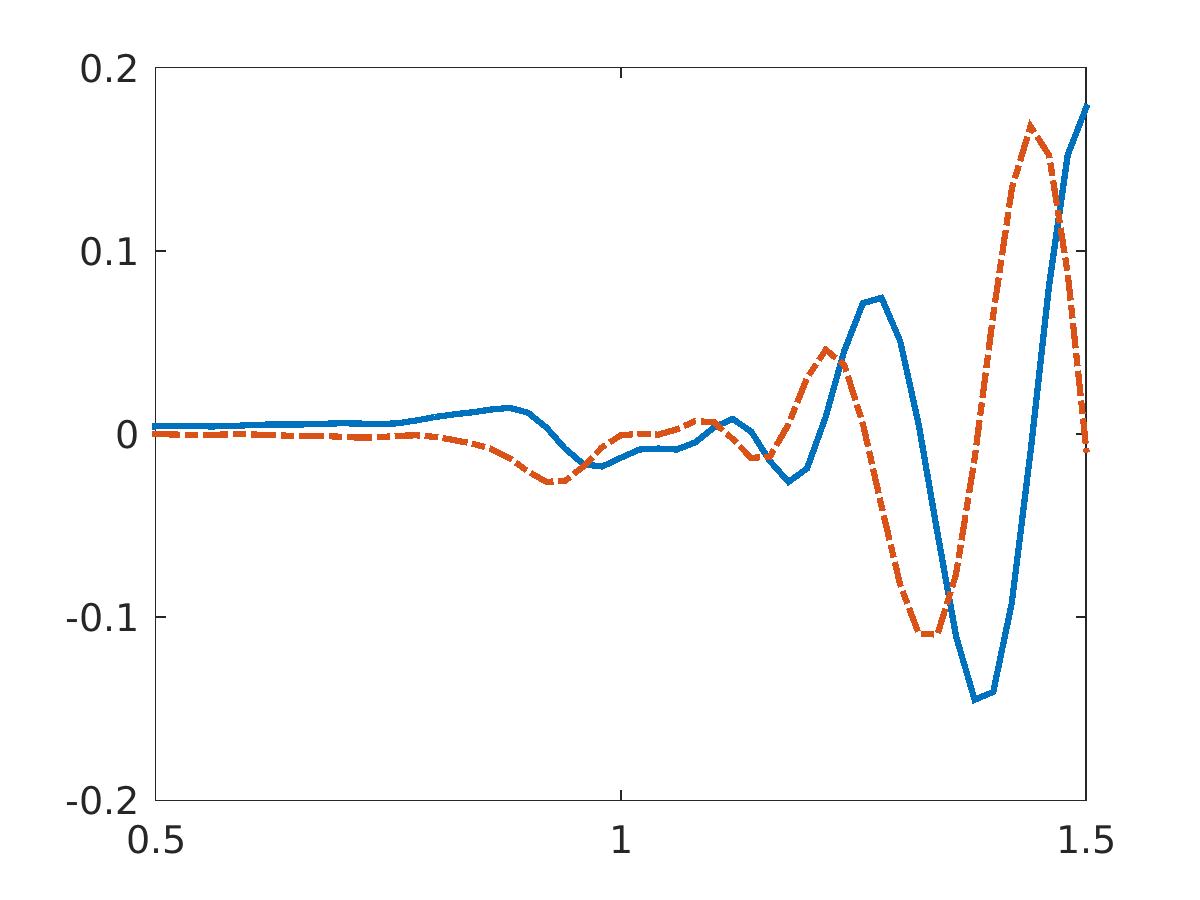} 
} \hfill 
\subfloat[\labeldd]{
		\includegraphics[width=\idth, height = \eight ]{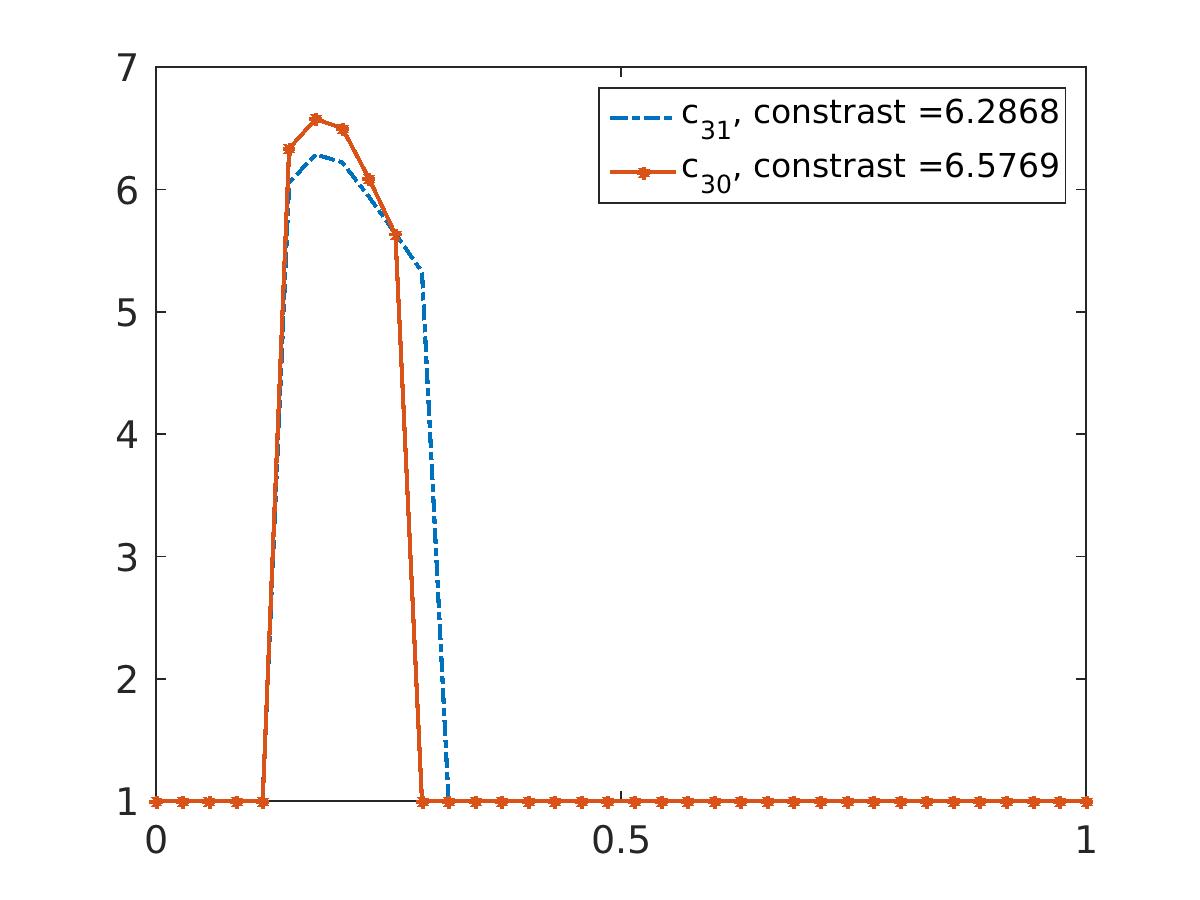}
}
\par
\hfill 
\subfloat[\labelaa]{
		\includegraphics[width=\idth, height = \eight ]{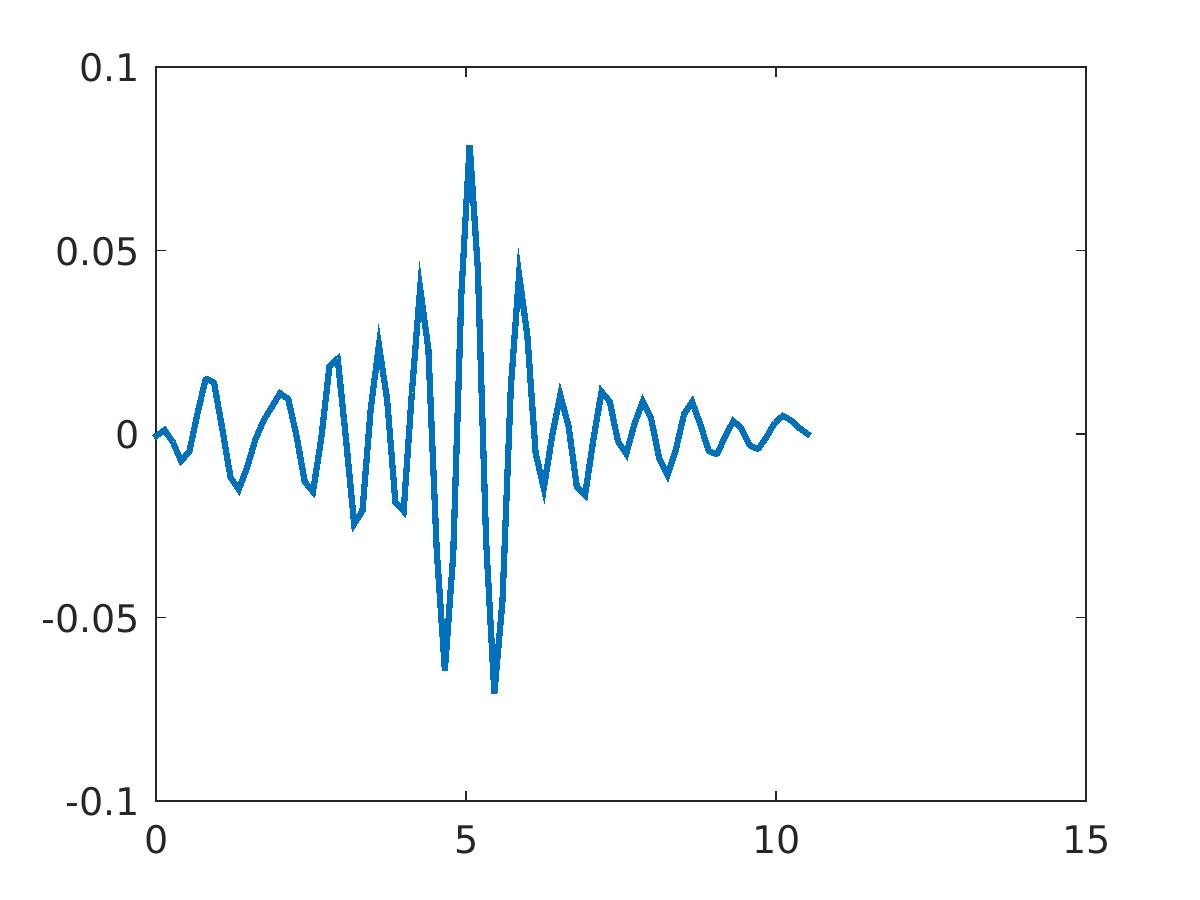} 
} \hfill%
\subfloat[\labelbb]{
		\includegraphics[width=\idth, height = \eight ]{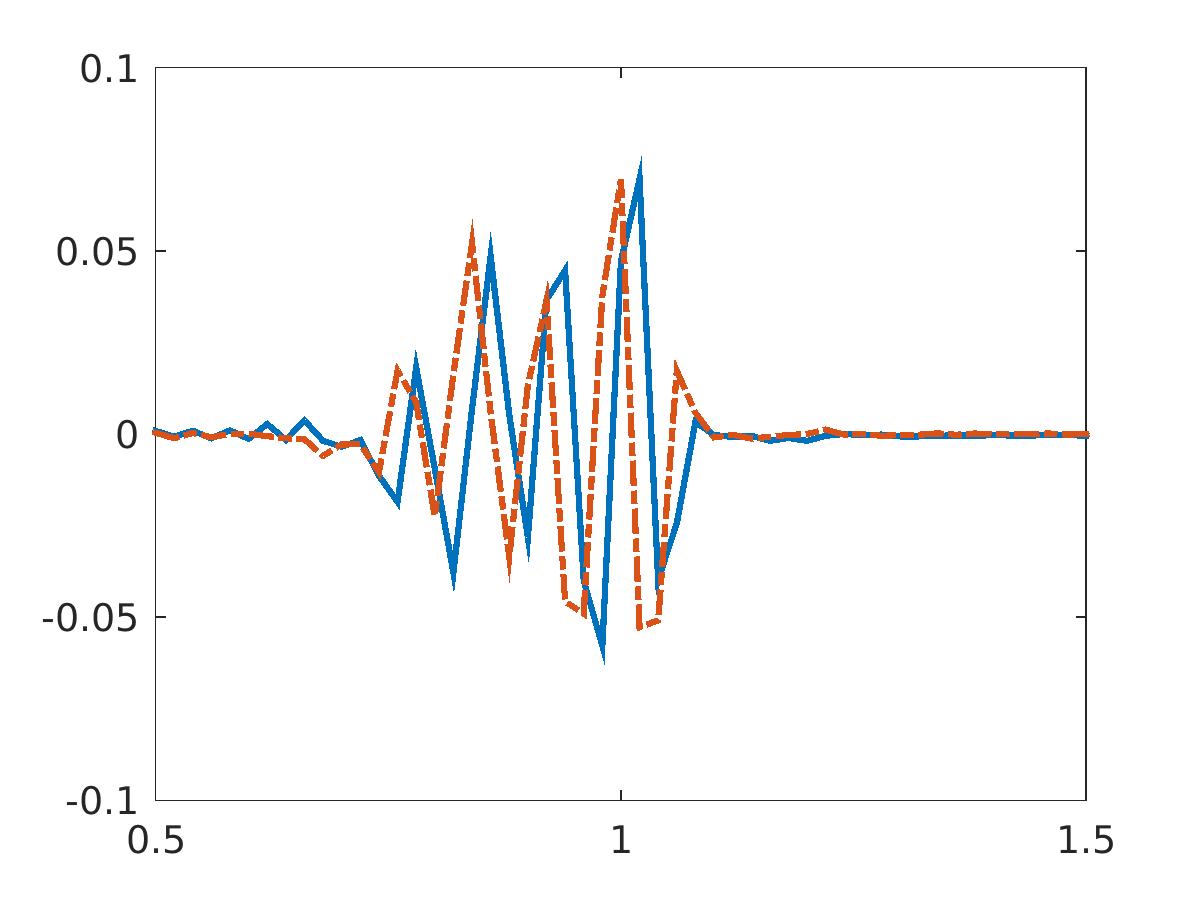} 
} \hfill%
\subfloat[\labelcc]{
		\includegraphics[width=\idth, height = \eight ]{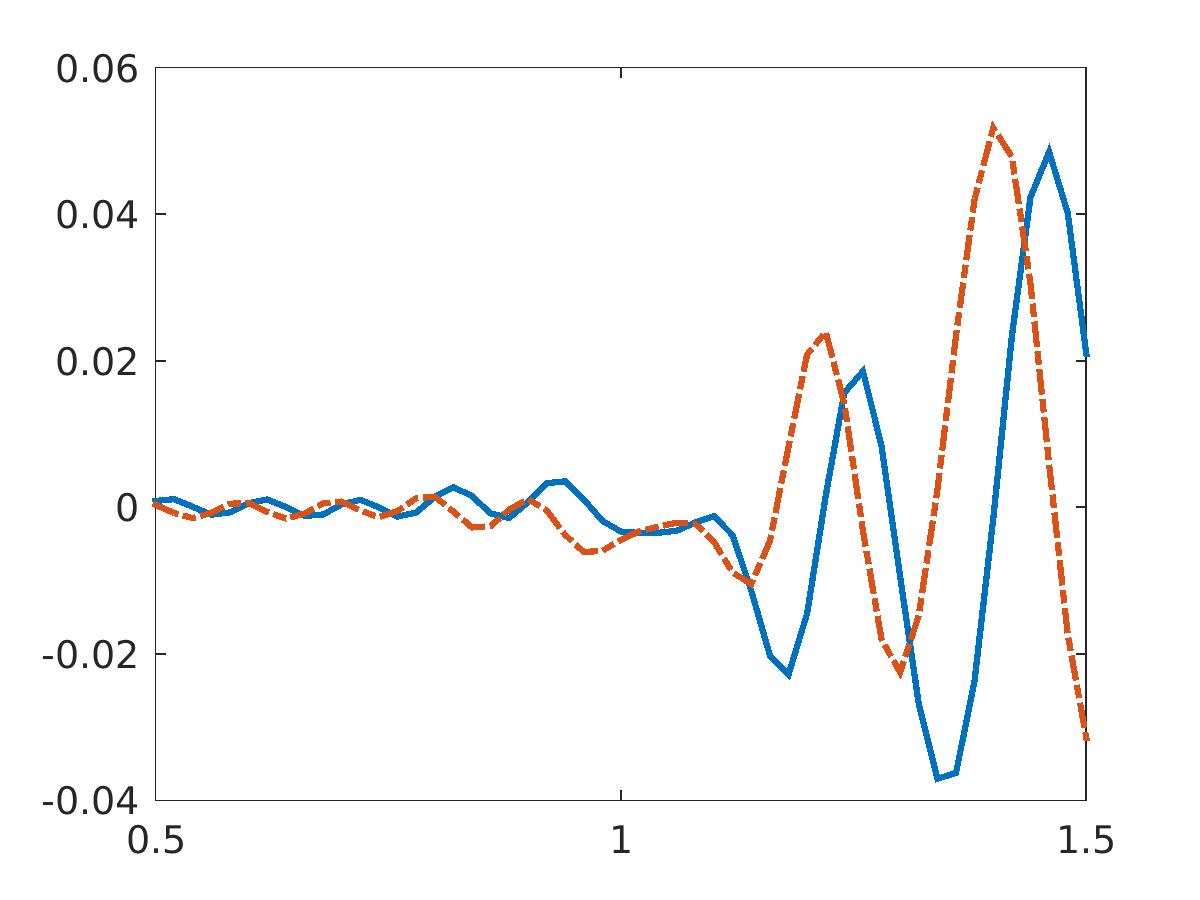} 
} \hfill%
\subfloat[\labeldd]{
		\includegraphics[width=\idth, height = \eight ]{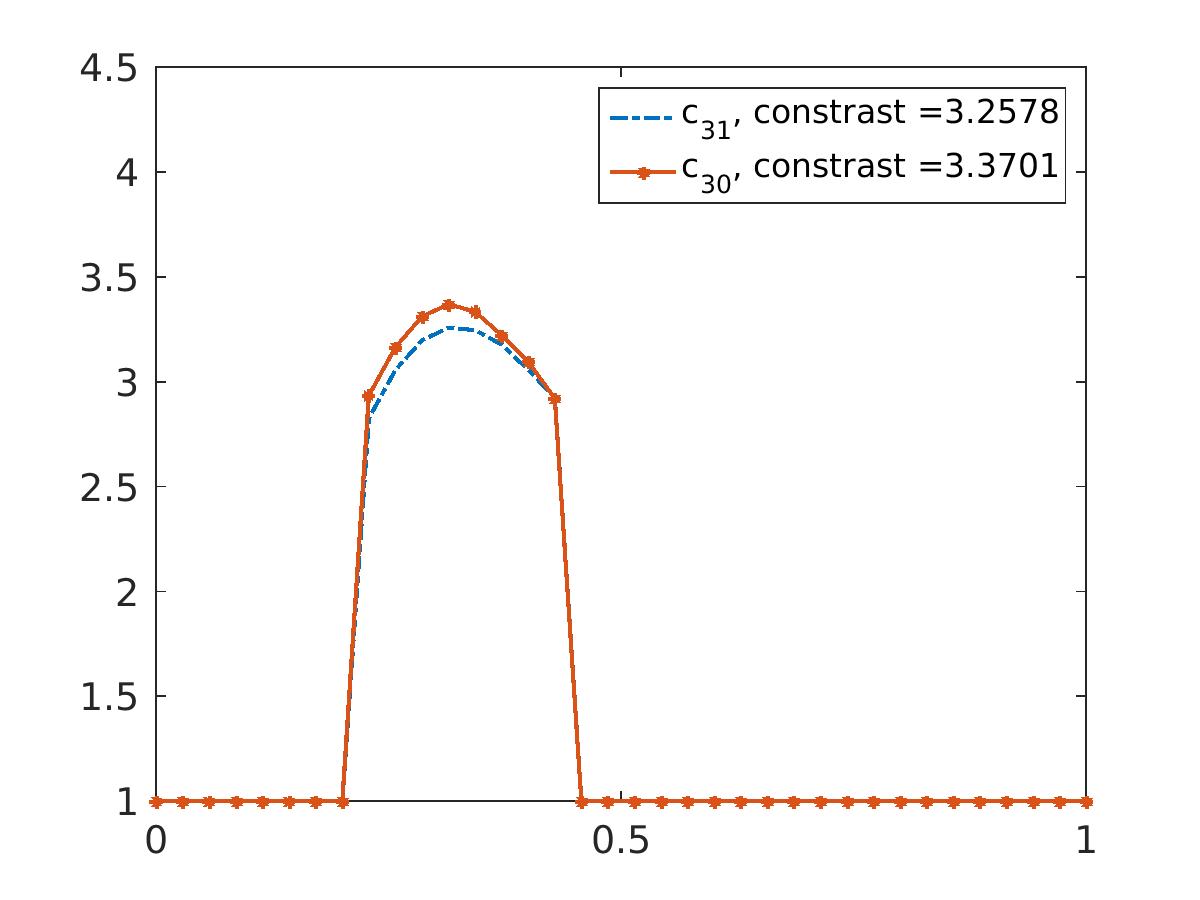}
}
\caption{\textit{The numerical test to evaluate the ``relative" dielectric
constant of bush (first row) and wood stake (second row) when they are put
in the air. Solid lines on b,c,f,g are real parts and dotted lines are
imaginary parts.} }
\label{Air}
\end{figure}

\begin{figure}[h!]
\hfill 
\subfloat[\labelaa]{
		\includegraphics[width=\idth, height = \eight ]{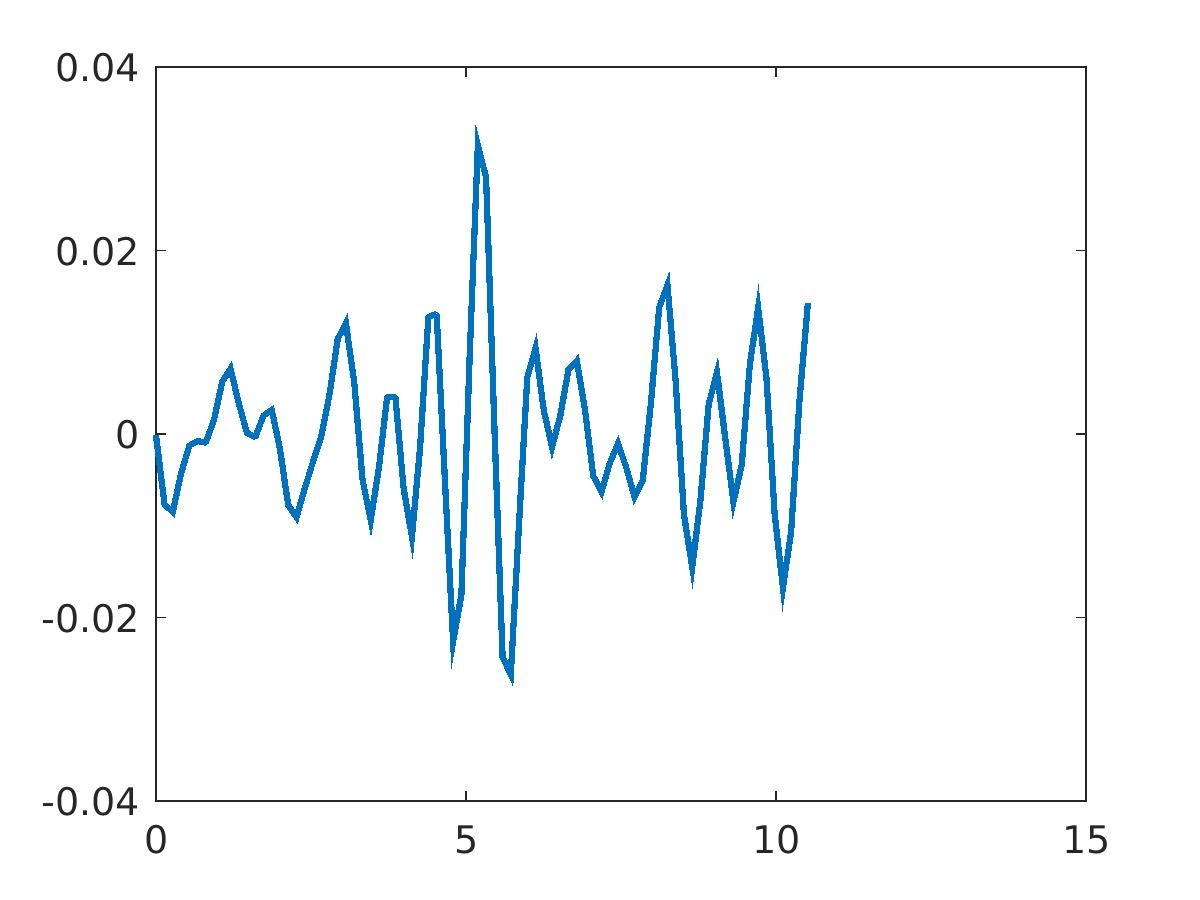} 
} \hfill%
\subfloat[\labelbb]{
		\includegraphics[width=\idth, height = \eight ]{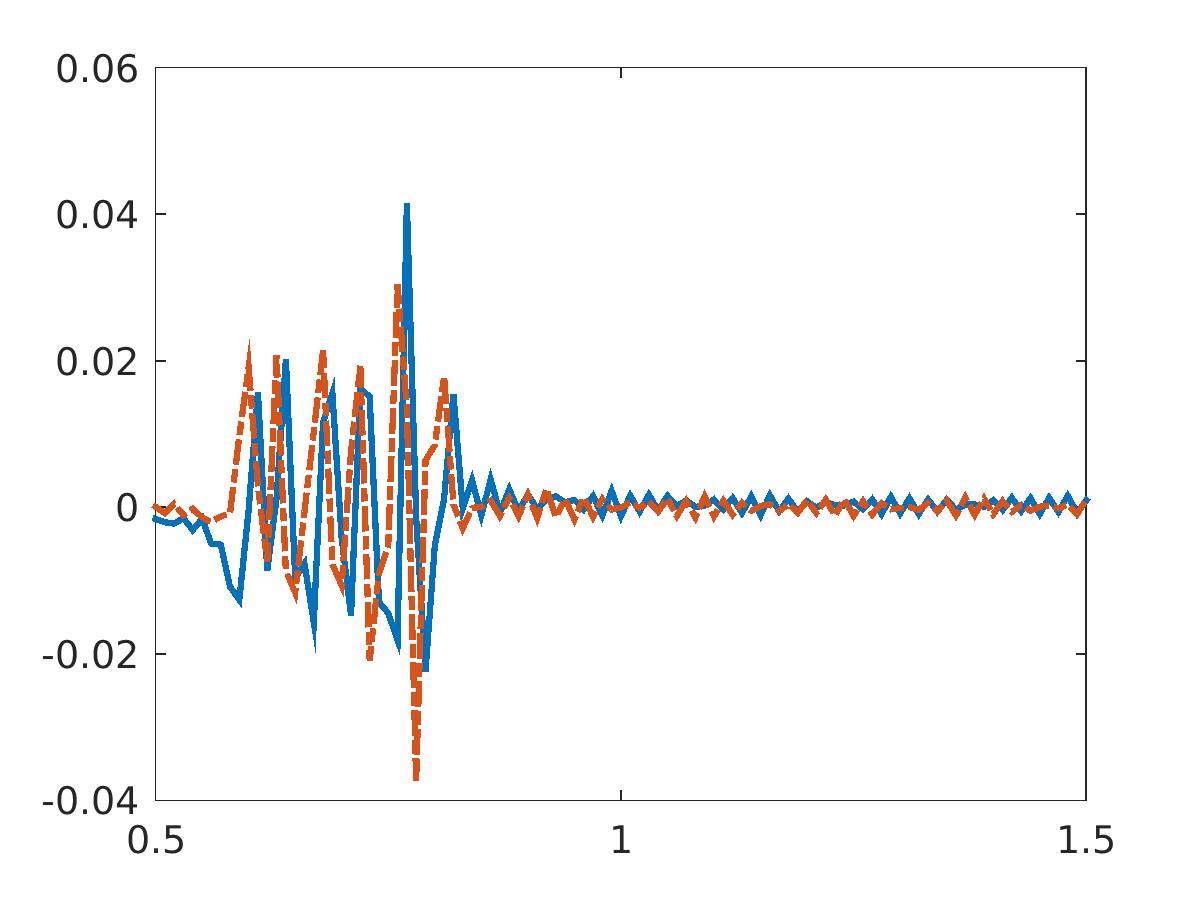} 
} \hfill%
\subfloat[\labelcc]{
		\includegraphics[width=\idth, height = \eight ]{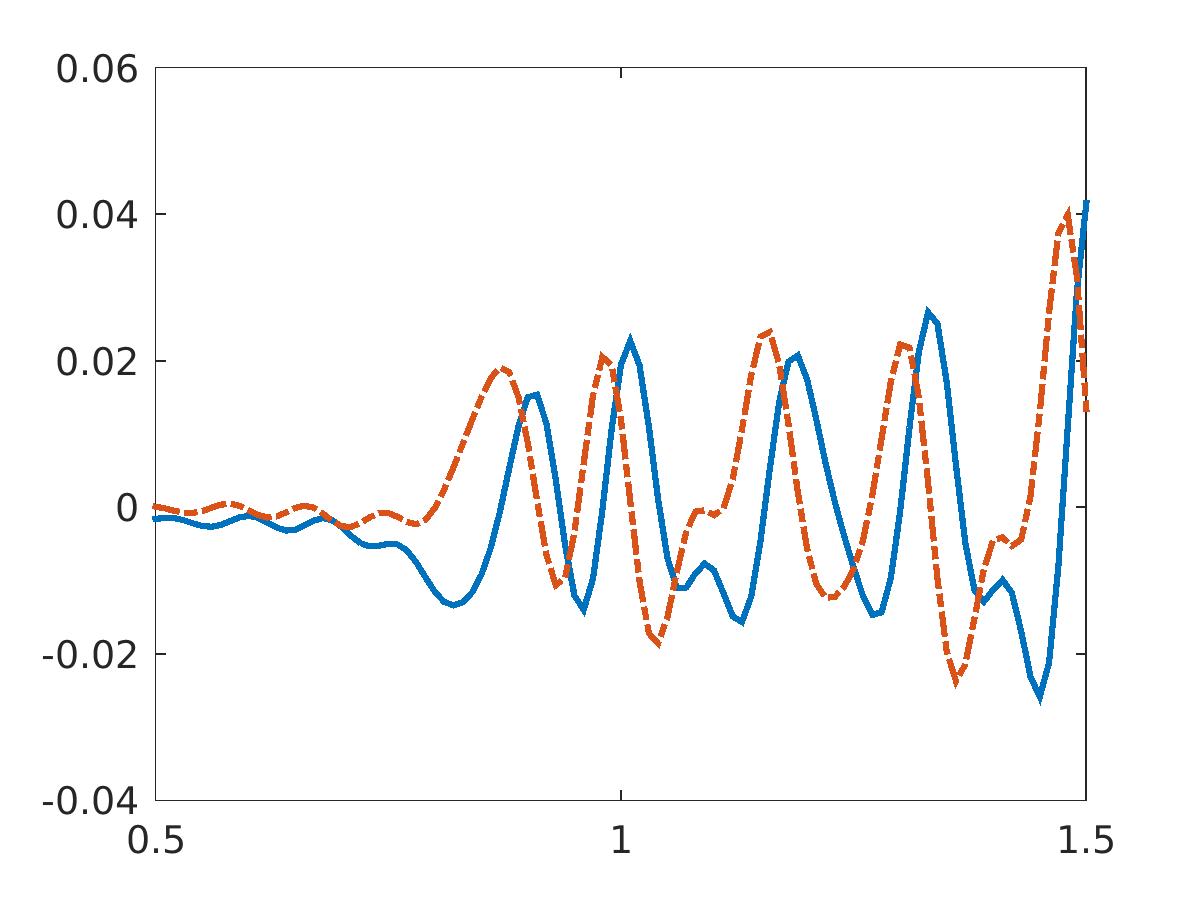} 
} \hfill%
\subfloat[\labeldd]{
		\includegraphics[width=\idth, height = \eight ]{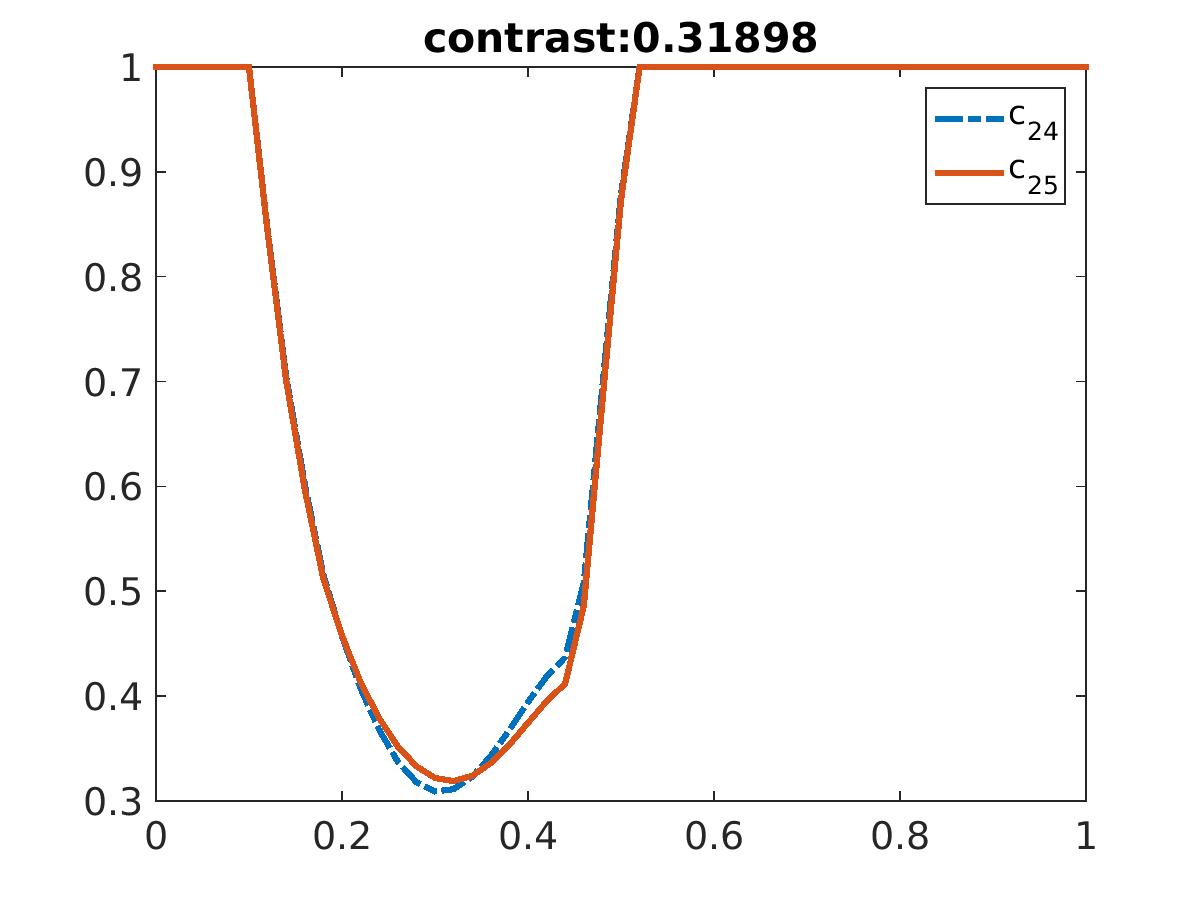}
}
\par
\hfill 
\subfloat[\labelaa]{
		\includegraphics[width=\idth, height = \eight ]{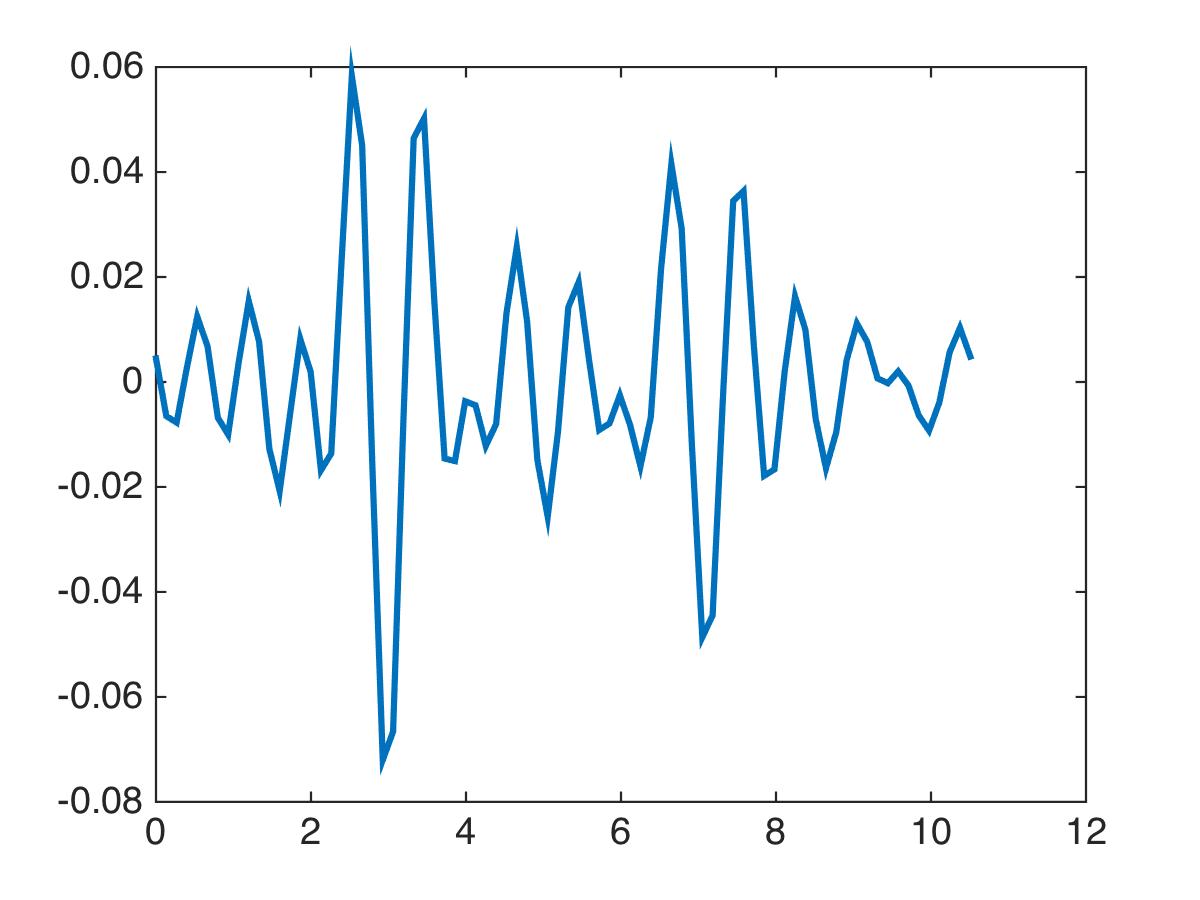} 
} \hfill%
\subfloat[\labelbb]{
		\includegraphics[width=\idth, height = \eight ]{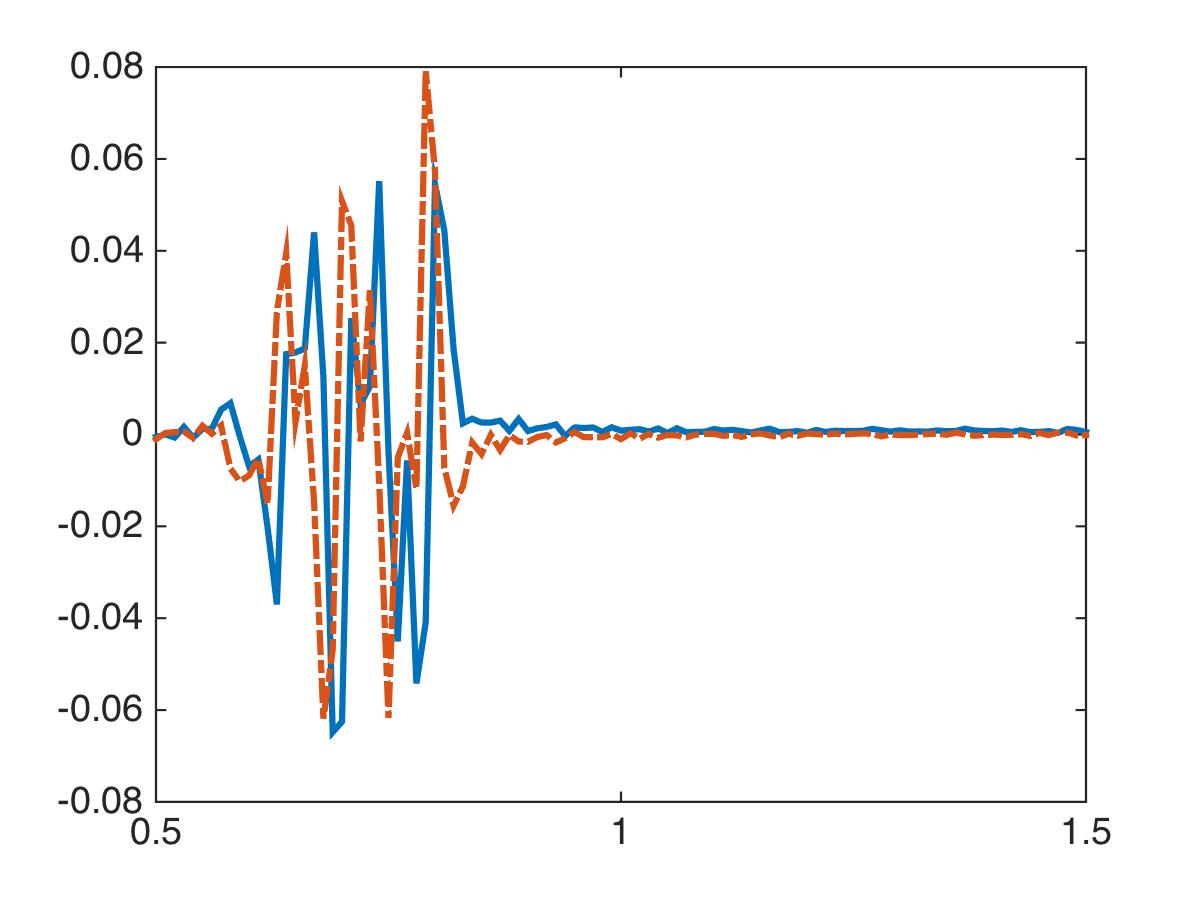} 
} \hfill%
\subfloat[\labelcc]{
		\includegraphics[width=\idth, height = \eight ]{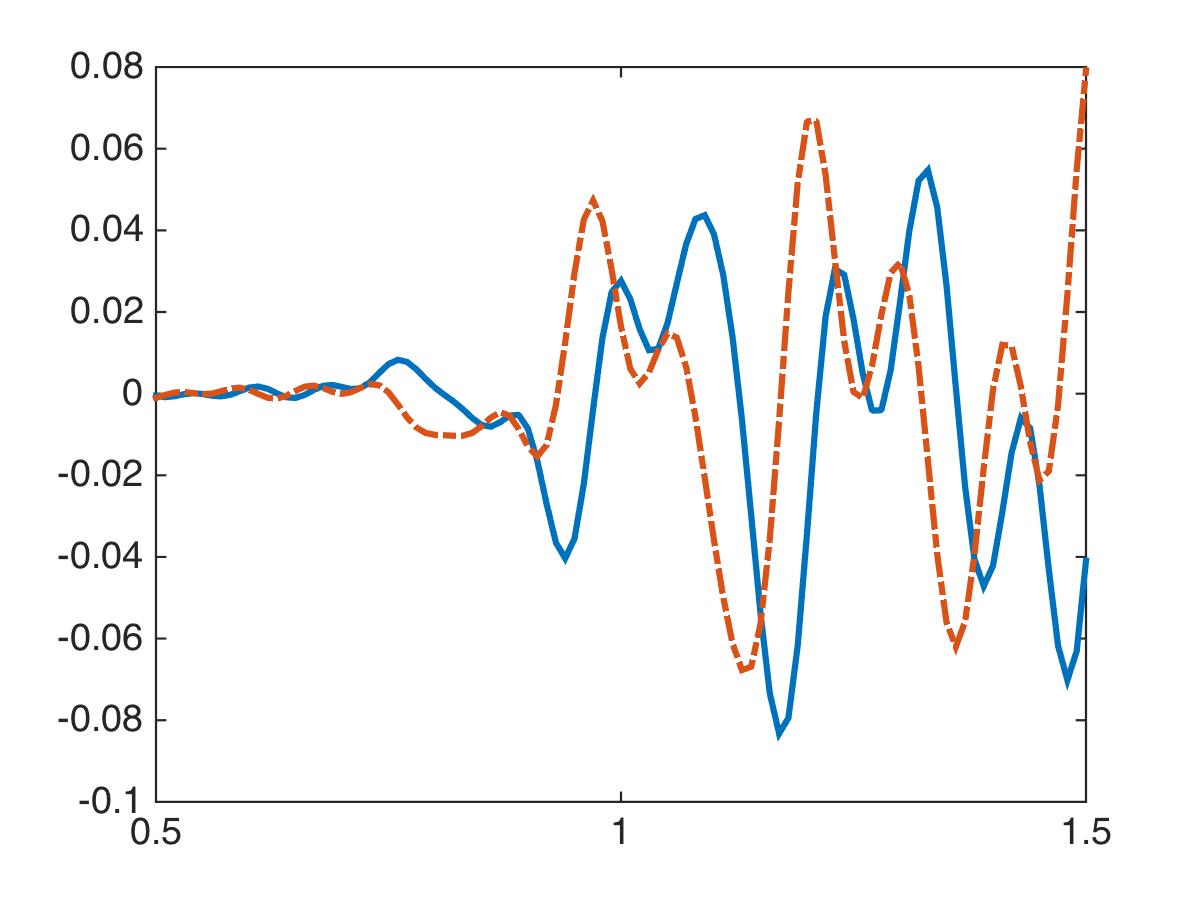} 
} \hfill%
\subfloat[\labeldd]{
		\includegraphics[width=\idth, height = \eight ]{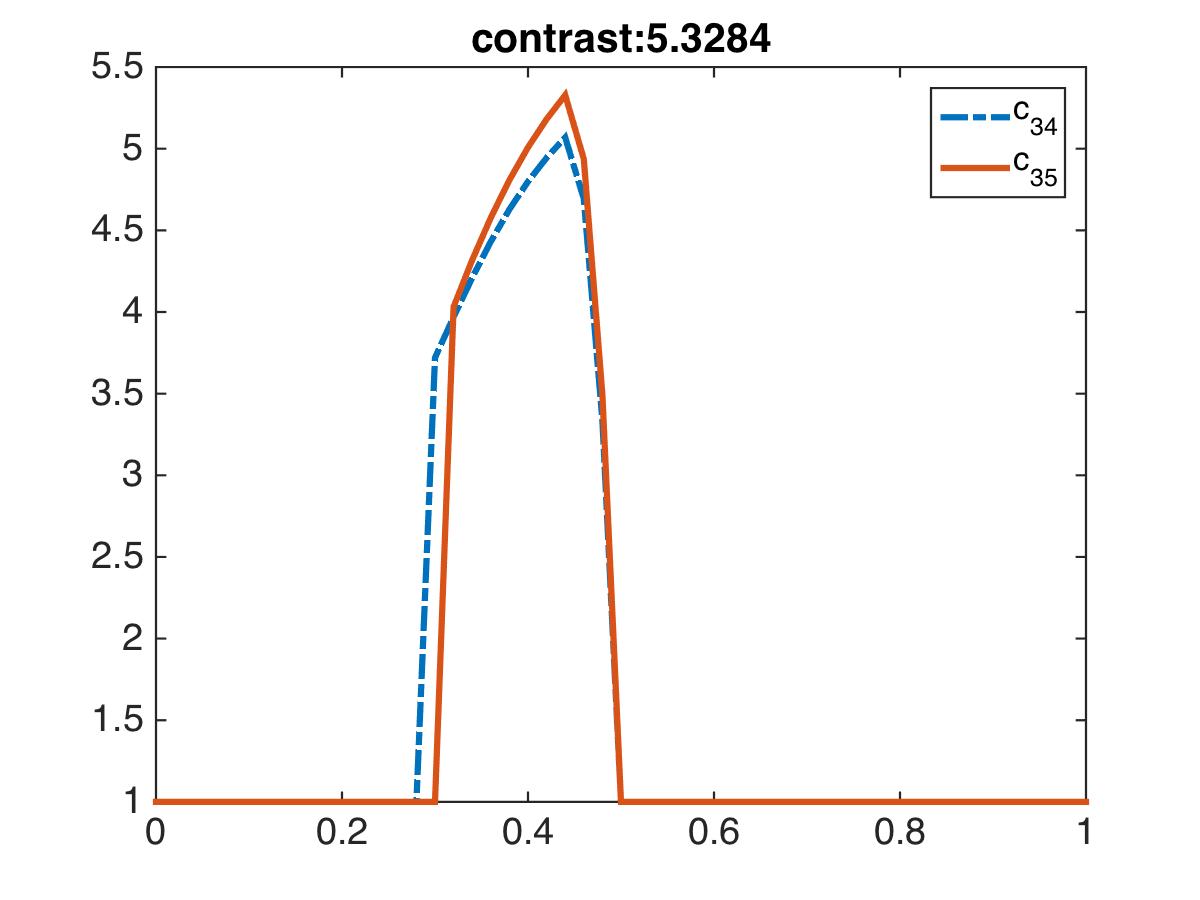}
}
\par
\hfill%
\subfloat[\labelaa]{
		\includegraphics[width=\idth, height = \eight ]{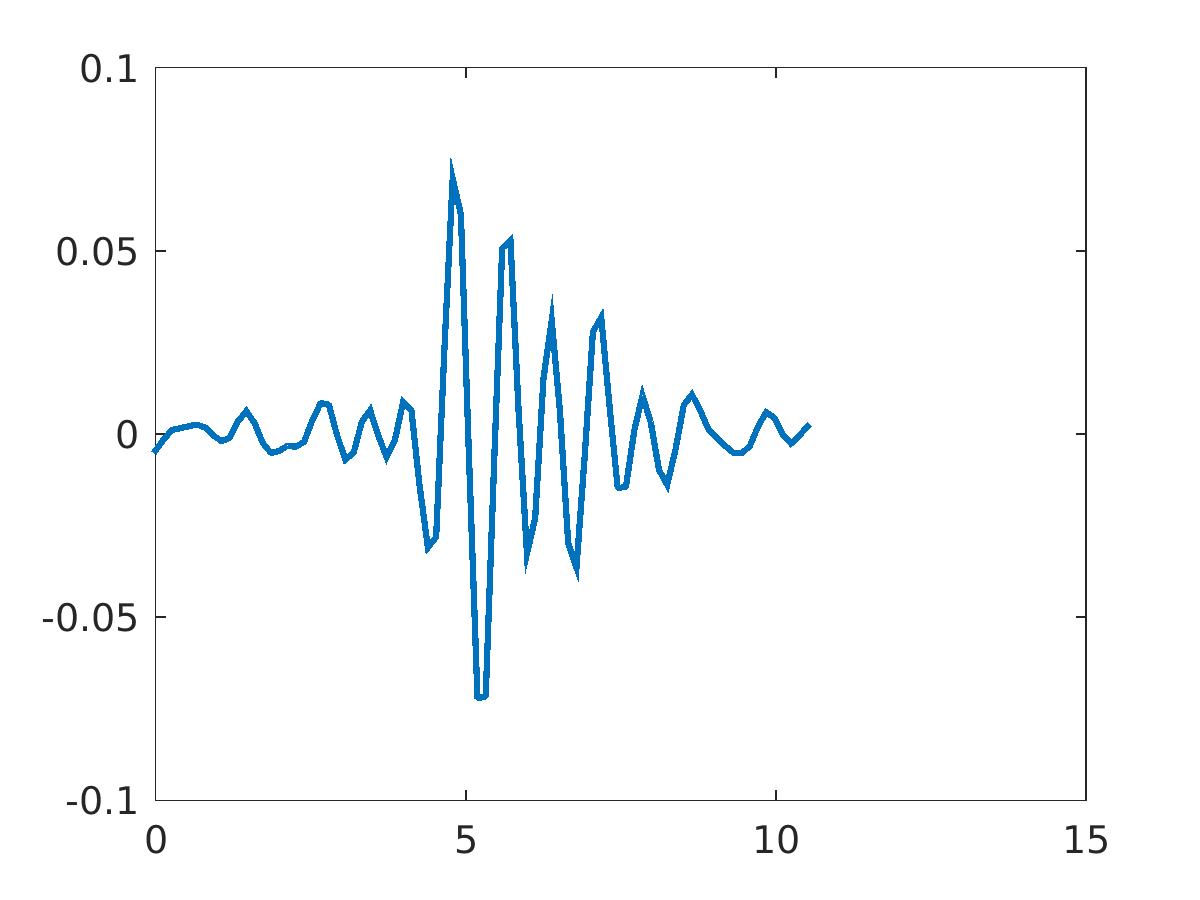} 
} \hfill%
\subfloat[\labelbb]{
		\includegraphics[width=\idth, height = \eight ]{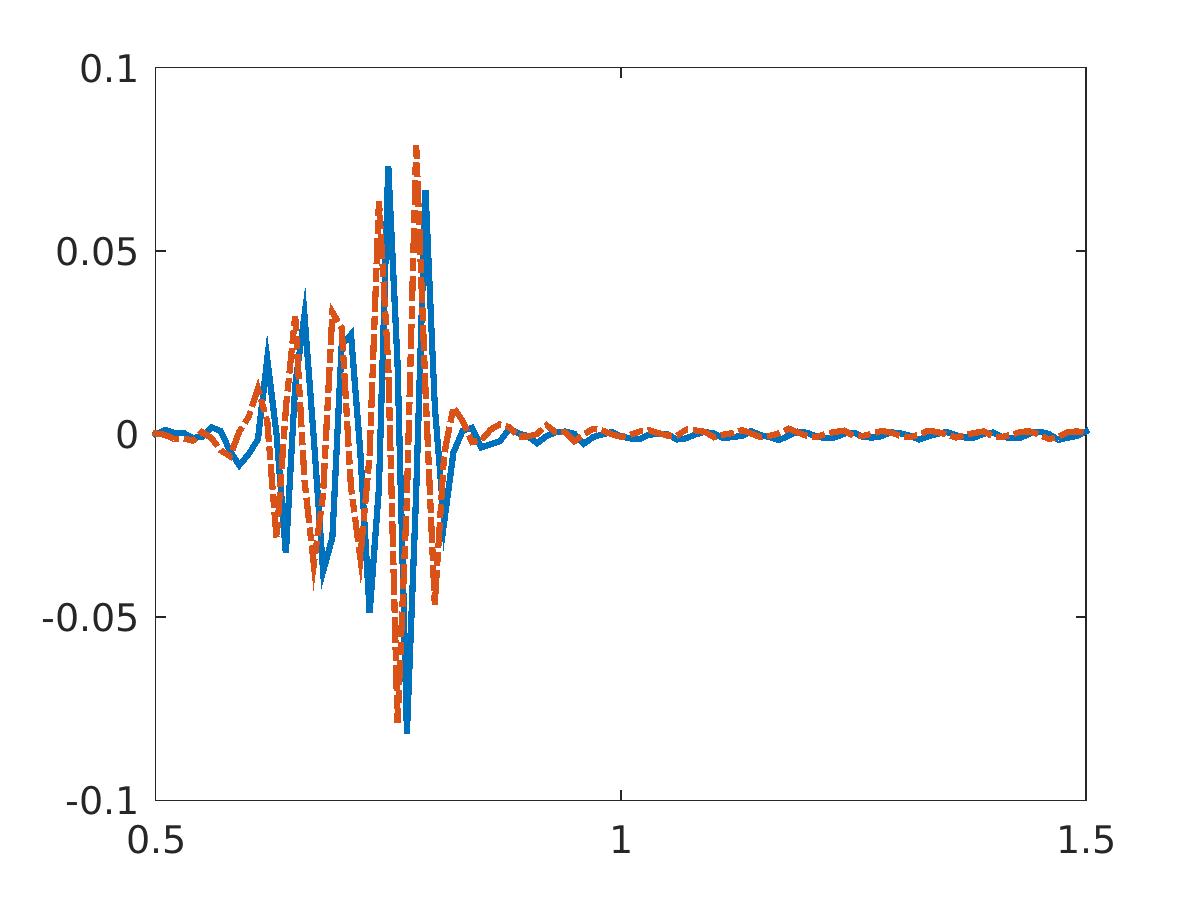} 
} \hfill%
\subfloat[\labelcc]{
		\includegraphics[width=\idth, height = \eight ]{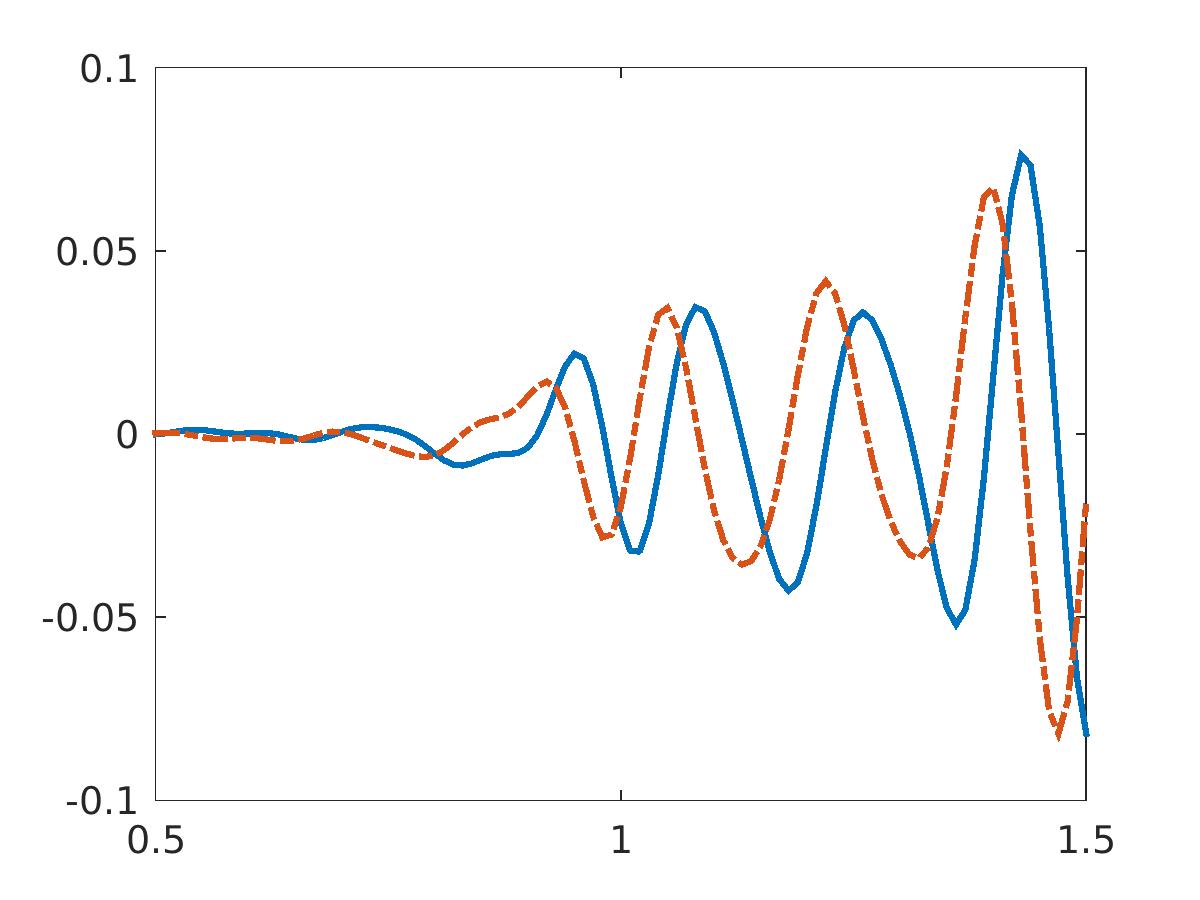} 
} \hfill%
\subfloat[\labeldd]{
		\includegraphics[width=\idth, height = \eight ]{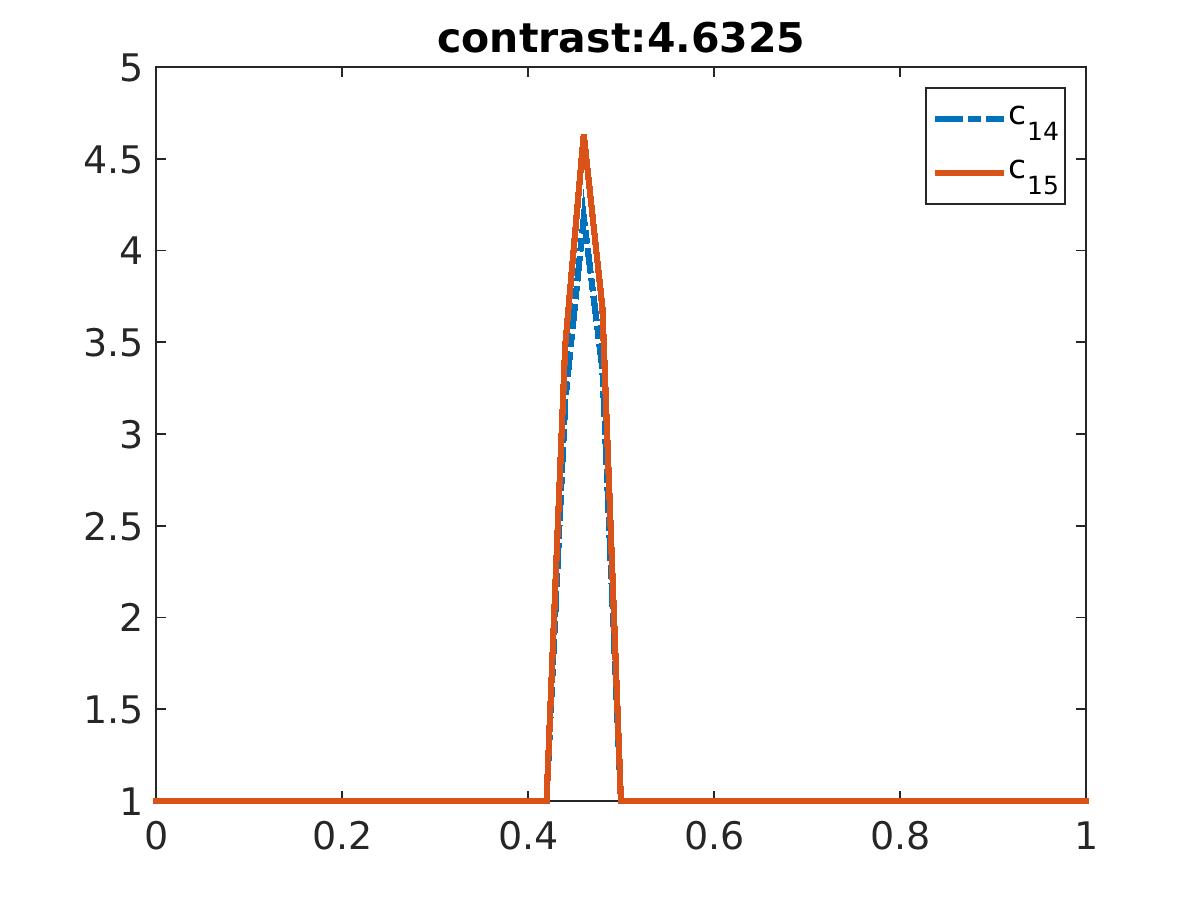}
}
\caption{\textit{The case when plastic (row 1), metal cylinder (row 2), and
metal box(row 3) are buried under the ground. Solid lines on b,c,f,g,j,k are
real parts and dotted lines are imaginary parts.}}
\label{Ground}
\end{figure}

We consider two cases: targets in air (Figure \ref{Air}) and targets buried
about a few centimeters under the ground (Figure \ref{Ground}). We had
experimental data for total of five (5) targets. The reconstructed
dielectric constants of these targets are summarized in Table 1. In this
table, computed $c_{\text{comp}}=\widetilde{R}\cdot c_{\text{bckgr}}$. In
tables of dielectric constants of materials, their values are usually given
within certain intervals \cite{Tabl}. Now about the intervals of the true $%
c:=c_{\text{true}}$ in the $6^{\text{th}}$ column of Table 1. In the cases
when targets were a wood stake and a plastic cylinder, we have taken those
intervals from a published table of dielectric constants \cite{Tabl}. The
interval of the true $c$ for the case when the target was bush, was taken
from \cite{Ch}. As to the metal targets, it was established in \cite{KSNF1}
that they can be considered as such targets whose dielectric constants
belong to the interval $[10,30].$

\begin{center}
\textbf{Table 1:}{\ Computed dielectric constants of five targets} 
\begin{tabular}{|c|c|c|c|c|c|}
\hline
Target & $c_{\mathrm{bckgr}}$ & Reconstructed $\widetilde{R}$ & $c_{\mathrm{%
bckgr}}$ & $c_{\text{comp}}$ & True $c_{\text{true}}$ \\ \hline
Bush & 1 & 6.5 & 1 & 6.5 & $[3,20]$ \\ 
Wood stake & 1 & 3.3 & 1 & 3.3 & $[2,6]$ \\ 
Metal box & 4 & 4.6 & $[3,5]$ & $[13.8,23]$ & $[10,30]$ \\ 
Metal cylinder & 4 & 5.3 & $[3,5]$ & $[15.9,26.5]$ & $[10,30]$ \\ 
Plastic cylinder & 4 & 0.3 & $[3,5]$ & $[0.9,1.5]$ & $\left[ 1.1,3.2\right] $
\\ \hline
\end{tabular}%
\\[0pt]
\end{center}

\section{Summary}

\label{sec:7}

In this paper, we have developed a frequency domain analog of the 1-d
globally convergent method of \cite{Karch,KSNF1,KSNF2}. We have tested this
analog on both computationally simulated and time resolved experimental
data. The experimental data are the same as ones used in \cite%
{Karch,KSNF1,KSNF2}. We have modeled the process of electromagnetic waves
propagation by the 1-d wave-like PDE. The reason why we have not used a 3-d
model, as in, e.g. earlier works of this group on experimental data \cite%
{BK1,TBKF1,TBKF2}, is that we had only one time resolved experimental curve
for each of our five targets.

Our numerical method has the global convergence property. In other words, we
have proven a theorem (Theorem 5.1), which claims that we obtain some points
in a sufficiently small neighborhood of the exact solution without any
advanced knowledge of this neighborhood. Our technique heavily relies on the
Quasi Reversibility Method (QRM). The proof of the convergence of the QRM is
based on a Carleman estimate. A significant modification of our technique,
as compared with \cite{Karch,KSNF1,KSNF2}, is due to two factors. First, we
use the Fourier transform of time resolved data instead of the Laplace
transform in \cite{Karch,KSNF1,KSNF2}. Second, when updating tail functions
via (\ref{4.17}), we solve the problem (\ref{4.3}), (\ref{4.4}) using the
QRM. On the other hand, in \cite{Karch,KSNF1,KSNF2} tail functions were
updated via solving the \textquotedblleft Laplace transform analog" of the
problem (\ref{2.4}), (\ref{2.6}) as a regular forward problem

Since the dielectric constants of targets were not measured in experiments,
the maximum what we can do to evaluate our results is to compare them with
published data in, e.g. \cite{Tabl}. Results of Table 1 are close to those
obtained in \cite{Karch,KSNF1,KSNF2}. One can see in Table 1 that our
reconstructed values of dielectric constants are well within published
limits. We consider the latter as a good result. This is achieved regardless
on a significant discrepancy between computationally simulated and
experimental data, regardless on a quite approximate nature of our
mathematical model and regardless on the presence of clutter at the data
collection site. That discrepancy is still very large even after the data
pre-processing, as it is evident from comparison of Figures 1b,e with
Figures 2c,g and Figures 3c,g,k. Besides, the source position $x_{0}$ was
unknown but rather prescribed by ourselves as $x_{0}=-1.$ Thus, our results
indicate a high degree of stability of our method.

\end{document}